# GLOBAL EXISTENCE FOR A QUASILINEAR WAVE EQUATION OUTSIDE OF STAR-SHAPED DOMAINS

MARKUS KEEL, HART F. SMITH, AND CHRISTOPHER D. SOGGE

*In memory of Tom Wolff*

1. **Introduction.**

The purpose of this paper is to establish global existence of small-amplitude solutions for certain quasilinear Dirichlet-wave equations outside of smooth, compact star-shaped obstacles $\mathcal{K} \subset \mathbb{R}^3$. Precisely, we shall consider smooth quasilinear systems of the form

$$\begin{cases} \partial_t^2 u - \Delta u = F(u, du, d^2u), \quad (t,x) \in \mathbb{R}_+ \times \mathbb{R}^3 \backslash \mathcal{K} \\ u(t, \cdot)|_{\partial \mathcal{K}} = 0 \\ u(0, \cdot) = f, \quad \partial_t u(0, \cdot) = g, \end{cases} \quad (1.1)$$

which satisfy the so-called null condition [15]. The global existence for such equations in the absence of obstacles was established by Christodoulou [2] and Klainerman [11] using different techniques. We begin by describing our assumptions in more detail.

We let $u$ denote a $N$-tuple of functions, $u = (u^1, u^2, \ldots, u^N)$. We assume that $\mathcal{K}$ is smooth and strictly star-shaped with respect to the origin. By this, we understand that in polar coordinates $x = r\omega$, $(r, \omega) \in [0, \infty) \times S^2$, we can write

$$\mathcal{K} = \{x = r\omega : \phi(\omega) - r \geq 0\}, \quad (1.2)$$

where $\phi$ is a smooth *positive* function on $S^2$. Thus,

$$0 \in \mathcal{K}, \text{ but } 0 \notin \partial\mathcal{K} = \{x : r = \phi(\omega)\}.$$

By quasilinearity, we mean that $F(u, du, d^2u)$ is linear in the second derivatives of $u$. We shall also assume that the highest order nonlinear terms are diagonal, by which we mean that, if we denote $\partial_0 = \partial_t$, then

$$F^I(u, du, d^2u) = G^I(u, du) + \sum_{0 \leq j,k \leq 3} \gamma^{I,jk}(u, du)\, \partial_j \partial_k u^I, \quad 1 \leq I \leq N. \quad (1.3)$$

A key assumption is that the nonlinear terms satisfy the null condition. Recall that even in the obstacle-free case there can be blowup in finite time for arbitrarily small data if this condition is not satisfied (see John [9]).

The first part of the null condition is that the nonlinear terms are free of linear terms,

$$F(0, 0, 0) = 0, \text{ and } F'(0, 0, 0) = 0. \quad (1.4)$$

The authors were supported in part by the NSF.

Additionally, we assume that the quadratic terms do not depend on $u$, which means that we can write

$$F(u, du, d^2u) = Q(du, d^2u) + R(u, du, d^2u), \tag{1.5}$$

where $Q$ is a quadratic form, and where the remainder term $R$ vanishes to third order at $(u, du, d^2u) = 0$; that is,

$$R(p, q, r) = O\big((p^2 + q^2)r\big) + O\big((|p| + |q|)^3\big). \tag{1.6}$$

The null condition concerns the quadratic term $Q$. To describe it, we split $Q$ into its semilinear and quasilinear parts:

$$Q(du, d^2u) = s(du, du) + k(du, d^2u).$$

Then in terms of the $N$ components of $u$ we can rewrite these terms as

$$s^I(du, du) = \sum_{1 \leq J, K \leq N} \sum_{0 \leq j, k \leq 3} s^{I,j,k}_{J,K} \, \partial_j u^J \, \partial_k u^K,$$

and

$$k^I(du, d^2u) = \sum_{J=1}^{N} \sum_{0 \leq i, j, k \leq 3} k^{I,i,j,k}_J \, \partial_i u^J \, \partial_j \partial_k u^I,$$

where the $s^{I,j,k}_{J,K}$ and $k^{I,i,j,k}_J$ are constants. The null condition can then be stated succinctly as requiring that, if $1 \leq I, J, K \leq N$,

$$\sum_{0 \leq j, k \leq 3} s^{I,j,k}_{J,K} \xi_j \xi_k = 0, \quad \text{if} \quad \xi_0^2 = \xi_1^2 + \xi_2^2 + \xi_3^2,$$

and, if $1 \leq I, J \leq N$,

$$\sum_{0 \leq i, j, k \leq 3} k^{I,i,j,k}_J \xi_i \xi_j \xi_k = 0 \quad \text{if} \quad \xi_0^2 = \xi_1^2 + \xi_2^2 + \xi_3^2.$$

For further discussion, we refer the reader to Christodoulou [2], p. 277–278.

As was shown in [2] and [11], this condition forces the semilinear terms $s^I(du, du)$ to be linear combinations of the basic null forms

$$q_0(du^J, du^K) = \partial_0 u^J \, \partial_0 u^K - \sum_{j=1}^{3} \partial_j u^J \, \partial_j u^K, \tag{1.7}$$

and

$$q_{ij}(du^J, du^K) = \partial_i u^J \, \partial_j u^K - \partial_j u^J \, \partial_i u^K, \quad 0 \leq i, j \leq 3. \tag{1.8}$$

The quasilinear term $k^I(du, d^2u)$ in turn must be a linear combination of terms of the form

$$q(du^J, d\partial_j u^I), \quad 1 \leq J \leq N, \quad 0 \leq j \leq 3, \tag{1.9}$$

where $q$ is a basic null form as in (1.7) and (1.8), along with terms of the form

$$\partial_j u^J \big(\partial_t^2 u^I - \Delta u^I\big), \quad 1 \leq J \leq N, \quad 0 \leq j \leq 3. \tag{1.10}$$

In addition to the null condition, we must assume that the Cauchy data $f, g$ satisfy certain compatibility conditions at the boundary. We leave the statement of these conditions to Definition 9.2.

As in Christodoulou's [2] results for the non-obstacle case, we shall not need to assume that the data has compact support. Instead, we make the assumption that $f$ and $g$ belong to certain weighted Sobolev spaces. To state our assumptions precisely, we recall the weighted Sobolev spaces used by Christodoulou [2], which are given by the norm

$$\|f\|_{H^{m,j}(\mathbb{R}^3)} = \sum_{|\alpha|\leq m} \left(\int_{\mathbb{R}^3} \left(1+|x|^2\right)^{|\alpha|+j}|\partial_x^\alpha f(x)|^2\, dx\right)^{1/2}.$$

The associated weighted Dirichlet-Sobolev spaces for $m = 1, 2\ldots$ are defined by

$$H_D^{m,j}(\mathbb{R}^3\setminus\mathcal{K}) = \{f \in H^{m,j}(\mathbb{R}^3\setminus\mathcal{K}) : f|_{\partial\mathcal{K}} = 0\}, \tag{1.11}$$

where $H^{m,j}(\mathbb{R}^3\setminus\mathcal{K})$ is the space of restrictions of elements of $H^{m,j}(\mathbb{R}^3)$. Hence,

$$\|f\|_{H_D^{m,j}}^2 = \sum_{|\alpha|\leq m} \int_{\mathbb{R}^3\setminus\mathcal{K}} \left(1+|x|^2\right)^{|\alpha|+j}|\partial_x^\alpha f(x)|^2\, dx, \tag{1.12}$$

gives the natural norm on $H_D^{m,j}(\mathbb{R}^3\setminus\mathcal{K})$. We can now state our main result.

**Theorem 1.1.** *Assume that $\mathcal{K}$ and $F(u, du, d^2u)$ are as above. Assume further that $(f, g) \in C^\infty(\mathbb{R}^3\setminus\mathcal{K})$ satisfies the compatibility conditions to infinite order (see Definition 9.2). Then there exists $\varepsilon_0 > 0$, such that if*

$$\|f\|_{H_D^{9,8}(\mathbb{R}^3\setminus\mathcal{K})} + \|g\|_{H_D^{8,9}(\mathbb{R}^3\setminus\mathcal{K})} < \varepsilon_0\,, \tag{1.13}$$

*then there is a unique solution $u \in C^\infty(\mathbb{R}_+ \times \mathbb{R}^3\setminus\mathcal{K})$ of (1.1). Furthermore, for all $\sigma > 0$, there exists $C_\sigma < \infty$, such that*

$$|u(t,x)| \leq C_\sigma \left(1+t\right)^{-1}\left(1+|t-|x||\right)^{-1+\sigma}. \tag{1.14}$$

We will actually establish existence of limited regularity solutions $u$ for data $f \in H_D^{9,8}$ and $g \in H_D^{8,9}$ satisfying compatibility conditions of order 8; see Theorem 7.1. The fact that $u$ is smooth if $f$ and $g$ are smooth and satisfy compatibility conditions of infinite order will follow by the local existence theorems of section 9.

It should be possible to relax the regularity assumptions in the smallness condition (1.13). In particular, our techniques should just require that $\|f\|_{H_D^{4,3}} + \|g\|_{H_D^{3,4}}$ be small, which would be the analog of Christodoulou's assumption in [2]. Additionally, the result should hold with $\sigma = 0$.

The authors [10] were able to show that if $\mathcal{K}$ is strictly convex then one has global existence for the semilinear case for data $f \in H_D^{2,1}$ and $g \in H_D^{1,2}$. The work was based on a variant of Christodoulou's method which involved weighted estimates, where, as in the present work, the weights on the derivatives compensate for the degeneracy of the image of $\mathbb{R}_+ \times \partial\mathcal{K}$ as $t \to +\infty$ under Penrose's conformal compactification of Minkowski space. The proof depended on results of the last two authors [24] which extended estimates of Klainerman and Machedon [15] to the setting of strictly convex obstacles. These results are not known in the setting of general star-shaped obstacles.

The special case of Theorem 1.1 in which one assumes spherical symmetry for $u$ and $\mathcal{K}$ was obtained by Godin in [4]. His proof involved an adaptation of Christodoulou's [2] method to this setting. If one drops the assumption of spherical symmetry, it does not appear that the arguments in [4] will apply in a straightforward way.

Also, results similar to those in Theorem 1.1 were announced in Datti [3], but there appears to be a gap in the argument which has not been repaired.

Previous work in higher dimensions applied Lorentz vector field techniques to the exterior problem. For general nonlinearities quadratic in $du$, global smooth solutions were shown by Shibata and Tsutsumi [20], [21] to exist for dimensions $n \geq 6$. In Hayashi [5], global existence of smooth solutions in the exterior of a sphere for $n \geq 4$ is shown for a restricted class of quadratic nonlinearities.

Let us give an overview of our proof of Theorem 1.1. First of all, as in Christodoulou [2], we shall use the so-called conformal method (see also [6]). Thus, we shall apply Penrose's conformal compactification of Minkowski space. Recall this is a map $\mathcal{P} : \mathbb{R} \times \mathbb{R}^3 \to (-\pi, \pi) \times S^3$, where the image is the so-called Einstein diamond

$$\overline{\mathbb{E}}^4 = \{(T, X) \in (-\pi, \pi) \times S^3 : |T| + R < \pi\} \subset \mathbb{E}^4 = (-\pi, \pi) \times S^3.$$

Here $R$ denotes the distance on $S^3$ from the north pole

$$\mathbf{1} = (1, 0, 0, 0)$$

measured in the standard metric. The Penrose map preserves the angular variable, while if $r$ is the radial variable and $t$ the time variable in Minkowski space then under $\mathcal{P}$ the corresponding variables in $\overline{\mathbb{E}}^4$ are related as follows

$$\begin{aligned} R &= \arctan(t+r) - \arctan(t-r), \\ T &= \arctan(t+r) + \arctan(t-r). \end{aligned} \tag{1.15}$$

Under this map the pushforward of the Minkowski metric $dt^2 - dx^2$ is the Lorentz metric $\tilde{g}$ in $\overline{\mathbb{E}}^4$ given by

$$dT^2 - g = \Omega^2 \tilde{g}, \tag{1.16}$$

where $dT^2 - g$ is the standard Lorentz metric on $\mathbb{R} \times S^3$, and where the conformal factor $\Omega$ is given by the formula

$$\Omega = \cos T + \cos R = \frac{2}{(1 + (t+r)^2)^{1/2}(1 + (t-r)^2)^{1/2}}, \tag{1.17}$$

with $(T, R)$ and $(t, r)$ being identified as above.

Continuing, let

$$\Box_g = \partial_T^2 - \Delta_g$$

be the D'Alembertian coming from the standard Laplace-Beltrami operator $\Delta_g$ on $S^3$. If we change our earlier notation a bit and let $\tilde{\Box}$ denote the D'Alembertian on $\mathbb{R}^{1+3}$ or $\overline{\mathbb{E}}^4$, depending on the context, that arises from the standard Lorentz metric $dt^2 - dx^2$, then a key fact for us is the way that the two D'Alembertians are related in $\overline{\mathbb{E}}^4$:

$$\Box_g + 1 = \Omega^{-3} \tilde{\Box} \Omega, \tag{1.18}$$

with the additive constant 1 arising because of the non-zero scalar curvature of $g$. Equivalently,

$$\tilde{\Box} \tilde{u} = F \iff (\Box_g + 1)v = G \text{ with } \tilde{u} = \Omega v \text{ and } G = \Omega^{-3} F. \tag{1.19}$$

On account of this, if
$$\mathcal{K}_* = \mathcal{P}\big([0, +\infty) \times \mathcal{K}\big) \tag{1.20}$$
is the pushforward of our obstacle in Minkowski space, then the task of showing that we can find small-amplitude solutions of (1.1) is equivalent to showing that we can find small-amplitude solutions of
$$\begin{cases} (\Box_g + 1)v = G(v, dv, d^2v), & (T,X) \in \overline{\mathbb{E}}^4_+ \setminus \mathcal{K}_* \\ v(T,X) = 0, & (T,X) \in \partial \mathcal{K}_* \end{cases}$$
with $\overline{\mathbb{E}}^4_+ = \{(T,X) \in \overline{\mathbb{E}}^4 : 0 \le T < \pi\}$, and $u, v, F$, and $G$ related as above.

Christodoulou [2] showed that the transformed nonlinear term $G$ extends to a nonlinear term with $C^\infty$ coefficients on the cylinder $\mathbb{R} \times S^3$ if and only if the null condition is satisfied. Indeed, the transformed quadratic terms coming from $Q$ in (1.5) extend analytically to the cylinder if and only if the null condition is verified (see [2], p. 277-278), while the transformed remainder term coming from $R$ trivially extends smoothly because of (1.6) and (1.19). Because of this, as was argued in [2], the assertion that there are small-amplitude global solutions for $(\partial_t^2 - \Delta)u = F(u, du, d^2u)$ verifying the null condition in the boundaryless Minkowski case just follows from a routine local existence theorem for $\mathbb{R}_+ \times S^3$.

This simple approach breaks down for obstacle problems due to the fact that the transformed obstacle $\mathcal{K}_*$ given by (1.20) is a time-dependent obstacle which collapses to a point as $T \to \pi$. Indeed, it follows from (1.15) that there must be a uniform constant $1 < C < \infty$ so that for $0 \le T < \pi$
$$C^{-1}(\pi - T)^2 \le \mathrm{dist}(X, \mathbf{1}) \le C(\pi - T)^2, \quad \text{if } (T,X) \in \partial \mathcal{K}_*, \tag{1.21}$$
with $\mathbf{1}$ as above being the north pole on $S^3$. Thus, if we let
$$P_0 = (\pi, \mathbf{1}), \tag{1.22}$$
it follows that $\mathcal{K}_*$ collapses to $P_0$ as $T \to \pi$.

Following the approach in our earlier work [10], we shall surmount this difficulty by modifying the usual existence arguments for the non-obstacle case. In our approach, we shall need to obtain and apply estimates that involve weighted derivatives because of the quadratic degeneracy of $\mathcal{K}_*$ at $P_0$.

To state our main estimates we need to introduce some more notation. We let $X_j$, $j = 0, 1, 2, 3$ be the coordinate functions on $\mathbb{R}^4$, and then let
$$\frac{\partial}{\partial T}, \quad X_j \frac{\partial}{\partial X_k} - X_k \frac{\partial}{\partial X_j}, \quad 0 \le j < k \le 3 \tag{1.23}$$
be the spanning set of vector fields on $\overline{\mathbb{E}}^4$, where we identify $S^3 = \{X \in \mathbb{R}^4 : |X| = 1\}$.

We arrange these vector fields as $\Gamma = \{\Gamma_0, \ldots, \Gamma_6\}$. Our main estimates will involve the weighted derivatives
$$Z^\alpha = \big[(\pi - T)^2 \Gamma\big]^\alpha = \big((\pi - T)^2 \Gamma_0\big)^{\alpha_0} \cdots \big((\pi - T)^2 \Gamma_6\big)^{\alpha_6}. \tag{1.24}$$
These turn out to be natural to use due to the fact that, near $\mathcal{K}_*$, $Z_j$ pulls back via $\mathcal{P}$ to a vector field in Minkowski space that essentially has unit length. As a side remark, this is not the case near the set where $T + R = \pi$, and this accounts for the importance of the null condition in three spatial dimensions.

To show that we can solve the transformed equation, and hence the original (1.1), we need certain $L^2$ estimates and pointwise estimates involving $Z^\alpha$. Special cases of the $L^2$ estimates state that if $v$ solves the Dirichlet-wave equation $(\Box_g + 1)v = G$, $v|_{\partial \mathcal{K}_*} = 0$, then under appropriate conditions on the data and forcing terms we have

$$\sum_{|\alpha| \leq k} \|Z^\alpha v'(T, \cdot)\|_2 \leq C_k \int_0^T \sum_{|\alpha| \leq k} \|Z^\alpha G(S, \cdot)\|_2 \, dS$$
$$+ C_k \sup_{0 < S < T} (\pi - S)^2 \sum_{|\alpha| \leq k-1} \|Z^\alpha G(S, \cdot)\|_2 + C_k \sum_{|\alpha| \leq k} \|Z^\alpha v'(0, \cdot)\|_2, \quad 0 < T < \pi, \tag{1.25}$$

where, for a given $T$, the norms are taken over $\{(X : (T, X) \in \mathbb{E}_+^4 \setminus \mathcal{K}_*\}$. The key step in the proof of this will be to show that the bounds hold when $k = 0$:

$$\|v'(T, \cdot)\|_2 \leq C\|v'(0, \cdot)\|_2 + C \int_0^T \|G(S, \cdot)\|_2 \, dS. \tag{1.26}$$

Here, as throughout this paper, $v'$ denotes the unweighted 4-gradient of $v$, or equivalently $v'$ denotes the collection $\{\Gamma_j v, \, 0 \leq j \leq 6\}$.

To prove (1.26) we shall adapt Morawetz's [18] proof of a related estimate in Minkowski space outside star-shaped obstacles. The proof of (1.26) is based on the fact that, when one applies standard arguments involving the energy-momentum tensor, the boundary integrals that arise have integrands with a favorable sign. Because of this, we can also obtain energy estimates for appropriate small variable coefficient perturbations of $\Box_g$. The fact that the analog of (1.25) remains valid in this setting is necessary to handle the nonlinear perturbations of the metric in (1.1). If $\mathcal{X} = \mathcal{P}_*(\partial/\partial t)$ is the pushforward of the Minkowski time derivative, then a key step in seeing that (1.25) follows from (1.26) is that a variant of (1.26) holds when $v$ is replaced by $\mathcal{X}v$, since $\mathcal{X}v$ also satisfies the Dirichlet boundary condition.

A special case of our pointwise estimate states that if $(\Box_g + 1)v = G$, $v|_{\mathcal{K}_*} = 0$, then under appropriate assumptions on the data and forcing term, if $p > 1$ is fixed then for $0 < T < \pi$ we have uniform bounds

$$|v(T, X)| \leq C \sup_{0 \leq S \leq T} \sum_{|\alpha| \leq 1} \left( \|Z^\alpha G(S, \cdot)\|_2 + (\pi - S)^{-2} \|Z^\alpha G(S, \cdot)\|_p \right)$$
$$+ C \sum_{|\alpha| \leq 1} \|Z^\alpha v'(0, \cdot)\|_2. \tag{1.27}$$

We shall also obtain analogous estimates for $Z^\alpha v$. These estimates imply that the solution of the transformed version of (1.1) to $\overline{\mathbb{E}}_+^4$ satisfies

$$Z^\alpha v(T, X) = O\big((\pi - T)^{-\sigma}\big) \quad \text{for } \sigma > 0. \tag{1.28}$$

For technical reasons, we do not obtain uniform bounds $\sigma = 0$ due to the fact that (1.27) only holds for Lebesgue exponents $p > 1$.

In our earlier work [10] on the semilinear case, we showed only that the solution of the transformed version of (1.1) satisfies (1.28) with $\sigma = 1$. As we shall see, the fact that we can now obtain bounds which blow up like $(\pi - T)^{-\sigma}$ for some $\sigma < 1$ plays a crucial role in our analysis. This is because the iterations we shall use in this paper would involve logarithmic terms if $\sigma = 1$, and hence be useless.

We also remark that the proof of (1.28) is modeled after the recent proof by the last two authors [25] of global Strichartz estimates outside of convex obstacles.

## 2. The conformal transformation and the transformed equation.

In this section we provide further details about the conformal method. In particular, we recall formulas which relate derivatives in Minkowski space to derivatives in the Einstein diamond. We also go over estimates for the nonlinear terms of the pushforward via $\mathcal{P}$ of equations such as (1.1) which satisfy the null condition. As we stressed in the introduction, it is important for our analysis that the nonlinear terms are small near the "tip" $P_0$ of the Einstein diamond defined by (1.22). Finally, we show how the weighted Dirichlet-Sobolev spaces $H_D^{m,j}$ in (1.12) are related to the usual Dirichlet-Sobolev spaces on the 3-sphere minus an obstacle.

We start by reviewing the way that derivatives transform under Penrose's conformal compactification of Minkowski space. For this it is convenient to use stereographic projection coordinates on $S^3$. We note that the south pole stereographic projection coordinates $U$ arise as the restriction of $\mathcal{P}^{-1}$ to the slice $T = 0$:

$$U = \mathcal{P}_0^{-1}(X) = \frac{\sin R}{1 + \cos R}\,\omega = \tan(\tfrac{R}{2})\,\omega\,. \tag{2.1}$$

The coordinates $V$ of the stereographic north pole projection are obtained by applying the Kelvin transform to the south pole stereographic coordinate,

$$V_j = |U|^{-2}\,U_j\,. \tag{2.2}$$

To compute the pushforwards of vector fields on $\overline{\mathbb{E}}^{1+3}$ it is convenient to use the vector fields $\Gamma_j$ defined by (1.23). We then have the following result (see [6]).

**Proposition 2.1.** *The pushforwards of $\partial_t$ and $\partial_{x_j}$ by $\mathcal{P}$ are given by*

$$\partial_t = \left(1 + \frac{1 - |U|^2}{1 + |U|^2}\cos T\right)\partial_T - \sin T\,\langle U, \partial_U\rangle \tag{2.3}$$

$$= \left(1 + \frac{|V|^2 - 1}{|V|^2 + 1}\cos T\right)\partial_T + \sin T\,\langle V, \partial_V\rangle \tag{2.4}$$

*and*

$$\partial_{x_j} = \frac{-2U_j}{1 + |U|^2}\sin T\,\partial_T + \tfrac{1}{2}\left((1 + |U|^2)\cos T + 1 - |U|^2\right)\partial_{U_j} + (1 - \cos T)\,U_j\,\langle U, \partial_U\rangle \tag{2.5}$$

$$= \frac{-2V_j}{1 + |V|^2}\sin T\,\partial_T + \tfrac{1}{2}\left((1 + |V|^2)\cos T + |V|^2 - 1\right)\partial_{V_j} + (1 + \cos T)\,V_j\,\langle V, \partial_V\rangle\,. \tag{2.6}$$

*The pushforwards via $\mathcal{P}^{-1}$ of the vector fields $\Gamma_j$ defined by (1.23) are given by*

$$X_j\,\partial_{X_k} - X_k\partial_{X_j} = x_j\,\partial_{x_k} - x_k\,\partial_{x_j}\,, \quad 1 \le j < k \le 3\,, \tag{2.7}$$

$$X_0\,\partial_{X_k} - X_k\,\partial_{X_0} = \tfrac{1}{2}(1 + t^2 - |x|^2)\,\partial_{x_k} + x_k\left(\partial_t + \langle x, \partial_x\rangle\right)\,, \quad 1 \le k \le 3\,, \tag{2.8}$$

and

$$\partial_T = \tfrac{1}{2}(1+t^2-|x|^2)\,\partial_t + t\,\langle x,\partial_x\rangle. \tag{2.9}$$

*Finally, if $\Omega = \cos T + \cos R$, then*

$$\partial_t \Omega = -\Omega\,\sin T\,\frac{1-|U|^2}{1+|U|^2} = -\Omega\,\sin T\,\cos R, \tag{2.10}$$

$$\partial_{x_j}\Omega = -\Omega\,\frac{2\cos T}{1+|U|^2}\,U_j. \tag{2.11}$$

Note that the coefficients of $\partial_T$ and $\partial_U$ in (2.3) and (2.4) are $O((\pi-T)^2)$ if $0 \le T < \pi$ and $R \le (\pi-T)/4$. Similarly, if $|x| \le t/4$, then the coefficients of $\partial_t$ and $\partial_x$ in (2.5) and (2.6) are $O(t^2+|x|^2)$. Hence we have the following useful result.

**Proposition 2.2.** *In the region where $|x| \le t/4$ we can write*

$$\partial_t = \sum a_{0k}(T,X)\,\Gamma_k, \quad \text{and} \quad \partial_{x_j} = \sum a_{jk}(T,X)\,\Gamma_k,$$

*where, if $P_0$ is as in (1.22), we have*

$$|\Gamma^\alpha a_{jk}| \le C\,\operatorname{dist}((T,X),P_0)^{2-|\alpha|}, \quad |\alpha| \le 2.$$

*Also, if $0 \le T < \pi$ and $R \le (\pi-T)/4$, then*

$$\partial_T = b_{00}(t,x)\,\partial_t + \sum b_{0k}(t,x)\,\partial_{x_k},$$

*and*

$$X_0\,\partial_{X_j} - X_j\,\partial_{X_0} = b_{0j}(t,x)\,\partial_t + \sum b_{jk}(t,x)\,\partial_{x_j},$$

*where if $\partial = (\partial_t, \partial_{x_1}, \ldots, \partial_{x_3})$ we have*

$$\bigl|\partial^\alpha b_{jk}(t,x)\bigr| \le C\bigl(1+|t|+|x|\bigr)^{2-|\alpha|}, \quad |\alpha| \le 2.$$

Using the above facts about the way that derivatives transform we shall be able to see how the nonlinear term in (1.1) transforms.

We begin by examining how the basic null forms $q_0$ and $q_{ij}$ defined by (1.7) and (1.8) transform. Let $q$ be such a null form in Minkowski space. Then if $u$ is a function on $\overline{\mathbb{E}}^4$, following (1.19), we shall let $\tilde u$ denote the function in Minkowski space [1] which is the pullback of $\Omega u$ via $\mathcal{P}$. Following (1.19) again we see that

$$\mathcal{Q}(u(T,X), du(T,X); v(T,X), dv(T,X)) = \Omega^{-3} q(d\tilde u(t,x), d\tilde v(t,x)), \quad \mathcal{P}(t,x) = (T,X)$$

is the null form transformed to $\overline{\mathbb{E}}^4$, in the sense that the following special case of (1.1)

$$\begin{cases} (\partial_t^2 - \Delta)\tilde u = q(d\tilde u, d\tilde u), & (t,x) \in \mathbb{R}_+ \times \mathbb{R}^3 \setminus \mathcal{K} \\ \tilde u(t,x)|_{\partial \mathcal{K}} = 0, \\ \tilde u(0,\cdot) = \tilde f, \quad \partial_t \tilde u(0,\cdot) = \tilde g, \end{cases}$$

---

[1] Here, as in the next several sections, we shall denote functions on Minkowski space with a tilde, while corresponding functions coming from (1.19) on the Einstein or $\mathbb{R} \times S^3$ will not have a tilde.

transforms via $\mathcal{P}$ to the following equation in $\overline{\mathbb{E}}^4 \setminus \mathcal{K}_*$

$$\begin{cases} (\Box_g + 1)u = \mathcal{Q}(u, du; u, du) \\ u(T, X) = 0, \quad (T, X) \in \partial \mathcal{K}_* \\ u(0, \cdot) = f, \quad \partial_T u(0, \cdot) = g, \end{cases}$$

if the data satisfies

$$\tilde{f} = \mathcal{P}_0^*(\Omega f), \quad \text{and} \quad \tilde{g} = \mathcal{P}_0^*(\Omega^2 g), \tag{2.12}$$

and if $\mathcal{K}_*$ is as in (1.20).

To proceed, we need the following result.

**Lemma 2.3.** *Let $q$ be any of the basic null forms defined by (1.7) or (1.8), and let $\mathcal{Q}$ be as above. Then $\mathcal{Q}$ extends to a bilinear function of $(u, du; v, dv)$ on $\mathbb{R} \times S^3$ with analytic coefficients. Moreover, if $\{\Gamma_j\}$ are defined by (1.23), we can write*

$$\mathcal{Q} = \sum_{j,k} a^{jk}(T, X) \, \Gamma_j u \, \Gamma_k v + v \sum_j b_1^j(T, X) \, \Gamma_j u + u \sum_j b_2^j(T, X) \, \Gamma_j v + c(T, X) \, u \, v, \tag{2.13}$$

*such that*

$$\Gamma^\alpha a^{jk}(P_0) = 0, \quad |\alpha| \le 1, \quad \text{and} \quad b_i^j(P_0) = 0, \tag{2.14}$$

*where $P_0 = (\pi, \mathbf{1})$ is as in (1.22).*

This result was used in the authors' earlier work [10]. The proof has two steps. The difficult step was carried out by Christodoulou [2], where it was shown that one can write $\mathcal{Q}$ as in (2.13) with the coefficients being analytic. Given this step we observe from Proposition 2.2 that if we restrict the coefficients to the region where $R \le (\pi - T)/4$, then the $a^{jk}$ must vanish to second order at $P_0$, while the $b_i^j$ must vanish there. By combining these two steps we get (2.14).

Lemma 2.3 provides the result we need for the transformation of the semilinear part $s(d\tilde{u}, d\tilde{u})$ of the nonlinear terms of our equation. We now consider the quasilinear part $k(d\tilde{u}, d^2\tilde{u})$. Recall that $k(d\tilde{u}, d^2\tilde{u})$ must be a combination of terms of the form (1.9) and (1.10). We first consider the term (1.9). If $0 \le j \le 3$ is fixed, then by Lemma 2.3 and equations (2.10) and (2.11), we can write

$$\Omega^{-3} q(d\partial_{x_j} \tilde{u}(t, x), d\tilde{v}(t, x)) = \kappa(u(T, X), du(T, X), d^2 u(T, X); v(T, X), dv(T, X)),$$

where $\kappa$, initially defined on $\overline{\mathbb{E}}^4$, extends to a bilinear function of $(u, du, d^2 u; v, dv)$ on $\mathbb{R} \times S^3$ with analytic coefficients. Moreover, we can write the extension of $\kappa$ in the form

$$\sum_{j,k} \gamma^{i,jk}(T, X) \, \Gamma_i v \, \Gamma_j \Gamma_k u + \sum_{j,k} \gamma_0^{jk}(T, X) \, v \, \Gamma_j \Gamma_k u$$
$$+ \sum_{j,k} a^{jk}(T, X) \, \Gamma_j v \, \Gamma_k u + v \sum_j b_1^j(T, X) \, \Gamma_j u + u \sum_j b_2^j(T, X) \, \Gamma_j v + c(T, X) \, u \, v,$$

where here the $a^{jk}$ and $b_i^j$ satisfy (2.14), and moreover

$$\Gamma^\alpha \gamma^{i,jk}(P_0) = 0, \quad |\alpha| \le 3, \quad \text{and} \quad \Gamma^\alpha \gamma_0^{jk}(P_0) = 0, \quad |\alpha| \le 2.$$

Lastly, to handle the quasilinear null form (1.10), we just use (1.19) and Proposition 2.1 to conclude that a term of the form $\partial_{x_j}\tilde{v}\,(\partial_t^2\tilde{u} - \Delta\tilde{u})$, also transforms into a term of the above form.

By these observations we have essentially proven the following.

**Proposition 2.4.** *Let $F(\tilde{u}, d\tilde{u}, d^2\tilde{u})$ be as in Theorem 1.1 and set*

$$\mathcal{F}(u(T,X), du(T,X), d^2u(T,X))$$
$$= \Omega^{-3} F(\tilde{u}(t,x), d\tilde{u}(t,x), d^2\tilde{u}(t,x)), \quad (T,X) = \mathcal{P}(t,x). \quad (2.15)$$

*Then $\mathcal{F}$ extends to a function of $(u, du, d^2u)$ on $\mathbb{R} \times S^3$ which is $C^\infty$ in all its variables. Moreover, if for a given $1 \leq I \leq N$ we let $\mathcal{F}^I$ be the $I$-th component of $\mathcal{F}$, then*

$$\mathcal{F}^I = \sum_{j,k} \gamma^{I,jk}(T,X;u,du)\,\Gamma_j\Gamma_k u^I + \mathcal{G}^I(T,X;u,du), \quad (2.16)$$

*where in the region $\{(T,X) : 0 \leq T < \pi,\ R \leq 2(\pi - T)\}$ if $\alpha$ is fixed one has the uniform bounds*

$$|Z^\alpha \gamma^{I,jk}| \leq C\,(\pi-T)^4 \sum_{|\gamma|\leq|\alpha|} |Z^\gamma u'| + C\,(\pi-T)^3 \sum_{|\gamma|\leq|\alpha|} |Z^\gamma u| \quad (2.17)$$
$$+ C \sum_{|\gamma_1|+|\gamma_2|\leq|\alpha|} \Big((\pi-T)^2|Z^{\gamma_1}u| + (\pi-T)^4|Z^{\gamma_1}u'|\Big)\Big((\pi-T)^2|Z^{\gamma_2}u| + (\pi-T)^4|Z^{\gamma_2}u'|\Big)$$

*and*

$$|Z^\alpha \mathcal{G}| \leq C(\pi-T)^2 \sum_{|\gamma_1|+|\gamma_2|\leq|\alpha|} |Z^{\gamma_1}u'|\,|Z^{\gamma_2}u'| + C \sum_{|\gamma_1|+|\gamma_2|\leq|\alpha|} (\pi-T)|Z^{\gamma_1}u'|\,|Z^{\gamma_2}u|$$
$$+ C \sum_{|\gamma_1|+|\gamma_2|\leq|\alpha|} |Z^{\gamma_1}u|\,|Z^{\gamma_2}u| + C \sum_{|\gamma_1|+|\gamma_2|+|\gamma_3|\leq|\alpha|} |Z^{\gamma_1}u|\,|Z^{\gamma_2}u|\,|Z^{\gamma_3}u|, \quad (2.18)$$

*assuming in both cases that*

$$(\pi-T)^2 \sum_{|\gamma|\leq 1+|\alpha|/2} |Z^\gamma u| \leq B, \quad (2.19)$$

*where $B$ is a fixed constant. Here, as before, $Z^\alpha = ((\pi-T)^2\Gamma)^\alpha$.*

*We also have the following bounds,*

$$\big|\mathcal{G}^I(T,X;u,du) - \mathcal{G}^I(T,X;v,dv)\big| \leq C\,(|u|+|v|)\Big((\pi-T)|u'-v'| + |u-v|\Big)$$
$$+ C\,(|u'|+|v'|)\Big((\pi-T)^2|u'-v'| + (\pi-T)|u-v|\Big), \quad (2.20)$$

*and*

$$\big|\gamma^{I,jk}(T,X;u,du) - \gamma^{I,jk}(T,X;v,dv)\big| \leq C\Big((\pi-T)^2|u'-v'| + |u-v|\Big), \quad (2.21)$$

*assuming that condition (2.19) holds with $|\alpha| = 0$.*

If $R(u, du, d^2u) \equiv 0$ in (1.5), then these results follow from our earlier bounds for the transformed semilinear and quasilinear quadratic terms. On the other hand, since $R(u, du, d^2u)$ is linear and diagonal in the second derivatives of $u$, and since it satisfies (1.6), it follows that if $Q(du, d^2u) \equiv 0$ in (1.5) then the above bounds must hold if (2.19) holds. Indeed, we can write

$$R^I(\tilde{u}, d\tilde{u}, d^2\tilde{u}) = O\big(|\tilde{u}|^3 + |d\tilde{u}|^3\big) + \sum \tilde{r}^{I,jk}(\tilde{u}, d\tilde{u})\, \partial_j \partial_k \tilde{u},$$

where $\tilde{r}^{I,jk} = O\big(|\tilde{u}|^2 + |\partial_{t,x}\tilde{u}|^2\big)$. The semilinear part of the remainder is controlled by the last term in (2.18), since its transformed version must be $O\big(|u|^3\big) + O\big((\pi-T)^6 \sum |\Gamma_j u|^3\big)$ on the region in the diamond $R \leq 2(\pi-T)$. Likewise, the quasilinear part of the remainder is controlled by the last term in (2.17), since it transforms to $\sum r^{I,jk}(u, du)\, \Gamma_j \Gamma_k u$ where

$$r^{I,jk}(u, du) = O\big((\pi-T)^4 |u|^2 + (\pi-T)^8 \sum |\Gamma_j u|^2\big) \quad \text{if} \quad R \leq 2(\pi-T). \quad \square$$

We now recall standard facts about how the Sobolev spaces in (1.11) transform under $\mathcal{P}_0$. Recall that the inverse of $\mathcal{P}_0$ is the south pole stereographic projection map, and so $\mathcal{P}_0(\mathcal{K}) \subset S^3$ is star-shaped with respect to the north pole and has smooth boundary. For $m = 1, 2, \ldots$ we then let

$$H_D^m(S^3 \backslash \mathcal{P}_0(\mathcal{K})) = \{f \in H^m(S^3 \backslash \mathcal{P}_0(\mathcal{K})) : f|_{\partial \mathcal{P}_0(\mathcal{K})=0}\},$$

with $H^m(S^3 \backslash \mathcal{P}_0(\mathcal{K}))$ being the Sobolev space of restrictions of elements of $H^m(S^3)$.

If then $\mathcal{P}_0^* f$ denotes the pullback of the function $f$ on $S^3 \backslash \mathcal{P}_0(\mathcal{K})$ via $\mathcal{P}_0$, and we relate $\tilde{f}$ to $f$ via

$$\tilde{f} = \Omega \mathcal{P}_0^* f$$

then the map $\tilde{f} \to f$ is continuous from $H^{m,m-1}(\mathbb{R}^3 \backslash \mathcal{K})$ to $H^m(S^3 \backslash \mathcal{P}_0(\mathcal{K}))$. That is, for fixed $m$ there is a constant $C_m$ so that

$$\|f\|_{H^m(S^3 \backslash \mathcal{P}_0(\mathcal{K}))} \leq C_m \|\Omega \mathcal{P}_0^* f\|_{H^{m,m-1}(\mathbb{R}^3 \backslash \mathcal{K})} = C_m \|\tilde{f}\|_{H^{m,m-1}(\mathbb{R}^3 \backslash \mathcal{K})}.$$

Thus, if $\tilde{u} = \Omega \mathcal{P}^* u$, then

$$\|u(0, \cdot)\|_{H^m(S^3 \backslash \mathcal{P}_0(\mathcal{K}))} + \|\partial_T u(0, \cdot)\|_{H^{m-1}(S^3 \backslash \mathcal{P}_0(\mathcal{K}))}$$
$$\leq C_m \|\tilde{u}(0, \cdot)\|_{H^{m,m-1}(\mathbb{R}^3 \backslash \mathcal{K})} + C_m \|\partial_t \tilde{u}(0, \cdot)\|_{H^{m-1,m}(\mathbb{R}^3 \backslash \mathcal{K})}, \quad (2.22)$$

since the pushforward of $\partial_t$ is $\Omega \partial_T$ if $t = 0$, and since $\Omega = 2/(1+|x|^2)$ if $t = 0$.

We close this section by presenting some of the notation that we shall use in the rest of the paper. First of all we shall let

$$Y = \big([0, \pi) \times S^3\big) \backslash \mathcal{K}_*, \quad (2.23)$$

where $\mathcal{K}_* = \mathcal{P}(\mathcal{K})$. Thus, $Y$ is the image of Minkowski space minus the obstacle. Also, for each fixed $0 \leq T < \pi$, we let

$$Y_T = \{X \in S^3 : (T, X) \in Y\} \quad (2.24)$$

be the $T$ cross section of $Y$.

Next, by dilating the Minkowski variables if necessary, we will assume that

$$\partial \mathcal{K} \subset \{x \in \mathbb{R}^3 : |x| < 1/4\}. \quad (2.25)$$

If we then let

$$r = r(T, X) = \frac{\sin(R)}{\cos T + \cos R}$$

be the spatial radial component of $\mathcal{P}^{-1}(T,X)$, it follows that

$$\partial Y \subset \mathcal{B}_{1/4}, \tag{2.26}$$

where we define

$$\mathcal{B}_r = \{(T,X) \in [0,\pi) \times S^3 : r(T,X) < r\}. \tag{2.27}$$

Thus, $\mathcal{B}_r$ is the pushforward of the cylinder $\mathbb{R}_+ \times \{x : |x| < r\}$ via $\mathcal{P}$. Equation (1.15) implies that, if $r_0 > 0$ is fixed, and if $0 \leq r \leq r_0$, then there is a uniform constant $C = C(r_0)$ so that

$$C^{-1} r (\pi - T)^2 \leq R \leq C r (\pi - T)^2, \quad \text{if} \quad (T,X) \in \partial \mathcal{B}_r. \tag{2.28}$$

As before, $R$ denotes the north pole distance on $S^3$ measured with respect to the standard metric. If $0 \leq T < \pi$, we shall let $\mathcal{B}_r^T$ denote the $T$ cross section of $\mathcal{B}_r$,

$$\mathcal{B}_r^T = \{X \in S^3 : (T,X) \in \mathcal{B}_R\}. \tag{2.29}$$

## 3. First order estimates.

Let $\Box_g = \partial_T^2 - \Delta_g$ be the wave operator on $\mathbb{R} \times S^3$, where $g_{jk}(X)\,dX_j\,dX_k$ is the standard metric on $S^3$. In this section we shall prove energy estimates for certain perturbations of $\Box_g$ in $Y$, where as in (2.23) we let $Y$ denote the image of Minkowski space minus the strictly star-shaped obstacle.

Before proving energy estimates for perturbations of $\Box_g$, we first handle $\Box_g$ itself since the arguments in this case are simpler and serve as a model for the more technical case involving perturbations. The argument that we shall use is similar to that of Morawetz [18] (see also [17], p. 261-264) for a related energy-decay estimate in Minkowski space minus a star-shaped obstacle. In particular, we shall see that when one goes through the standard proof of energy estimates the (variable) boundary contributes a term with the "correct" sign if $\mathcal{K}$ is strictly star-shaped with respect to the origin. When we handle perturbations of $\Box_g$ there will be additional boundary terms coming from the perturbation, but these will be absorbed by the Morawetz term under smallness assumptions for the perturbation.

We introduce the energy momentum 4-vector $e$ associated to a function $u(T,X)$ on $\mathbb{R} \times S^3$, defined by

$$e_0 = |\partial_T u|^2 + \|\operatorname{grad} u\|^2$$

$$e_j = -2(\partial_T u)\operatorname{grad} u, \quad 1 \leq j \leq 3$$

where grad and $\|\cdot\|$ are associated to the metric $g$. In local coordinates,

$$(\operatorname{grad} u)_j = \sum_{k=1}^{3} g^{jk}(X)\,\partial_k u(T,X),$$

$$\|\operatorname{grad} u\|^2 = \sum_{k=1}^{3} g^{jk}(X)\,\partial_j u(T,X)\,\partial_k u(T,X).$$

For convenience in future use, we will use the abbreviation

$$|u'(T,X)|^2 \equiv e_0(T,X).$$

**Theorem 3.1.** *Suppose that $\mathcal{K} \subset \mathbb{R}^3$ is strictly star-shaped with respect to the origin. Suppose also that $u \in C^2$ and $u(T, X) = 0$ if $(T, X) \in \partial Y$, and let*
$$F = \Box_g u.$$
*Then for $0 < T < \pi$*
$$\|u'(T, \cdot)\|_{L^2(Y_T)} \leq \|u'(0, \cdot)\|_{L^2(Y_0)} + \int_0^T \|F(S, \cdot)\|_{L^2(Y_S)} \, dS. \tag{3.1}$$

Here the $L^2$-norms are taken with respect to the volume element arising from the metric $g$.

We prove Theorem 3.1 by applying the divergence theorem to the vector field $e$ on $Y \cap [0, T] \times S^3$. Precisely, we consider $\mathbb{R} \times S^3$ as a Riemannian manifold with metric $dT^2 + g$. In local coordinates, the divergence of $e$ then equals
$$\partial_T e_0 + \frac{1}{\sqrt{|g|}} \sum_{j=1}^3 \partial_j \big( \sqrt{|g|} \, e_j \big) = 2 (\partial_T u) F \,.$$

The divergence theorem yields:
$$\int_{Y_T} e_0(T, \cdot) \, dX - \int_{Y_0} e_0(0, \cdot) \, dX + \int_{\partial Y \cap [0,T] \times S^3} \langle \nu, e \rangle \, d\sigma = 2 \int_{Y \cap [0,T] \times S^3} (\partial_T u) F \, dT \, dX \,.$$

Here, $\nu$ denotes the outer unit normal to $\partial Y$ in the metric $dT^2 + g$. We write
$$\nu = (\nu_T, \nu_X), \qquad (\nu_T)^2 + \|\nu_X\|^2 = 1 \,. \tag{3.2}$$

The vector $\big(1, -\nu_T \nu_X / \|\nu_X\|^2\big)$ is tangent to $\partial Y$, and by the Dirichlet conditions for $u$ we thus have
$$\partial_T u = \nu_T \|\nu_X\|^{-2} \partial_{\nu_X} u \quad \text{for} \quad (T, X) \in \partial Y \,. \tag{3.3}$$

Combining (3.2) and (3.3) shows that
$$(\partial_\nu u)^2 = \big(\nu_T \partial_T u + \partial_{\nu_X} u\big)^2 = \|\nu_X\|^{-4} (\partial_{\nu_X} u)^2 \,. \tag{3.4}$$

By Dirichlet conditions on $u$,
$$e_0(T, X) = (\partial_\nu u)^2, \qquad \text{for} \quad (T, X) \in \partial Y \,.$$

Combining (3.2), (3.3), and (3.4) yields
$$\langle \nu, e \rangle = \nu_T (\partial_\nu u)^2 - 2(\partial_T u)(\partial_{\nu_X} u) = -\nu_T \big(\|\nu_X\|^2 - \nu_T^2\big)(\partial_\nu u)^2 \,. \tag{3.5}$$

Consequently,
$$\int_{Y_T} e_0(T, \cdot) \, dX - \int_{Y_0} e_0(0, \cdot) \, dX - \int_{\partial Y \cap [0,T] \times S^3} \nu_T \big(\|\nu_X\|^2 - \nu_T^2\big) e_0(T, X) \, d\sigma$$
$$= 2 \int_{Y \cap [0,T] \times S^3} (\partial_T u) F \, dT \, dX \,.$$

The important observation now is that $\nu_T$ is strictly negative. Indeed, working in polar coordinates about the north pole, by (1.2) we can write
$$\partial \mathcal{K}_* = \big\{ \big( T, \Phi(T, \omega) \big) : 0 \leq T < \pi \big\},$$

with $\Phi$ smooth. The crucial fact is that

$$\partial_T \Phi < 0, \quad 0 < T < \pi. \tag{3.6}$$

Indeed,

$$\partial_T \Phi = -4t\phi(r,\omega) \, / \, [(t + \phi(r,\omega))^2 + (t - \phi(r,\omega))^2],$$

and so

$$\partial_T \Phi \leq -c\min\{(\pi - T), T\}, \quad 0 \leq T < \pi, \tag{3.7}$$

for some fixed constant $c > 0$. These facts follow from our strict star-shaped hypothesis (1.2) and an elementary calculation.

From equation (3.7) and the fact that $\partial Y$ is timelike, there is a uniform constant $c > 0$ so that

$$-\nu_T\bigl(\|\nu_X\|^2 - \nu_T^2\bigr) \geq c\min\{(\pi - T), T\}, \quad 0 \leq T < \pi. \tag{3.8}$$

An application of the Gronwall inequality completes the proof of the theorem.

**Energy estimates for perturbed operators**

In this section we work with a Lorentzian metric $h$ which we shall assume to be a small perturbation of the standard Lorentz metric $\eta$ defined by

$$\eta = dT^2 - g.$$

We let $\square_h$ denote the associated D'Alembertian, which, in local coordinates takes the form

$$\square_h u = |h|^{-1/2} \sum_{j,k=0}^{3} \partial_j(h^{jk}|h|^{1/2}\partial_k u), \tag{3.9}$$

where $(h^{jk}) = (h_{jk})^{-1}$.

We will assume that $h$ is uniformly close to the standard metric,

$$\bigl|h(V,W) - \eta(V,W)\bigr| \leq \delta, \tag{3.10}$$

for all pairs of vectors $V$, $W$ of norm one in the metric $dT^2 + g$. We will take $\delta$ sufficiently small (to be determined.) We shall also assume that we have the following bounds for the covariant derivatives of $h - \eta$ with respect to $dT^2 + g$,

$$\|\nabla(h - \eta)\|_{L^1_T L^\infty_X} \leq C_0. \tag{3.11}$$

It will be convenient to use local coordinates in our calculations; we thus cover the sphere with two compact coordinate patches using north pole and south pole projective coordinates. We then write

$$h^{jk} = \eta^{jk} + \gamma^{jk}.$$

Expressed in these coordinate systems, our conditions are equivalent (up to constants) to the following

$$\sum_{j,k=0}^{3} \|\gamma^{jk}(T,X)\|_{L^\infty_{T,X}} \leq \delta\,,$$

$$\sum_{i,j,k=0}^{3} \|\partial_i \gamma^{jk}(T,X)\|_{L^1_T L^\infty_X} \leq C_0\,. \tag{3.12}$$

We then have the following

**Theorem 3.2.** *Assume $h$ is as above and let $u \in C^2$ satisfy*

$$\begin{cases} \Box_h u(T,X) = F(T,X)\,, & (T,X) \in Y \\ u(T,X) = 0\,, & (T,X) \in \partial Y\,. \end{cases}$$

*Then if (3.10) holds for $\delta > 0$ sufficiently small, and if (3.11) holds, then*

$$\|u'(T,\cdot)\|_{L^2(Y_T)} \leq C\|u'(0,\cdot)\|_{L^2(Y_0)} + C\int_0^T \|F(S,\cdot)\|_{L^2(Y_S)}\,dS\,, \quad 0 < T < \pi\,, \tag{3.13}$$

*for a uniform constant $C$ (depending on $C_0$).*

To proceed, we introduce the energy-momentum vector

$$\tilde{e} = 2(\partial_T u)\operatorname{grad}_h u - \langle \operatorname{grad}_h u, \operatorname{grad}_h u \rangle \partial_T\,,$$

where $\operatorname{grad}_h$ and $\langle \cdot, \cdot \rangle_h$ denote the Lorentz gradient and inner product for $h$. In local coordinates,

$$\tilde{e}_0 = 2(\partial_T u)\sum_{k=0}^{3} h^{0k}(T,X)\,\partial_k u(T,X) - \sum_{j,k=0}^{3} h^{jk}(T,X)\partial_j u(T,X)\,\partial_k u(T,X)\,,$$

$$\tilde{e}_j = 2(\partial_T u)\sum_{k=0}^{3} h^{jk}(T,X)\,\partial_k u(T,X)\,, \quad j = 1,2,3.$$

We now apply the divergence theorem on $Y \cap [0,T] \times S^3$ using the divergence $(\partial_T, \operatorname{div})$ associated to the standard Riemannian metric $dT^2 + g$. We first claim that

$$\partial_T \tilde{e}_0 + \operatorname{div} \tilde{e}_X = 2(\partial_T u)\Box_h u + R(u',u')\,,$$

where $R(\cdot,\cdot)$ is a quadratic form whose coefficients (in any orthonormal frame) belong to $L^1_T L^\infty_X$ with norm bounded by some fixed multiple of $C_0$.

To see this, we work in local coordinates. There, we may write

$$\operatorname{div} \tilde{e}_X = \frac{1}{\sqrt{|g|}}\sum_{j=1}^{3} \partial_j\bigl(\sqrt{|g|}\,\tilde{e}_j\bigr) = \sum_{j=1}^{3} \partial_j \tilde{e}_j + \sum_{j=1}^{3} r_j \tilde{e}_j\,,$$

where the $r_j$ are uniformly bounded functions.

Next, a simple calculation shows that

$$\sum_{j=0}^{3} \partial_j \tilde{e}_j = 2(\partial_T u)\sum_{j,k=0}^{3} \partial_j\Bigl(h^{jk}\partial_k u\Bigr) + \sum_{j,k=0}^{3} r_{jk}\partial_j u \partial_k u$$

where the coefficients $r_{jk}(T, X)$ involve first derivatives of the $h^{jk}$, and hence belong to $L_T^1 L_X^\infty$. The last expression may be written in the form

$$2(\partial_T u)\Box_h u + R(u', u').$$

Next, we claim that

$$\left|\langle \nu, \tilde{e} - e\rangle\right| \leq C\,\delta\,|\nu_T|\,|u'|^2.$$

(The inner product is with respect to $dT^2 + g$.) To see this, we write

$$\langle \nu, \tilde{e}-e\rangle = -\nu_T\left(\sum_{j,k=0}^{3} \gamma^{jk}\partial_j u\,\partial_k u\right) + 2\nu_T\,(\partial_T u)\sum_{k=0}^{3}\gamma^{0k}\partial_k u + 2(\partial_T u)\sum_{i,j=1}^{3}\sum_{k=0}^{3} g_{ij}\,\gamma^{jk}\partial_k u.$$

The first two terms clearly have the desired bounds; to handle the last, we use (3.3).

It follows from (3.5) and (3.8) that for $\delta$ sufficiently small, the boundary term $\langle \nu, \tilde{e}\rangle$ is positive, and by the divergence theorem we have

$$\int_{Y_T} \tilde{e}_0(T, X)\,dX$$
$$\leq \int_{Y_0} \tilde{e}_0(0, X)\,dX + 2\int_{Y\cap[0,T]\times S^3}(\partial_T u)F\,dT\,dX + \int_{Y\cap[0,T]\times S^3} R(u', u')\,dT\,dX.$$

The proof of the theorem now follows from the Gronwall lemma by noting that, for $\delta$ small, we have

$$C^{-1}|u'(T, X)|^2 \leq \tilde{e}_0(T, X) \leq C\,|u'(T, X)|^2.$$

In Section 5 we shall use the fact that (3.13) holds for equations of the form

$$\Bigl(\Box_g + \sum_{j,k=0}^{6}\gamma^{jk}\,\Gamma_j\,\Gamma_k + 1\Bigr)u = F,$$

provided that

$$\sum_{j,k=0}^{6}\|\gamma^{jk}(T, X)\|_{L_{T,X}^\infty} \leq \delta, \tag{3.14}$$

$$\sum_{i,j,k=0}^{6}\|\Gamma_i\gamma^{jk}(T, X)\|_{L_T^1 L_X^\infty} \leq C_0. \tag{3.15}$$

To see this, we note that there is a unique metric $h$ such that the operator

$$\Box_h - \Box_g - \sum_{j,k=0}^{6}\gamma^{jk}\,\Gamma_j\,\Gamma_k$$

is of first order. Furthermore, the metric $h$ satisfies the conditions (3.10) and (3.11) (with possibly different constants.) We now just observe that the proof of Theorem 3.2 goes through if $F$ is modified by first order derivatives in $u$.

4. **Sobolev Estimates.**

For our applications in later sections, in addition to controlling the $L^2$ norm of $u'(T, \cdot)$, we also need to control the $L^6$ norm of $u(T, \cdot)$. For this we will make use of the following elementary result.

**Lemma 4.1.** *Suppose that $u \in C^1$ vanishes on $\partial Y$. Then there is a uniform constant $C$ so that for $0 \leq T < \pi$*
$$\|u(T, \cdot)\|_{L^6(Y_T)} \leq C\|u'(T, \cdot)\|_{L^2(Y_T)} + C\|u(T, \cdot)\|_{L^2(Y_T)}. \tag{4.1}$$

The lemma follows from the Sobolev embedding $H^1 \hookrightarrow L^6$ for the sphere $S^3$ by noting that the extension of $u(T, \cdot)$ to the entire sphere, obtained by setting it equal to 0 on the complement of $Y_T$, belongs to $H^1(S^3)$ with norm controlled by the right hand side of (4.1).

Using the lemma we can strengthen (3.13) somewhat.

**Corollary 4.2.** *Let $u$ be as in Theorem 3.2. Then*
$$\|u'(T, \cdot)\|_{L^2(Y_T)} + \|u(T, \cdot)\|_{L^6(Y_T)}$$
$$\leq C\|u'(0, \cdot)\|_{L^2(Y_0)} + C \int_0^T \|F(S, \cdot)\|_{L^2(Y_S)} \, dS, \quad 0 < T < \pi. \tag{4.2}$$

This follows from (3.13) and Lemma 4.1 above by noting that, since $u$ satisfies Dirichlet conditions on $\partial Y$, we can bound
$$\|u(T, \cdot)\|_{L^2(Y_T)} \leq \int_0^T \|u'(S, \cdot)\|_{L^2(Y_S)} \, dS + \|u(0, \cdot)\|_{L^2(Y_0)}$$
$$\leq \sup_{0 \leq S \leq T} \|u'(S, \cdot)\|_{L^2(Y_S)} + C \|u'(0, \cdot)\|_{L^2(Y_0)}.$$

In the arguments of the next section where we control higher derivatives of the solution $u$, we will need an elliptic regularity result for a perturbation of the operator $\Delta_g$ on the 3-sphere. Precisely, we will work with the operator
$$\Delta_g - a(T, X)\partial_R^2,$$
where
$$a(T, X) = \frac{\sin^2 T \sin^2 R}{(1 + \cos T \cos R)^2} = \frac{\sin^2 T \sin^2 R}{(1 + \cos(T + R) + \sin T \sin R)^2}. \tag{4.3}$$

We use the fact that
$$a(T, X) = \frac{\sin^2(\pi - T) \sin^2 R}{\left(1 - \cos(\pi - T) \cos R\right)^2} \leq \left(\frac{2\delta}{1 + \delta^2}\right)^2$$

if $\sin R = \delta \sin(\pi - T)$, for $\delta \leq 1$. Consequently, $\Delta_g - a(t, X)\partial_R^2$ is uniformly elliptic on the set $R < \delta(\pi - T)$ if $\delta < 1$. Also, from the fact that
$$1 - \cos(\pi - T) \cos R \approx \frac{1}{2}(\pi - T)^2 + \frac{1}{2}R^2,$$
for $(T, R)$ near $(\pi, 0)$, it is easy to see that
$$\left|\Gamma^\alpha a(T, X)\right| \leq C_\alpha \left[(\pi - T)^2 + R^2\right]^{-|\alpha|/2}. \tag{4.4}$$

Our estimate will involve the weighted derivatives $Z = (\pi - T)^2 \Gamma$ as in (1.24); we let

$$\{Z_X\} = \{Z_{jk} : 1 \le j < k \le 6\}$$

be the set of weighted derivatives that do not involve $\partial_T$, and similarly define $\Gamma_X$.

**Proposition 4.3.** *Let $k = 0, 1, 2, \ldots$. Then there is a constant $C = C(k)$, independent of $T$, so that whenever $h \in C^\infty(Y)$ vanishes on $\partial Y$, then for $0 \le T < \pi$,*

$$\sum_{|\alpha| \le k+1} \Big( \|Z_X^\alpha \Gamma_X h(T, \cdot)\|_{L^2(R < (\pi-T)/2)} + \|Z_X^\alpha h\|_{L^6(R < (\pi-T)/2)} \Big)$$

$$\le C (\pi - T)^2 \sum_{|\alpha| \le k} \|Z_X^\alpha (\Delta_g - a(T, X)\partial_R^2) h(T, \cdot)\|_{L^2(R < \pi - T)}$$

$$+ C \sum_{|\alpha| \le k} \|Z_X^\alpha h(T, \cdot)\|_{L^6(R < \pi - T)}. \quad (4.5)$$

We first show that the estimate holds if the norms on the left hand side are taken over a set of the form $R \le C_0 (\pi - T)^2$. To do this, we work in south pole stereographic coordinates $U$ on $S^3$, which map the north pole to the origin, and dilate in these variables by $(\pi - T)^2$. After dilation, the boundary $\partial Y_T$ is mapped to a surface $M_T \subset \mathbb{R}^3$ contained in the set $c \le R \le c^{-1}$, where $c > 0$ is independent of $T$, such that there are uniform bounds on the surface $M_T$ independent of $T$.

We next write $\Delta_g - a(T, X)\partial_R^2$ as $L_T(U, D_U)$ in the stereographic coordinates. Then the operator $P_T = L_T\big((\pi-T)^2 U, \partial_U\big)$ is seen to be a uniformly elliptic operator on the image of the set $R < \delta(\pi - T)$ for any $\delta < 1$, and furthermore there are uniform bounds on the derivatives of the coefficients of $P_T$ which are independent of $0 < T < \pi$. In fact, this statement is true for $L_T\big((\pi - T)U, \partial_U\big)$, which follows from (4.4), together with the fact that $a(T, X)$ vanishes quadratically at $R = 0$, and the fact that $\partial_R$ is mapped under stereographic coordinates to a smooth multiple of the radial vector field.

Because $Z$ scales to a unit vector field under this dilation, the desired estimate is a result of the following estimate in the scaled coordinates, for functions $f \in C^\infty(M_T^{\text{ext}})$ which vanish on $M_T$,

$$\sum_{|\alpha| \le k+1} \Big( \|\partial^\alpha (\nabla f)\|_{L^2(r \le C_0)} + \|\partial^\alpha f\|_{L^6(r \le C_0)} \Big)$$

$$\le C \sum_{|\alpha| \le k} \Big( \|\partial^\alpha (P_T f)\|_{L^2(r \le 2C_0)} + \|\partial^\alpha f\|_{L^6(r \le 2C_0)} \Big),$$

and this estimate holds by standard elliptic regularity theory.

We remark that this proof in fact shows that, to control the left hand side of (4.5) over the set $R < C_0(\pi - T)^2$, it suffices to take the norms on the right hand side over the set $R \le 2C_0(\pi - T)^2$, a fact we will use in the proof of Proposition 6.2.

Now let $\phi$ be a cutoff to the set $R \geq c\,(\pi - T)^2$, where $c$ is chosen so that $\phi = 0$ on $\partial Y$. From the fact that $\|Z^\alpha \nabla_X \phi(T, \cdot)\|_{L^3} \leq C$, and the estimates (4.4), it follows that

$$(\pi - T)^2 \sum_{|\alpha| \leq k} \|Z_X^\alpha (\Delta_g - a(T,X)\partial_R^2)(\phi h)(T,\cdot)\|_{L^2(R \leq \pi - T)}$$
$$\leq C\,(\pi - T)^2 \sum_{|\alpha| \leq k} \|Z_X^\alpha (\Delta_g - a(T,X)\partial_R^2) h(T,\cdot)\|_{L^2(R < \pi - T)}$$
$$+ C \sum_{|\alpha| \leq k+1} \|Z_X^\alpha h(T,\cdot)\|_{L^6(R < c(\pi - T)^2)}.$$

By the preceding steps the last term is controlled by the right hand side of (4.5); consequently we are reduced to the case of establishing (4.5) in the absence of a boundary.

To see that (4.5) holds in the absence of a boundary, we again work in south pole stereographic coordinates, and now dilate by $(\pi - T)$, so that the set $R \leq (\pi - T)/2$ is mapped to a ball of radius close to $1/2$. We now use the fact that $P_T = L_T\big((\pi - T)U, \partial_U\big)$ is uniformly elliptic on the region of interest, with smooth coefficients that have uniform bounds on $0 < T < \pi$.

Next, by an induction argument we may consider just the terms on the left hand side of (4.5) where $|\alpha| = k + 1$. Then, after scaling, we are led to the estimate

$$(\pi - T)^k \sum_{|\alpha| = k+1} \Big( \|\partial^\alpha(\nabla f)\|_{L^2(r \leq \frac{1}{2})} + \|\partial^\alpha f\|_{L^6(r \leq \frac{1}{2})} \Big)$$
$$\leq C \sum_{|\alpha| \leq k} (\pi - T)^{|\alpha|} \Big( \|\partial^\alpha(P_T f)\|_{L^2(r \leq 1)} + (\pi - T)^{-1} \|\partial^\alpha f\|_{L^6(r \leq 1)} \Big).$$

Since the powers of $(\pi - T)$ on the right are less than or equal to $k$, and $(\pi - T)$ is bounded above, this estimate follows as before by elliptic regularity theory. $\square$

## 5. Higher Order Estimates.

In this section, we establish a priori estimates on higher order weighted derivatives of the solution $u$, in terms of weighted derivatives of the coefficients $\gamma^{jk}$. For convenience, we assume that $\gamma^{jk}$ and $u$ belong to $C^\infty(Y)$, where we recall that $Y = ([0, \pi) \times S^3) \backslash \mathcal{K}_*$.

We shall also assume that the $\gamma^{jk}$ satisfy the hypotheses (3.14) and (3.15) so that by Theorem 3.2 we have control of the $L^2$-norm of $u'$.

**Theorem 5.1.** *Suppose that $\gamma^{jk}, u \in C^\infty(Y)$, and that $u(T, X) = 0$ if $(T, X) \in \partial Y$. Suppose that the $\gamma^{jk}$ satisfy (3.14) and (3.15), where $\delta > 0$ in (3.14) is small enough so that (3.13) holds. Let*

$$F = \big(\Box_g + \sum \gamma^{jk} \Gamma_j \Gamma_k + 1\big) u\,.$$

Then, given $N = 0, 1, 2, \ldots$, there is a constant $C$ depending only on $N$, $\delta$, and $C_0$, so that for $0 \leq T < \pi$,

$$\sum_{|\alpha| \leq N} \left( \|Z^\alpha u'(T, \cdot)\|_2 + \|Z^\alpha u(T, \cdot)\|_6 \right) \leq C \sum_{|\alpha| \leq N} \left( \|Z^\alpha u'(0, \cdot)\|_2 + \int_0^T \|Z^\alpha F(S, \cdot)\|_2 \, dS \right)$$

$$+ C \sup_{0 \leq S \leq T} (\pi - S)^2 \sum_{|\alpha| \leq N-1} \|Z^\alpha F(S, \cdot)\|_2$$

$$+ C \int_0^T (\pi - S)^{-2} \sum_{j,k} \sum_{\substack{|\alpha_1| + |\alpha_2| \leq N+1 \\ |\alpha_1|, |\alpha_2| \leq N}} \left( \|(Z^{\alpha_1} \gamma^{jk}) Z^{\alpha_2} u'(S, \cdot)\|_2 + \|(Z^{\alpha_1} \gamma^{jk}) Z^{\alpha_2} u(S, \cdot)\|_6 \right) dS$$

$$+ C \sup_{0 \leq S \leq T} \sum_{j,k} \sum_{\substack{|\alpha_1| + |\alpha_2| \leq N \\ |\alpha_1| \leq N-1}} \|(Z^{\alpha_1} \gamma^{jk}) Z^{\alpha_2} u'(S, \cdot)\|_2. \quad (5.1)$$

**Remark.** Before proving Theorem 5.1, we point out that if we fix $N$ and assume that $\delta > 0$ is small enough so that $C \sum_{j,k=0}^3 \|\gamma^{jk}\|_\infty < 1/2$ then we can strengthen (5.1) somewhat. Specifically, the part of the last summand in the right side of (5.1) where $\alpha_1 = 0$ and $|\alpha_2| = N$ can be absorbed in the left side of (5.1). As a result, under this additional smallness assumption, we have

$$\sum_{|\alpha| \leq N} \left( \|Z^\alpha u'(T, \cdot)\|_2 + \|Z^\alpha u(T, \cdot)\|_6 \right) \leq C \sum_{|\alpha| \leq N} \left( \|Z^\alpha u'(0, \cdot)\|_2 + \int_0^T \|Z^\alpha F(S, \cdot)\|_2 \, dS \right)$$

$$+ C \sup_{0 \leq S \leq T} (\pi - S)^2 \sum_{|\alpha| \leq N-1} \|Z^\alpha F(S, \cdot)\|_2$$

$$+ C \int_0^T (\pi - S)^{-2} \sum_{j,k} \sum_{\substack{|\alpha_1| + |\alpha_2| \leq N+1 \\ |\alpha_1|, |\alpha_2| \leq N}} \left( \|(Z^{\alpha_1} \gamma^{jk}) Z^{\alpha_2} u'(S, \cdot)\|_2 + \|(Z^{\alpha_1} \gamma^{jk}) Z^{\alpha_2} u(S, \cdot)\|_6 \right) dS$$

$$+ C \sup_{0 \leq S \leq T} \sum_{j,k} \sum_{\substack{|\alpha_1| + |\alpha_2| \leq N \\ |\alpha_1|, |\alpha_2| \leq N-1}} \|(Z^{\alpha_1} \gamma^{jk}) Z^{\alpha_2} u'(S, \cdot)\|_2.$$

We will establish Theorem 5.1 by induction. The inequality holds with $N = 0$, by Theorem 3.2 and the remark at the end of Section 3. For $N = 0$ only the first two terms on the right are needed. We thus make the following

**Induction hypothesis:** Inequality (5.1) is valid if $N$ is replaced by $N - 1$.

We then show that this implies (5.1) for $N = N$. To do this, we write the norm on the left hand side of (5.1) as a sum of two terms, by separately considering the regions $R < (\pi - T)/2$ and $R > (\pi - T)/2$. We begin by considering $R < (\pi - T)/2$. To estimate this term, we will make use of the vector field $\mathcal{X}$ on $Y$ obtained by pushing forward the Minkowski time derivative via the Penrose compactification,

$$\mathcal{X} = \mathcal{P}_*(\partial_t).$$

We note that $\mathcal{X}$ is tangent to $\partial Y$, so that if $u$ vanishes on $\partial Y$ then so does $\mathcal{X}u$, that is,

$$\mathcal{X}u(T, X) = 0 \text{ if } (T, X) \in \partial Y. \quad (5.2)$$

GLOBAL EXISTENCE FOR NONLINEAR WAVE EQUATIONS        21We use the following formulae involving $\mathcal{X}$:

**Lemma 5.2.** *As above, let $R$ denote the polar distance from the north pole in $S^3$. Then*

$$\mathcal{X} = (1 + \cos T \cos R)\,\partial_T - \sin T \sin R\,\partial_R \qquad (5.3)$$
$$= (1 + \cos(T+R) + \sin T \sin R)\,\partial_T - \sin T \sin R\,\partial_R\,.$$

*Moreover,*

$$[\Box_g, \mathcal{X}] = -2\cos R \sin T\,\Box_g + 2\cos R \cos T\,\partial_T + 2\sin T \sin R\,\partial_R\,. \qquad (5.4)$$

The proof of this lemma is a straightforward calculation, and will be postponed until the end of this section.

If we let $L = \Box_g + \sum \gamma^{jk}\Gamma_j\Gamma_k + 1$ then (5.4) yields

$$L\mathcal{X}u = \mathcal{X}F - 2\cos R \sin T\,(F - \sum_{j,k}\gamma^{jk}\Gamma_j\Gamma_k u - u)$$
$$+ 2\cos R \cos T\,\partial_T u + 2\sin T \sin R\,\partial_R u + \sum_{j,k}[\gamma^{jk}\Gamma_j\Gamma_k, \mathcal{X}]u\,. \qquad (5.5)$$

Now let $\phi$ be a cutoff function such that $\phi(T, X) = 1$ for $R \leq 2(\pi - T)$, and $\phi(T, X) = 0$ for $R \geq 3(\pi - T)$, and

$$Z^\alpha \phi(T, X) \leq C_\alpha\,. \qquad (5.6)$$

Let $w$ solve the equation

$$Lw = \phi\Big(\mathcal{X}F - 2\cos R \sin T\,(F - \sum_{j,k}\gamma^{jk}\Gamma_j\Gamma_k u - u)$$
$$+ 2\cos R \cos T\,\partial_T u + 2\sin T \sin R\,\partial_R u + \sum_{j,k}[\gamma^{jk}\Gamma_j\Gamma_k, \mathcal{X}]u\Big)\,. \qquad (5.7)$$

Then, by finite propagation velocity, $w = \mathcal{X}u$ for $R \leq (\pi - T)$. We will show that, if one takes the right hand side of (5.1) with $N$ replaced by $N-1$, and $F$ replaced by the right hand side of (5.7), then the result is bounded by the right hand side of (5.1) with $N = N$. The induction hypothesis, using the fact that $w = 0$ on $\partial Y$ as a result of (5.2), will then show that the following quantity is bounded by the right hand side of (5.1),

$$\sum_{|\alpha|\leq N-1} \|Z^\alpha(\mathcal{X}u)'(T,\cdot)\|_{L^2(R<(\pi-T))} + \|Z^\alpha \mathcal{X}u(T,\cdot)\|_{L^6(R<(\pi-T))}\,.$$

Since $(\mathcal{X}u)' = \mathcal{X}u' + O(u')$, we conclude that

$$\sum_{|\alpha|\leq N-1} \big(\|Z^\alpha \mathcal{X}u'(T,\cdot)\|_{L^2(R<(\pi-T))} + \|Z^\alpha \mathcal{X}u(T,\cdot)\|_{L^6(R<(\pi-T))}\big) \qquad (5.8)$$

is also bounded by the right hand side of (5.1).

To bound the right hand side of (5.1) with $N$ replaced by $N-1$ and $F$ replaced by the right hand side of (5.7), we first use (5.6) to note that it suffices to bound the same quantity with $F$ replaced by the right hand side of (5.5), but with the norms taken over the set $R \leq 3(\pi - T)$.

Next, we notice that for $R \leq 3(\pi - T)$, one can write $\mathcal{X}$ as a combination of the vector fields $\{Z\}$, with coefficients that satisfy the same estimates (4.4) as $a(T, X)$. Based on this, one can see that

$$\sum_{j,k} \sum_{|\alpha| \leq N-1} \Big( \big\| Z^\alpha [\gamma^{jk} \Gamma_j \Gamma_k, \mathcal{X}] u \big\|_{L^2(R<3(\pi-T))} + \big\| Z^\alpha (\gamma^{jk} \Gamma_j \Gamma_k\, u) \big\|_{L^2(R<3(\pi-T))} \Big)$$
$$\leq C\,(\pi - T)^{-2} \sum_{j,k} \sum_{\substack{|\alpha_1|+|\alpha_2| \leq N+1 \\ |\alpha_1|,|\alpha_2| \leq N}} \big\| (Z^{\alpha_1}\gamma^{jk}) Z^{\alpha_2} u'(T, \cdot) \big\|_2. \quad (5.9)$$

Similarly,

$$\sum_{|\alpha| \leq N-1} \| Z^\alpha \mathcal{X} F \|_{L^2(R<3(\pi-T))} + \| Z^\alpha F \|_{L^2(R<3(\pi-T))} \leq C \sum_{|\alpha| \leq N} \| Z^\alpha F \|_2.$$

These terms are thus bounded by the right hand side of (5.1).

To handle the remaining terms, which involve $u'$ and $u$, we note that

$$\sum_{|\alpha| \leq N-1} \big\| Z^\alpha \big( \cos R \cos T \partial_T u + \sin T \sin R\, \partial_R u \big) \big\|_2 \leq C \sum_{|\alpha| \leq N-1} \| Z^\alpha u' \|_2,$$

while

$$\sum_{|\alpha| \leq N-1} \big\| Z^\alpha \big( \cos R \sin T\, u \big) \big\|_2 \leq C \sum_{|\alpha| \leq N-1} \| Z^\alpha u \|_6$$

By the induction hypothesis, the norms on the right hand side of these two equations are in fact bounded by the right hand side of (5.1) with $N$ replaced by $N-1$, and thus with $N = N$.

Thus, we have shown that the quantity in (5.8) is bounded by the right hand side of (5.1). To proceed, we use (5.3) to write

$$(1 + \cos T \cos R)\big(\partial_T^2 u - a(T,X)\partial_R^2 u\big) = \mathcal{X}\partial_T u - \frac{\sin T \sin R}{1 + \cos T \cos R}\, \mathcal{X}\, \partial_R u,$$

where $a(T, X)$ is as in (4.3). Since $R \geq c\,(\pi - T)^2$ on $Y_T$, we may write $\partial_R u = b \cdot u'$ where $|Z^\alpha b| \leq C_\alpha$. Consequently, we may bound

$$\sum_{|\alpha| \leq N-1} \big| Z^\alpha \big( (1 + \cos T \cos R)(\partial_T^2 u - a(T,X)\partial_R^2 u) \big) \big| \leq C \sum_{|\alpha| \leq N-1} |Z^\alpha \mathcal{X} u'| + |Z^\alpha u'|.$$

The $L^2$ norm of the right hand side over $R < (\pi - T)$ is bounded by the quantity (5.8), and consequently the following quantity is bounded by the right hand side of (5.1),

$$\sum_{|\alpha| \leq N-1} \big\| Z^\alpha \big( (1 + \cos T \cos R)(\partial_T^2 u(T, \cdot) - a(T,X)\partial_R^2 u(T, \cdot)) \big) \big\|_{L^2(R<(\pi-T))}.$$

We next write

$$\Delta_g u - a(T,X)\partial_R^2 u = \partial_T^2 u - a(T,X)\partial_R^2 u + \sum_{j,k} \gamma^{jk}\Gamma_j\Gamma_k u + u - F.$$

Note that $1 + \cos R \cos T \approx (\pi - T)^2$ for $R < 3(\pi - T)$. Therefore
$$\sum_{j,k} \sum_{|\alpha| \leq N-1} |Z^\alpha (1 + \cos T \cos R) \gamma^{jk} \Gamma_j \Gamma_k u(T, \cdot)|$$
$$\leq C \sum_{j,k} \sum_{\substack{|\alpha_1| + |\alpha_2| \leq N \\ |\alpha_1| \leq N-1}} \left| (Z^{\alpha_1} \gamma^{jk}) Z^{\alpha_2} u'(T, \cdot) \right|.$$

Based on this and the induction hypothesis we deduce that the following quantity is bounded by the right side of (5.1)
$$\sum_{|\alpha| \leq N-1} \| Z^\alpha \left( (1 + \cos T \cos R)(\Delta_g u(T, \cdot) - a(T, X) \partial_R^2 u(T, \cdot)) \right) \|_{L^2(R < (\pi - T))},$$
and thus by Proposition 4.3 so is the following quantity
$$\sum_{|\alpha| \leq N} \| Z_X^\alpha \Gamma_X u \|_{L^2(R \leq (\pi - T)/2)} + \| Z_X^\alpha u \|_{L^6(R \leq (\pi - T)/2)}.$$

We write
$$\mathcal{X} = (1 + \cos T \cos R)(\partial_T - b \cdot \Gamma_X).$$
Again from the fact that $(1 + \cos T \cos R) \approx (\pi - T)^2$ and the fact that (5.8) is bounded, the following quantity is bounded by the right hand side of (5.1),
$$\sum_{|\alpha| \leq N} (\pi - T)^{2|\alpha|} \left( \| \Gamma^\alpha (\partial_T - b \cdot \Gamma_X) u \|_{L^2(R \leq (\pi - T)/2)} + \| \Gamma_X^\alpha \Gamma_X u \|_{L^2(R \leq (\pi - T)/2)} \right)$$
$$+ \sum_{|\alpha| \leq N-1} (\pi - T)^{2|\alpha|} \left( \| \Gamma^\alpha (\partial_T - b \cdot \Gamma_X) u \|_{L^6(R \leq (\pi - T)/2)} + \| \Gamma_X^\alpha \Gamma_X u \|_{L^6(R \leq (\pi - T)/2)} \right).$$

A simple induction in the number of $T$ derivatives shows that this in turn bounds the following quantity,
$$\sum_{|\alpha| \leq N} (\pi - T)^{2|\alpha|} \left( \| \Gamma^\alpha \Gamma u \|_{L^2(R \leq (\pi - T)/2)} + \| \Gamma^\alpha u \|_{L^6(R \leq (\pi - T)/2)} \right),$$
which is comparable to
$$\sum_{|\alpha| \leq N} \left( \| Z^\alpha u'(T, \cdot) \|_{L^2(R \leq (\pi - T)/2)} + \| Z^\alpha u(T, \cdot) \|_{L^6(R \leq (\pi - T)/2)} \right).$$

To finish, we need to show that the norm on the left hand side of (5.1), taken over the set $R \geq (\pi - T)/2$, is bounded by the right hand side of (5.1). As before let
$$\{\Gamma\} = \{\partial/\partial T, X_j \partial_k - X_k \partial_j, 1 \leq j < k \leq 4\},$$
and recall that $\Gamma^{\alpha_0}$ commutes with $\Box_g$. Therefore if
$$(\Box_g + \sum \gamma^{jk} \Gamma_j \Gamma_k + 1) v = G,$$
and if $v$ vanishes in a neighborhood of $\partial Y$, the first order energy estimate Theorem 3.2 yields
$$\| (\Gamma^{\alpha_0} v)'(T, \cdot) \|_2 + \| \Gamma^{\alpha_0} v(T, \cdot) \|_6 \leq C \| (\Gamma^{\alpha_0} v)'(0, \cdot) \|_2$$
$$+ C \int_0^T \left( \| \Gamma^{\alpha_0} G(S, \cdot) \|_2 + \sum_{j,k} \| [\Gamma^{\alpha_0}, \gamma^{jk} \Gamma_j \Gamma_k] v(S, \cdot) \|_2 \right) dS. \quad (5.10)$$

To apply this we shall let $\eta \in C^\infty$ satisfy $\eta = 1$ for $R \geq (\pi - T)/2$, and $\eta = 0$ for $R \leq (\pi - T)/3$, such that $|\Gamma^\alpha \eta| \leq C_\alpha (\pi - T)^{-|\alpha|}$ for all $\alpha$. We then will apply (5.10) to $v = \eta u$, in which case

$$G = \eta F - \sum \left[\gamma^{jk} \Gamma_k \Gamma_k, \eta\right] u + 2\partial_T \eta \, \partial_T u - 2\nabla_X \eta \cdot \nabla_X u + (\Box_g \eta) u \,.$$

We need to show that $(\pi - T)^{2|\alpha_0|}$ times the right hand side of (5.10) is bounded by the right hand side of (5.1), where $N = |\alpha_0|$. We first consider the term $G$, and take $|\alpha_0| \geq 1$, since for $\alpha_0 = 0$ the result holds by our energy estimate. To begin, note that

$$(\pi - T)^{2|\alpha_0|} \int_0^T \|\Gamma^{\alpha_0}(\eta F)(S,\cdot)\|_2 \leq C_{\alpha_0} \int_0^T \sum_{|\alpha| \leq |\alpha_0|} \|Z^\alpha F(s,\cdot)\|_2 \,.$$

We now note that the remaining terms in $G$ are supported on the set $R \leq (\pi - T)/2$. Thus, by Hölder's inequality, we may bound

$$(\pi - T)^{2|\alpha_0|} \int_0^T \left\|\Gamma^{\alpha_0}\left(2\partial_T \eta \, \partial_T u - 2\nabla_X \eta \cdot \nabla_X u + (\Box_g \eta)u\right)(S,\cdot)\right\| dS$$
$$\leq (\pi - T)^{2|\alpha_0|} \int_0^T (\pi - S)^{-2|\alpha_0|-1} \times$$
$$\sum_{|\alpha| \leq |\alpha_0|} \left(\|Z^\alpha u'(S,\cdot)\|_{L^2(R \leq (\pi-T)/2)} + \|Z^\alpha u(S,\cdot)\|_{L^6(R \leq (\pi-T)/2)}\right) dS$$

We have already shown that for each $S \leq T$ the summand is bounded by the right hand side of (5.1), and consequently the integral is bounded by the right hand side of (5.1). The final term in $G$ is similarly bounded.

To bound the last term on the right hand side of (5.10), we observe that

$$(\pi - T)^{2|\alpha_0|} \int_0^T \|[\Gamma^{\alpha_0}, \gamma^{jk} \Gamma_j \Gamma_k] v(S,\cdot)\|_2 \, dS$$
$$\leq C \int_0^T (\pi - S)^{-2} \sum_{|\alpha_1| + |\alpha_2| \leq |\alpha_0|} \left(\|(\Gamma^{\alpha_1} \gamma^{jk}) \Gamma^{\alpha_2} u'(S,\cdot)\|_2 + \|(\Gamma^{\alpha_1} \gamma^{jk}) \Gamma^{\alpha_2} u(S,\cdot)\|_6\right) dS \,,$$

which is contained in the right hand side of (5.1), completing the proof of Theorem 5.1.

**Proof of Lemma 5.2.** Formula (5.4) follows immediately from (2.3). We next recall that

$$\Box_g u = u_{TT} - u_{RR} - \frac{1}{\sin^2(R)} \left(u_{\phi_1 \phi_1} + u_{\phi_2 \phi_2}\right) - \frac{2\cos(R)}{\sin(R)} u_R - \frac{\cos(\phi_1)}{\sin^2(R) \sin(\phi_1)} u_{\phi_1}$$

Recalling (5.4), we obtain

$$\begin{aligned}
[\Box_g, \mathcal{X}] &\equiv [\partial_T^2, X] - [\Delta_{S^3}, X] \\
&= [\partial_T^2, \mathcal{X}] - \left[\partial_R^2 + \frac{2\cos R}{\sin R}\partial_R, (1+\cos R\,\cos T\,\partial_T)\right] \\
&\quad + \left[\partial_R^2 + \frac{2\cos R}{\sin R}\partial_R, \sin T\,\sin R\,\partial_R\right] \\
&\quad + \left[\frac{1}{\sin^2 R}(\partial_{\phi_1} + \partial_{\phi_2}) + \frac{\cos\phi_1}{\sin^2 R\sin\phi_1}\partial_{\phi_1}, \sin T\,\sin R\,\partial_R\right] \\
&= I + II + III + IV.
\end{aligned}$$

We first notice that

$$I = -2\cos R\,\sin T\,\partial_T^2 - \cos R\,\cos T\,\partial_T + \sin T\,\sin R\,\partial_R - 2\cos T\,\sin R\,\partial_R\,\partial_T,$$

and

$$II = 2\sin R\,\cos T\,\partial_R\,\partial_T + \cos R\,\cos T\,\partial_T + \frac{2\cos R}{\sin R}\sin R\,\cos T\,\partial_T.$$

Also,

$$IV = 2\cos R\,\sin T\,\frac{1}{\sin^2 R}(\partial_{\phi_1}^2 + \partial_{\phi_2}^2) + 2\cos R\,\sin T\,\frac{\cos\phi_1}{\sin^2 R\sin\phi_1}\partial_{\phi_1}.$$

Thus,

$$\begin{aligned}
I + II + IV = -2\cos R\,\sin T\,\Big(\partial_T^2 &- \frac{1}{\sin^2 R}(\partial_{\phi_1}^2 + \partial_{\phi_2}^2) - \frac{\cos\phi_1}{\sin^2 R\sin\phi_1}\partial_{\phi_1}\Big) \\
&+ \sin T\,\sin R\,\partial_R + 2\cos R\,\cos T\,\partial_T.
\end{aligned}$$

The remaining term $III$ equals

$$\begin{aligned}
&2\sin T\,\cos R\,\partial_R^2 - \sin T\,\sin R\,\partial_R + \Big(\frac{2\cos R}{\sin R}\sin T\,\cos R - \sin T\,\sin R\frac{-2}{\sin^2 R}\Big)\partial_R \\
&= 2\sin T\,\cos R\,\partial_R^2 - \sin T\,\sin R\,\partial_R + 2\frac{\sin T}{\sin R}(\cos^2 R + 1)\partial_R \\
&= 2\sin T\,\cos R\,\partial_R^2 - \sin T\,\sin R\,\partial_R + 2\sin T\,\cos R\,\frac{2\cos R}{\sin R}\partial_R + 2\sin T\,\sin R\,\partial_R.
\end{aligned}$$

If we combine the last two steps we obtain the equality

$$\begin{aligned}
[\Box_g, \mathcal{X}] = -2\cos R\,\sin T\,\Big(\partial_T^2 - \partial_R^2 &- \frac{1}{\sin^2 R}(\partial_{\phi_1}^2 + \partial_{\phi_2}^2) - \frac{\cos\phi_1}{\sin^2 R\sin\phi_1}\partial_{\phi_1} - \frac{2\cos R}{\sin R}\partial_R\Big) \\
&+ 2\cos R\,\cos T\,\partial_T + 2\sin T\,\sin R\,\partial_R,
\end{aligned}$$

as claimed. $\square$

## 6. Pointwise estimates.

To prove our global existence theorem for quasilinear equations we shall need pointwise estimates for solutions of the unperturbed Dirichlet-wave equation on $Y$.

**Theorem 6.1.** *Let $u \in C^\infty$ vanish on $\partial Y$, and let $(\Box_g + 1)u = F$. Then for every fixed $p > 1$ and $k = 0, 1, 2, \ldots$, there is a constant $C = C_{p,k}$ so that for $0 < T < \pi$*

$$\sum_{|\alpha|\leq k} |Z^\alpha u(T,X)| \leq C \sup_{0\leq S\leq T} \sum_{|\alpha|\leq k+1} \Big(\|Z^\alpha F(S,\cdot)\|_2 + (\pi - S)^{-2}\|Z^\alpha F(S,\cdot)\|_p\Big)$$
$$+ C \sum_{|\alpha|\leq k+2} \|Z^\alpha u(0,\cdot)\|_2. \tag{6.1}$$

We shall use separate arguments to establish (6.1) on the set $\mathcal{B}_1$ and on its complement. Recall that $\mathcal{B}_1$ is the pushforward of the set $|x| \leq 1$ from Minkowski space, and is essentially a set of the form $R \leq c(\pi - T)^2$. Recall also that $\partial Y \subset \mathcal{B}_{1/4}$. To handle $\mathcal{B}_1$ we shall make use of the following

**Proposition 6.2.** *Let $u$ be as above. Then for fixed $k = 0, 1, 2, \ldots$, there is a constant $C$ so that for $0 < T < \pi$*

$$(\pi - T)^{-1} \sum_{|\alpha|\leq k+1} \Big(\|Z^\alpha u'(T,\cdot)\|_{L^2(\mathcal{B}_2^T)} + \|Z^\alpha u(T,\cdot)\|_{L^6(\mathcal{B}_2^T)}\Big)$$
$$\leq C \sup_{0\leq S\leq T} \sum_{|\alpha|\leq k+1} \Big(\|Z^\alpha F(S,\cdot)\|_2 + (\pi - S)^{-2}\|Z^\alpha F(S,\cdot)\|_p\Big)$$
$$+ C \sum_{|\alpha|\leq k+1} \|Z^\alpha u'(0,\cdot)\|_2. \tag{6.2}$$

The inequality (6.2) shows that on a suitable neighborhood of $\partial Y$, one can improve upon Theorem 5.1 by one power of $(\pi - T)$. The arguments needed to do this are similar to those in our previous paper [10]; we shall postpone the proof until the end of this section.

To apply (6.2) to obtaining pointwise estimates, we shall make use of the following estimate, which follows from the standard Sobolev lemma for $\mathbb{R}^3$ and a scaling argument.

**Lemma 6.3.** *There is a constant $C$ so that, for $0 < T < \pi$ and $h \in C^\infty(Y_T)$,*

$$\|h\|_{L^\infty(\mathcal{B}_1^T)} \leq C(\pi - T) \sum_{|\alpha|\leq 2} \|\Gamma^\alpha h\|_{L^2(\mathcal{B}_2^T)} + C(\pi - T)^{-1}\|h\|_{L^6(\mathcal{B}_2^T)}. \tag{6.3}$$

Note that $h$ does not have to vanish on $\partial Y_T$. If we apply this estimate to $h = Z^\alpha u(T, X)$, for $|\alpha| \leq k$, then we conclude that

$$\sum_{|\alpha|\leq k} \|Z^\alpha u(T,\cdot)\|_{L^\infty(\mathcal{B}_1^T)}$$
$$\leq C(\pi - T) \sum_{|\alpha|\leq k} \Big(\|Z^\alpha \Gamma u'(T,\cdot)\|_{L^2(\mathcal{B}_2^T)} + \|Z^\alpha u'(T,\cdot)\|_{L^2(\mathcal{B}_2^T)} + \|Z^\alpha u(T,\cdot)\|_{L^2(\mathcal{B}_2^T)}\Big)$$
$$+ C(\pi - T)^{-1} \sum_{|\alpha|\leq k} \|Z^\alpha u(T,\cdot)\|_{L^6(\mathcal{B}_2^T)}. \tag{6.4}$$

Since $S^3 \setminus Y_T$ is star-shaped, and since $u(T, X)$ vanishes when $X \in \partial Y_T$, a simple calculus argument, using the fact that $R \leq C(\pi - T)^2$ if $(T, X) \in \mathcal{B}_2$, yields

$$(\pi - T)^{-1} \|u(T, \cdot)\|_{L^6(\mathcal{B}_2^T)} \leq C(\pi - T) \sum_{|\alpha|=1} \|\Gamma u(T, \cdot)\|_{L^6(\mathcal{B}_2^T)},$$

and since

$$(\pi - T)^{-1} \sum_{0 < |\alpha| \leq k} \|Z^\alpha u(T, \cdot)\|_{L^6(\mathcal{B}_2^T)} \leq C(\pi - T) \sum_{0 < |\alpha| \leq k-1} \|Z^\alpha u'(T, \cdot)\|_{L^6(\mathcal{B}_2^T)}$$

we conclude that the terms involving the $L^6$-norms in the right side of (6.4) are dominated by

$$(\pi - T) \sum_{|\alpha| \leq k} \|Z^\alpha u'(T, \cdot)\|_{L^6(\mathcal{B}_2^T)}.$$

Similar arguments give

$$(\pi - T) \sum_{|\alpha| \leq k} \|Z^\alpha u(T, \cdot)\|_{L^2(\mathcal{B}_2^T)} \leq (\pi - T)^2 \sum_{|\alpha \leq k} \|Z^\alpha u'(T, \cdot)\|_{L^2(\mathcal{B}_2^T)}.$$

Thus, (6.2) and (6.4) and Hölder's inequality imply that

$$\sum_{|\alpha| \leq k} \|Z^\alpha u(T, \cdot)\|_{L^\infty(\mathcal{B}_1^T)}$$
$$\leq C \sup_{0 \leq S \leq T} \sum_{|\alpha| \leq k+1} \|Z^\alpha F(S, \cdot)\|_2 + C \sum_{|\alpha| \leq k+1} \|Z^\alpha u'(0, \cdot)\|_2. \quad (6.5)$$

To prove the bounds for $(T, X) \notin \mathcal{B}_1$, we shall use the following estimate for the free wave equation.

**Proposition 6.4.** *Suppose that $u_f \in C^\infty([0, \pi) \times S^3)$ and that $(\Box_g + 1) u_f = F$. Then for $1 < p \leq 2$ and $T \in (\pi/2, \pi)$,*

$$|u_f(T, X)| \leq C \int_{\pi/2}^T (T - S)^{2-3/p} \big(\|F'(S, \cdot)\|_p + \|F(S, \cdot)\|_p\big) \, dS$$
$$+ C \sum_{|\alpha| \leq 2} \|Z^\alpha u_f(\pi/2, \cdot)\|_2. \quad (6.6)$$

**Proof.** When $F = 0$, the estimate holds by the energy inequality and Sobolev embedding. We will thus assume that the Cauchy data of $u_f$ vanishes when $T = \pi/2$. The proof is then a consequence of Duhamel's formula together with the following estimate, where $\Delta_g$ is the standard Laplacian on $S^3$

$$\left\| \frac{\sin(T-S)(1-\Delta_g)^{1/2}}{(1-\Delta_g)} \right\|_{L^p \to L^\infty} = O((T-S)^{2-3/p}), \quad 1 < p \leq 2.$$

which is valid for $|T - S| \leq \pi/2$. This estimate in turn is a consequence of the following dyadic estimate, where we take $\beta \in C_0^\infty((1/2, 2))$,

$$\left\| \beta(\sqrt{-\Delta_g}/\lambda) \frac{\sin(T-S)(1-\Delta_g)^{1/2}}{(1-\Delta_g)} \right\|_{L^p \to L^\infty}$$
$$\leq C \min\{(T-S)^{3-4/p} \lambda^{1-1/p}, (T-S)^{1-2/p} \lambda^{1/p-1}\},$$

which is valid for $1 \leq p \leq 2$. The dyadic estimate follows by interpolation from the endpoint $p = 1$, where the bounds are $O((T-S)^{-1})$ independent of $\lambda$, and the endpoint $p = 2$, where the bounds are $O(\varepsilon\lambda^{1/2})$ for $\lambda \leq \varepsilon^{-1}$ and $O(\lambda^{-1/2})$ for $\lambda \geq \varepsilon^{-1}$. □

We shall apply Proposition 6.4 to estimate $|u(T,X)|$ for $(T,X) \notin \mathcal{B}_1$. We need to show that, if $|\gamma| \leq k$, then

$$(\pi - T)^{2|\gamma|}\big|\Gamma^\gamma u(T,X)\big| \leq C \sup_{0 \leq S \leq T} \Big( \sum_{|\alpha| \leq k+1} \|Z^\alpha F(S,\cdot)\|_2 + (\pi - S)^{-2}\|Z^\alpha F(S,\cdot)\|_p \Big)$$
$$+ C \sum_{|\alpha| \leq k+1} \|Z^\alpha u'(0,\cdot)\|_2, \quad \text{if } (T,X) \notin \mathcal{B}_1.$$

We fix a cutoff function $\rho(T,X)$ satisfying $\rho = 0$ when $r(T,X) \leq 1/2$ and $\rho = 1$ when $r(T,X) \geq 1$ (recall that $r(T,X)$ denotes the Euclidean $r$ via the Penrose transformation), as well as the natural size esimates on its derivatives, $|Z^\alpha \rho(T,X)| \leq C_\alpha$.

We fix $\gamma$ with $|\gamma| \leq k$, and let $u_f = \rho \Gamma^\gamma u$. Since we are assuming that $\partial Y$ is contained in the set $r(T,R) \leq 1/4$, it follows that $u_f$ solves the free (no obstacle) wave equation

$$(\Box_g + 1)u_f = \rho \Gamma^\gamma F + 2\partial_T\rho\partial_T\Gamma^\gamma u - 2\nabla_X\rho \cdot \nabla_x\Gamma^\gamma u + (\Box_g\rho)\Gamma^\gamma u.$$

We next decompose

$$u_f = u_f^0 + u_f^1,$$

where $(\Box_g + 1)u_f^0 = \rho\Gamma^\gamma F$, with $u_f^0(0,\cdot) = u_f(0,\cdot)$, $\partial_T u_f^0(0,\cdot) = \partial_T u_f(0,\cdot)$. It then follows from (6.6) that

$$(\pi - T)^{2|\gamma|}|u_f^0(T,X)|$$
$$\leq C(\pi - T)^{2|\gamma|} \int_{\pi/2}^T (T-S)^{2-3/p}\Big(\|(\rho\Gamma^\gamma F)'(S,\cdot)\|_p + \|\rho\Gamma^\gamma F(S,\cdot)\|_p\Big)\, dS$$
$$+ C \sum_{|\alpha| \leq 2} \|Z^\alpha u_f(\pi/2, \cdot)\|_2$$
$$\leq C \sup_{0 \leq S \leq T} \Big( \sum_{|\alpha| \leq k+1} \|Z^\alpha F(S,\cdot)\|_2 + (\pi - S)^{-2}\|Z^\alpha F(S,\cdot)\|_p \Big)$$
$$+ C \sum_{|\alpha| \leq k+2} \|Z^\alpha u(0,\cdot)\|_2.$$

To finish the estimate for $(T,X) \notin \mathcal{B}_1$, we must show that $(\pi - T)^{2|\gamma|}|u_f^1(T,X)|$ can also be bounded by the right side of (6.1), where

$$(\Box_g + 1)u_f^1 = G = 2\,\partial_T\rho\,\partial_T\big(\Gamma^\gamma u\big) - 2\,\nabla_X\rho \cdot \nabla_X\big(\Gamma^\gamma u\big) + (\Box_g\rho)\big(\Gamma^\gamma u\big),$$

and $u_f^1$ has zero initial data. Note that $G$ is supported in $\mathcal{B}_1 \subset \{R \leq C\,(\pi-T)^2\}$.

We decompose $[0,\pi) = \cup_{j>0} I_j$ where $I_j$ are intervals $[a_j, b_j]$ with $a_{j+1} = b_j$ and $|I_j| \approx (\pi - b_j)^2$. We then fix a partition of unity $\chi_j$ on $[0,\pi)$ with $\chi_j$ supported in $I_{j-1} \cup I_j \cup I_{j+1}$ and $\chi_j^{(m)} \leq C_m |I_j|^{-m}$, and set

$$G_j(T,X) = \chi_j(T)\,G(T,X).$$

It follows that $G_j$ is supported in a cube of size $(\pi - b_j)^2$ centered at $T = b_j$, $R = 0$, and by Hölder's inequality we have the bound

$$\|G'_j(S,\cdot)\|_p + \|G_j(S,\cdot)\|_p \leq C\,(\pi - S)^{-5+6/p}\Big(\big\|(\Gamma^\gamma u)''\big\|_{L^2(\mathcal{B}_1^S)} + \big\|(\Gamma^\gamma u)'\big\|_{L^6(\mathcal{B}_1^S)}\Big)$$
$$+ C\,(\pi - S)^{-7+6/p}\|\Gamma^\gamma u\|_{L^6(\mathcal{B}_1^S)} \qquad (6.7)$$

Now let $\Lambda_j^+$ be the set of $(T, X)$ such that $T - R \in I_j$. By the sharp Huygen's principle for $\mathbb{R} \times S^3$, there is a constant $B$ independent of $j$ such that $u_f^1(T, X)$ depends only on $\sum_{|i-j| \leq B} G_i$.

Therefore (6.6) implies that, for $(T, X) \in \Lambda_j^+$, provided $I_j \subset (\pi/2, \pi)$, we have

$$(\pi - T)^{2|\gamma|}|u_f^1(T, X)|$$
$$\leq C_p(\pi - T)^{2|\gamma|} \sum_{|j-k| \leq B} \int_0^T (T - S)^{2-3/p}\Big(\|G'_j(S,\cdot)\|_p + \|G_j(S,\cdot)\|_p\Big)\, dS$$
$$\leq C_p \sum_{|j-k| \leq B} \sup_{0 \leq S \leq T} (\pi - S)^{6-6/p+2|\gamma|}\Big(\|G'_j(S,\cdot)\|_p + \|G_j(S,\cdot)\|_p\Big)$$

where we have used the fact that $G_j$ is supported in an interval of size $(\pi - b_j)^2 \approx (\pi - S)^2$. By (6.7), this is in turn bounded by

$$C \sup_{0 \leq S \leq T} (\pi - S)^{1+2|\gamma|}\Big(\|(\Gamma^\gamma u)''(S,\cdot)\|_{L^2(\mathcal{B}_1^S)}$$
$$+ \|(\Gamma^\gamma u)'(S,\cdot)\|_{L^6(\mathcal{B}_1^S)} + (\pi - S)^{-2}\|\Gamma^\gamma u(S,\cdot)\|_{L^6(\mathcal{B}_1^S)}\Big).$$

Since $(\pi - T)^{2|\gamma|}|u_f^1(T, X)| = |Z^\gamma u(T, X)|$ for $(T, X) \notin \mathcal{B}_1$, we can use Proposition 6.2 to conclude that

$$|Z^\gamma u(T, X)| \leq C \sup_{0 \leq S \leq T} \sum_{|\alpha| \leq k+1} \|Z^\alpha F(S,\cdot)\|_2,$$

provided that $(T, X) \notin \mathcal{B}_1$, and $|T - R| \leq \pi/2$.

For $(T, X) \notin \mathcal{B}_1$ and $|T - R| \geq \pi/2$, we modify the above procedure by using energy estimates to bound

$$|u_f^1(T, X)| \leq \sup_{0 \leq S \leq \pi/2}\Big(\|G'(S,\cdot)\|_2 + \|G(S,\cdot)\|_2\Big)$$
$$\leq \sup_{0 \leq S \leq \pi/2}\Big(\|\Gamma^\gamma u''(S,\cdot)\|_{L^2(\mathcal{B}_1^S)} + \|\Gamma^\gamma u'(S,\cdot)\|_{L^6(\mathcal{B}_1^S)} + \|\Gamma^\gamma u(S,\cdot)\|_{L^6(\mathcal{B}_1^S)}\Big)$$

which completes the proof of Theorem 6.1 for $(T, X) \notin \mathcal{B}_1$.

**Proof of Proposition 6.2.**

We begin by showing that Proposition 6.2 is a consequence of the following estimate,

$$(\pi - T)^{-1}\Big(\|\mathcal{X}^{k+1}u'(T,\cdot)\|_{L^2(\mathcal{B}_3^T)} + \|u(T,\cdot)\|_{L^6(\mathcal{B}_3^T)}\Big)$$
$$\leq C \sup_{0\leq S\leq T} \sum_{|\alpha|\leq k+1} \Big(\|Z^\alpha F(S,\cdot)\|_2 + (\pi-S)^{-2}\|Z^\alpha F(S,\cdot)\|_p\Big)$$
$$+ C \sum_{|\alpha|\leq k+1} \|Z^\alpha u'(0,\cdot)\|_2, \qquad (6.8)$$

where, as before,

$$\mathcal{X} = (1 + \cos T \cos R)\, \partial_T - \sin T \sin R\, \partial_R$$

is the pushforward $\mathcal{P}_*(\partial/\partial t)$ of the time derivative from Minkowski space. We will use the fact that, absent the term $(\pi-T)^{-1}$ on the left hand side, the estimate (6.2) would be immediate from Theorem 5.1. Consequently, error terms in our calculations that involve an extra power of $(\pi - T)$ can be dominated by the right side of (6.8) using Theorem 5.1. Terms involving the commutator $[\Box, \mathcal{X}]$ fall in this category. We also make use of the fact that on $\mathcal{B}_3$

$$\mathcal{X} = \frac{1}{2}(\pi - T)^2\, \partial_T + O\big((\pi-T)^3\big)\,\Gamma,$$

and consequently, using the equation $\partial_T^2 u = \Delta_g u + F$,

$$(\pi - T)\|\Delta_g \mathcal{X}^k u\|_{L^2(\mathcal{B}_3)} \leq C (\pi - T)^{-1} \|\mathcal{X}^{k+1}u'\|_{L^2(\mathcal{B}_3)} + \cdots$$

where $\cdots$ indicates terms that can be dominated by the right side of (6.8) using Theorem 5.1. By the remark in the proof of Proposition 4.3, we conclude that

$$(\pi - T)^{-1} \sum_{|\alpha|\leq 1}\Big(\|Z^\alpha \mathcal{X}^k u'\|_{L^2(\mathcal{B}_{3-\varepsilon})} + \|Z^\alpha \mathcal{X}^k u\|_{L^6(\mathcal{B}_{3-\varepsilon})}\Big)$$
$$\leq (\pi - T)^{-1}\Big(\|\mathcal{X}^{k+1}u'\|_{L^2(\mathcal{B}_3)} + \|u\|_{L^6(\mathcal{B}_3)}\Big) + \cdots.$$

Repeating this procedure shows that Proposition 6.2 follows from (6.8).

We now show that (6.8) is a consequence of the following lemma, which states that better estimates hold if the data and forcing term vanish outside of $\mathcal{B}_8$.

**Lemma 6.5.** *Let $u$ be as above, Assume further that if*

$$(\Box_g + 1)u = F$$

*then $F(T,X) = 0$ in $\mathcal{B}_8^c$ Suppose also that $0 = \partial_T u(0,X) = u(0,X)$ when $(0,X) \notin \mathcal{B}_8$. Then, for each $k = 0, 1, 2, \ldots$, there is a constant $C$ so that for $0 < T < \pi$*

$$(\pi - T)^{-1}\Big(\|\mathcal{X}^{k+1}u'(T,\cdot)\|_{L^2(\mathcal{B}_8^T)} + \|u(T,\cdot)\|_{L^6(\mathcal{B}_8^T)}\Big)$$
$$\leq C(\pi - T) \sup_{0\leq S\leq T} \Big(\|\mathcal{X}^{k+1}F(S,\cdot)\|_2 + \|F(S,\cdot)\|_2 + (\pi - T)\sum_{|\alpha|\leq k}\|Z^\alpha F(S,\cdot)\|_2\Big)$$
$$+ C \sum_{|\alpha|\leq k+1} \|Z^\alpha u'(0,\cdot)\|_2. \qquad (6.9)$$

Before proving Lemma 6.5, we show that it implies (6.8) as a consequence. The argument uses techniques from [25].

We fix $\eta \in C^\infty$ so that $\eta = 1$ in $\mathcal{B}_4$ and $\eta = 0$ in $\mathcal{B}_8^c$, and such that $|Z^\alpha \eta| \leq C_\alpha$ for each $\alpha$. By taking $\eta(T, X) = \beta(r(T, X))$ for appropriate $\beta \in C^\infty(\mathbb{R})$, we may assume $\mathcal{X}\eta = 0$. We then split

$$u = v + w,$$

where

$$(\Box_g + 1)v = \eta F, \qquad (\Box_g + 1)w = (1 - \eta)F,$$

and

$$v(0, X) = (\eta u)(0, X), \qquad \partial_T v(0, X) = \partial_T(\eta u)(0, X).$$

The estimate (6.9) then yields the following estimate for $v$ that is even stronger than (6.8),

$$(\pi - T)^{-1}\Big(\|\mathcal{X}^{k+1} v'(T, \cdot)\|_{L^2(\mathcal{B}_3^T)} + \|v(T, \cdot)\|_{L^6(\mathcal{B}_3^T)}\Big)$$
$$\leq C(\pi - T) \sup_{0 \leq S \leq T} \Big(\|\mathcal{X}^{k+1} F(S, \cdot)\|_2 + \|F(S, \cdot)\|_2 + (\pi - T) \sum_{|\alpha| \leq k} \|Z^\alpha F(S, \cdot)\|_2\Big)$$
$$+ C \sum_{|\alpha| \leq k+1} \|Z^\alpha u'(0, \cdot)\|_2. \quad (6.10)$$

In the last step, we use that fact that $u$ satisfies Dirichlet conditions, which shows that $\sum_{|\alpha| \leq k+1} \|Z^\alpha v'(0, \cdot)\|_2$ is dominated by the last term in (6.10).

To handle the term $w$, we fix $\rho \in C^\infty([0, \pi) \times S^3)$ satisfying $\rho = 1$ in $\mathcal{B}_3$, $\rho = 0$ in $\mathcal{B}_4^c$, and

$$|Z^\alpha \rho| \leq C_\alpha \text{ for each } \alpha, \quad \text{and} \quad \partial_T^k \rho = O((\pi - T)^{-k}), \; k = 0, 1, 2, \ldots \quad (6.11)$$

This can be achieved by setting $\rho(T, X) = \beta(r(T, X))$, for appropriate $\beta \in C^\infty(\mathbb{R})$. We then write

$$w = w_f + w_r,$$

where $w_f$ solves the free (no obstacle) wave equation on $S^3 \times [0, \pi)$ with the same data as $w$,

$$(\Box_g + 1)w_f = (1 - \eta)F,$$
$$w_f(0, X) = \big((1-\eta)u\big)(0, X), \quad \partial_T w_f(0, X) = \partial_T\big((1-\eta)u\big)(0, X).$$

(Recall that $\eta$ vanishes near $\partial Y$.)

For $(T, X) \in \mathcal{B}_3$, the function $w$ agrees with the function $w_0$ defined by

$$w_0 = \rho w_f + w_r.$$

Note that $w_0$ solves the Dirichlet-wave equation

$$(\Box_g + 1)w_0 = G = 2(\partial_T \rho)(\partial_T w_f) - 2(\nabla_X \rho) \cdot (\nabla_X w_f) + (\Box_g \rho)w_f,$$

since $\eta = 1$ on the support of $\rho$. Applying (6.9), we thus obtain

$$(\pi - T)^{-1}\Big(\|\mathcal{X}^{k+1}w'(T,\cdot)\|_{L^2(\mathcal{B}_3^T)} + \|w(T,\cdot)\|_{L^6(\mathcal{B}_3^T)}\Big)$$
$$\leq C(\pi - T) \sup_{0 \leq S \leq T}\Big(\|\mathcal{X}^{k+1}G(S,\cdot)\|_2 + \|G(S,\cdot)\|_2 + (\pi - T)\sum_{|\alpha| \leq k}\|Z^\alpha G(S,\cdot)\|_2\Big)$$
$$+ C\sum_{|\alpha| \leq k+1}\|Z^\alpha u'(0,\cdot)\|_2.$$

Since $G(S, X) = 0$ when $R \geq C(\pi - S)^2$, applying Hölder's inequality and (6.11) yields

$$(\pi - S)\|\mathcal{X}^{k+1}G(S,\cdot)\|_2 + (\pi - S)^2\sum_{|\alpha| \leq k}\|Z^\alpha G(S,\cdot)\|_2$$
$$\leq C\sum_{1 \leq |\alpha| \leq k+1}(\pi - S)^{2|\alpha|-1}\Big(\|\Gamma^\alpha w'_f(S,\cdot)\|_2 + \|\Gamma^\alpha w_f(S,\cdot)\|_6\Big)$$
$$+ C\Big(\|w'_f(S,\cdot)\|_2 + \|w_f(S,\cdot)\|_6\Big). \tag{6.13}$$

Recall that $(\Box_g + 1)w_f = (1 - \eta)F$, and $|Z^\alpha \eta| \leq C_\alpha$. Energy estimates for the free (no obstacle) wave equation on $S^3 \times [0, \pi)$, together with the fact that $\Gamma$ commutes with $\Box_g$, show that the right side of (6.13) is dominated by

$$\int_0^S \|F(s,\cdot)\|_2\, ds + \sum_{1 \leq |\alpha| \leq k+1}(\pi - S)^{2|\alpha|-1}\int_0^S \|\Gamma^\alpha F(s,\cdot)\|_2\, ds + \sum_{|\alpha| \leq k+1}\|\Gamma^\alpha w'_f(0,\cdot)\|_2$$
$$\leq C \sup_{0 \leq S \leq T}\sum_{|\alpha| \leq k+1}\|Z^\alpha F(S,\cdot)\|_2 + \sum_{|\alpha| \leq k+1}\|Z^\alpha u'(0,\cdot)\|_2,$$

and the last terms are contained in the right side of (6.8).

To finish the desired estimates for $w$, it remains to show that we can bound the quantity $(\pi - S)\|G(S,\cdot)\|_2$ by the right hand side of (6.8). To do this, we use Hölder's inequality and the fact that $G = 0$ for $R \geq C(\pi - T)^2$ to bound

$$(\pi - S)\|G(S,\cdot)\|_2 \leq C\|w_f(S,\cdot)\|_\infty + C(\pi - S)\|w'_f(S,\cdot)\|_6.$$

The last term is contained in the right hand side of (6.13). On the other hand, by Proposition 6.4, we may bound

$$\|w_f(S,\cdot)\|_{L^\infty} \leq C \sup_{0 \leq S \leq T}(\pi - S)^{-2}\sum_{|\alpha| \leq 1}\|Z^\alpha F(S,\cdot)\|_p + C\sum_{|\alpha| \leq 1}\|Z^\alpha u'(0,\cdot)\|_2,$$

and the right hand side here is also contained in (6.13). This completes the reduction of Proposition 6.2 to Lemma 6.5. $\square$

To prove Lemma 6.5, we will pull things back to Minkowski space in order to exploit the energy decay estimates of Morawetz, Lax and Phillips.

To begin, we note that $\mathcal{P}^{-1}(\mathcal{B}_{r_1}^T) = \Pi_{r_1}^T$, where

$$\Pi_{r_1}^T = \{(t, x) \in \mathbb{R}_+ \times \mathbb{R}^3 \setminus \mathcal{K} : |x| \leq r_1,\ (t + \lambda)^2 = |x|^2 + 1 + \lambda^2,\ \lambda = \cot T\}. \tag{6.14}$$

The manifolds $\Pi^T_{r_1}$ form a uniform family of timelike hypersurfaces as $\lambda$ varies over $(-\infty, \infty)$, as can be seen by expressing $\Pi^T_{r_1}$ in the form

$$t = \tan(T/2) + \sqrt{1 + \lambda^2 + |x|^2} - \sqrt{1 + \lambda^2}\,.$$

In particular, it follows from this that

$$(t, x) \in \Pi^T_{r_1} \implies t \in [\tan(T/2), \tan(T/2) + r_1]\,. \tag{6.15}$$

Let $\Pi^T_{r_1}$ be endowed with the induced Lebesgue measure. We will use the following consequence of the Morawetz, Lax and Phillips energy decay estimates for star-shaped obstacles $\mathcal{K}$.

**Lemma 6.6.** *Suppose that $\tilde{u} \in C^\infty(\mathbb{R}_+ \times \mathbb{R}^3 \setminus \mathcal{K})$ vanishes for $x \in \partial \mathcal{K}$, and $(\partial_t^2 - \Delta)\tilde{u} = \tilde{F}$. Suppose also that $\tilde{F}$, $\tilde{u}(0, \cdot)$, and $\partial_t \tilde{u}(0, \cdot)$ vanish for $|x| > r_1$. Then there are constants $c > 0$ and $C < \infty$, depending on $r_1$, so that for all $T$,*

$$\|\partial_t^k \tilde{u}'\|_{L^2(\Pi^T_{r_1})} + \|\tilde{u}\|_{L^6(\Pi^T_{r_1})}$$
$$\leq C \sup_{S \leq T} \Big( \|\partial_t^k \tilde{F}\|_{L^2(\Pi^S_{r_1})} + \|\tilde{F}\|_{L^2(\Pi^S_{r_1})} \Big) + C\, e^{-c/(\pi - T)} \|\partial_t^k \tilde{u}'(0, \cdot)\|_2\,.$$

**Proof of Lemma 6.6.** It suffices to consider the case $k = 0$, since $\partial_t$ commutes with $\Box$ and preserves the Dirichlet conditions. We also use the fact that

$$\|\tilde{u}\|_{L^6(\Pi^T_{r_1})} \leq C \, \|\tilde{u}'\|_{L^2(\Pi^T_{2r_1})}\,.$$

We now use the energy decay estimate of Morawetz, Lax and Phillips, which says that for star-shaped obstacles $\mathcal{K}$, for given fixed $r_0$ there are constants $c > 0$ and $C < \infty$ so that

$$\|\tilde{u}'(t_0, \cdot)\|_{L^2(|x| \leq r_0)} \leq C \int_0^{t_0} e^{-c(t_0 - s)} \|\tilde{F}(s, \cdot)\|_2 \, ds + C\, e^{-c t_0} \|\tilde{u}'(0, \cdot)\|_2\,.$$

Also, by (6.15) and energy estimates, we have

$$\|\tilde{u}'\|_{L^2(\Pi^T_{2r_1})} \leq C \, \|\tilde{u}'(\tan(T/2), \cdot)\|_{L^2(|x| \leq 4r_1)} + C \sup_{S \leq T} \|\tilde{F}\|_{L^2(\Pi^S_{r_1})}\,.$$

Taking $t_0 = \tan(T/2) \approx 1/(\pi - T)$, the result now follows from the simple estimate

$$\int_0^{t_0} e^{-c(t_0 - s)} \|\tilde{F}(s, \cdot)\|_2 \, ds \leq C \sup_{S \leq T} \|\tilde{F}\|_{L^2(\Pi^S_{r_1})}\,. \quad \Box$$

**Proof of Lemma 6.5.** We will identify points $(t, x)$ in Minkowski space with points $(T, X)$ in the Einstein diamond via the Penrose transform $\mathcal{P}$. Let $\Box_g$ denote the wave operator on $S^3 \times \mathbb{R}$, and $\Box$ the wave operator on Minkowski space $\mathbb{R}^4$. The map $\mathcal{P}$ is conformal relative to the respective Lorentzian metrics, and if we let

$$\tilde{u} = \Omega u\,, \qquad \tilde{F} = \Omega^3 F\,,$$

then

$$(\Box_g + 1)u = F \iff \Box \tilde{u} = \tilde{F}\,.$$

As a map from $\{|x| \leq r_1\}$ to $\mathcal{B}_{r_1}$, the map $\mathcal{P}$ is also essentially conformal in the respective Riemannian metrics, in the sense that if $X_j$ are projective coordinates on $S^3$ near the north pole, then for $|x| \leq r_1$,

$$dX_j = \frac{2}{1+t^2}\,dx_j + O(t^{-3})(dt, dx)\,, \qquad dT = \frac{2}{1+t^2}\,dt + O(t^{-3})(dt, dx)\,.$$

Also, $(1+t^2)^{-1} \approx (\pi - T)^2 \approx \Omega$ on $\mathcal{B}_{r_1}$.

Consequently, if $d\sigma_T$ denotes the measure induced on $\Pi_{r_1}^T$ by $dx\,dt$, and $dX$ denotes the volume form on $S^3$, then $dX \approx (\pi - T)^6\,d\sigma_T$ for $|x| < r_1$. Together with the fact that $\nabla_{T,X}\Omega = O(\pi - T)$, this implies

$$(\pi - T)^{-1}\|\mathcal{X}^{k+1}u'(T,\cdot)\|_{L^2(\mathcal{B}_{r_1}^T)} \leq C(\pi - T)^{-2}\left\|\partial_t^{k+1}\nabla_{t,x}\tilde{u}\right\|_{L^2(\Pi_{r_1}^T)}$$
$$+ C(\pi - T)^{-1}\sum_{0\leq j\leq k}\|\partial_t^j\nabla_{t,x}\tilde{u}\|_{L^2(\Pi_{r_1}^T)} + C\|\tilde{u}\|_{L^2(\Pi_{r_1}^T)}\,.$$

Since $\tilde{u}$ vanishes on $\partial\mathcal{K}$, the last term is dominated by the $L^2$ norm of $\nabla_{t,x}\tilde{u}$ over the same set. Also,

$$(\pi - T)^{-1}\|u(T,\cdot)\|_{L^6(\mathcal{B}_{r_1}^T)} \leq C(\pi - T)^{-2}\|\tilde{u}\|_{L^6(\Pi_{r_1}^T)}\,.$$

By Lemma 6.6, we conclude that

$$(\pi - T)^{-1}\Big(\|\mathcal{X}^{k+1}u'(T,\cdot)\|_{L^2(\mathcal{B}_{r_1}^T)} + \|u(T,\cdot)\|_{L^6(\mathcal{B}_{r_1}^T)}\Big)$$
$$\leq C\,\|\partial_t^{k+1}\tilde{u}'(0,\cdot)\|_2 + C\,(\pi-T)^{-2}\sup_{S\leq T}\Big(\|\partial_t^{k+1}\tilde{F}\|_{L^2(\Pi_{r_1}^S)} + \|\tilde{F}\|_{L^2(\Pi_{r_1}^S)}\Big)$$
$$+ C\,(\pi-T)^{-1}\sum_{j\leq k}\sup_{S\leq T}\Big(\|\partial_t^j\tilde{F}\|_{L^2(\Pi_{r_1}^S)} + \|\tilde{F}\|_{L^2(\Pi_{r_1}^S)}\Big)\,.$$

Since $\tilde{F} = \Omega^3 F \approx (\pi-T)^6 F$, this is in turn dominated by the right hand side of (6.9). □

## 7. Iteration Argument.

The purpose of this section is to show that we can solve certain Dirichlet-wave equations of the form

$$\begin{cases}(\Box_g + 1)u^I = \mathcal{F}^I(T, X; u, du, d^2u)\,, & I = 1, \ldots, N \\ u|_{\partial Y} = 0 \\ u|_{T=0} = f\,, \quad \partial_T u|_{T=0} = g\,,\end{cases} \qquad (7.1)$$

provided that the data is small and satisfies the appropriate compatibility conditions.

Regarding the nonlinear term, we shall assume that $\mathcal{F}$ satisfies (2.16)-(2.21). To apply our estimates, we shall also assume that $\mathcal{F}$ vanishes when $R \geq 2(\pi - T)$, that is, if $\gamma^{I,jk}$ and $\mathcal{G}$ are as in (2.14), then

$$\gamma^{I,jk} = 0 \quad \text{and} \quad \mathcal{G} = 0 \text{ if } R \geq 2(\pi - T)\,. \qquad (7.2)$$

This will not affect the existence results for Minkowski space since we may multiply the nonlinear term $\mathcal{F}$ in Proposition 2.4 by a cutoff that equals one on the image $R \leq (\pi - T)$ of Minkowski space under the Penrose transform.

The compatibility conditions for (7.1) are the pushforwards of the conditions in Definition 9.2. Specifically, if $(\tilde{f}, \tilde{g})$ denote the pullbacks of the data to $\mathbb{R}^3 \backslash \mathcal{K}$ given by (2.12), then we shall say that the data $(f, g)$ satisfies the compatibility condition of order $k$ for (7.1) if the Minkowski data $(\tilde{f}, \tilde{g})$ satisfies the compatibility conditions for (1.1) with the nonlinear term $F$ there given by (2.15).

We now state the existence result in $Y$ which will be used to prove our main result, Theorem 1.1.

**Theorem 7.1.** *Assume that $\mathcal{F}$ satisfies (2.16)–(2.18), as well as (7.2). Assume further that the Cauchy data $(f, g)$ is in $H_D^9(S^3 \backslash \mathcal{P}_0(\mathcal{K})) \times H_D^8(S^3 \backslash \mathcal{P}_0(\mathcal{K}))$ and that $(f, g)$ satisfies the compatibility condition of order $8$. Then there exists $\delta_0 > 0$, so that if*

$$\|f\|_{H^9(S^3 \backslash \mathcal{P}_0(\mathcal{K}))} + \|g\|_{H^8(S^3 \backslash \mathcal{P}_0(\mathcal{K}))} \leq \delta_0, \tag{7.3}$$

*then (7.1) has a solution in $Y$ verifying*

$$\sup_{0 \leq T < \pi} \sum_{|\alpha| \leq 8} \left( \|Z^\alpha u'(T, \cdot)\|_2 + \|Z^\alpha u(T, \cdot)\|_6 \right)$$
$$+ \sup_{0 \leq T < \pi} (\pi - T)^\sigma \sum_{|\alpha| \leq 5} \|Z^\alpha u(T, \cdot)\|_\infty < \infty \tag{7.4}$$

*for all $\sigma > 0$.*

Before turning to the proof of Theorem 7.1, we state a few simple consequences of the assumptions (2.16)-(2.18). To begin, assuming that (2.19) holds, simple bookkeeping and (2.18) implies that for a given $\alpha$

$$\begin{aligned}|Z^\alpha \mathcal{G}| &\leq C \sum_{|\gamma| \leq |\alpha|} |Z^\gamma u'| \cdot (\pi - T)^2 \sum_{|\gamma| \leq |\alpha|/2} |Z^\gamma u'| \\&+ C \sum_{|\gamma| \leq |\alpha|} |Z^\gamma u'| \cdot (\pi - T) \sum_{|\gamma| \leq (1+|\alpha|)/2} |Z^\gamma u| \\&+ C(\pi - T)^2 \sum_{|\gamma| \leq |\alpha|} |Z^\gamma u'| \Big( \sum_{|\gamma| \leq |\alpha|/2} |Z^\gamma u| \Big)^2 + C|u|^3 + C|u|^2\,. \end{aligned} \tag{7.5}$$

Similarly, if $\alpha$, $I$, $j$, and $k$ are fixed, then (2.17) implies

$$|Z^\alpha(\gamma^{I,jk}(T, X; v, v')\Gamma_j \Gamma_k u^I)| \leq C \sum_{|\gamma| \leq |\alpha|+1} |Z^\gamma u'| \sum_{|\gamma| \leq 1+|\alpha|/2} |Z^\gamma v|$$
$$+ C \sum_{|\gamma| \leq |\alpha|+1} |Z^\gamma v'| \sum_{|\gamma| \leq 2+|\alpha|/2} |Z^\gamma u|\,, \tag{7.6}$$

and, if $N = 0, 2, \ldots$ is an even integer, then (2.17) implies

$$\sum_{\substack{|\alpha_1|+|\alpha_2| \leq N+1 \\ |\alpha_1|, |\alpha_2| \leq N}} \left| Z^{\alpha_1}\big(\gamma^{I,jk}(T, X; v, v')\big) Z^{\alpha_2} u' \right| \leq C(\pi - T)^2 \sum_{|\alpha| \leq 1+N/2} |Z^\alpha v| \sum_{|\alpha| \leq N} |Z^\alpha u'|$$
$$+ C(\pi - T)^2 \sum_{|\alpha| \leq N} |Z^\alpha v'| \sum_{|\alpha| \leq 1+N/2} |Z^\alpha u|\,. \tag{7.7}$$

For the first step in the proof of Theorem 7.1, we simplify the task at hand by reducing (7.1) to an equivalent inhomogeneous equation with zero Cauchy data. By making this reduction we shall not have to worry about the role of the compatibility conditions in the iteration argument to follow.

To make this reduction we shall use the fact that there is a local solution to (7.1). Specifically, given data $(f, g)$ as above satisfying the compatibility conditions, there exists a time $0 < T_0 < \pi$ and a solution $u$ of (7.1) verifying

$$\sum_{|\alpha| \leq 9} \sup_{0 \leq T \leq T_0} \|\partial^\alpha u(T, \cdot)\|_{L^2(Y_T)} \leq C\,\delta_0\,, \tag{7.8}$$

if (7.3) holds. The existence of $u$ follows from Theorem 9.4. To see this, we pull back the data and the equation to Minkowski space, and use Theorem 9.4 to show existence of $u$ on a neighborhood of the boundary. Away from the boundary, the existence of $u$ follows by applying Theorem 9.4 in the Einstein diamond.

To use this, let us fix a cutoff $\eta \in C^\infty(\mathbb{R})$ satisfying

$$\eta(T) = 1 \text{ if } T \leq T_0/2\,, \text{ and } \eta(T) = 0,\ T \geq T_0\,. \tag{7.9}$$

We then set

$$u_0 = \eta u\,, \tag{7.10}$$

and note that

$$(\Box_g + 1)u_0 = \eta \mathcal{F}(T, X; u, du, d^2 u) + [\Box_g, \eta] u\,. \tag{7.11}$$

Therefore, if we put

$$w = u - u_0\,,$$

then $u$ will solve (7.1) if and only if $w$ solves

$$\begin{cases} (\Box_g + 1)w = (1 - \eta)\mathcal{F}(T, X; u_0 + w, d(u_0 + w), d^2(u_0 + w)) - [\Box_g, \eta](u_0 + w) \\ w|_{\partial Y} = 0 \\ w(0, X) = \partial_T w(0, X) = 0\,, \quad X \in S^3 \backslash \mathcal{P}_0(\mathcal{K})\,. \end{cases} \tag{7.12}$$

Note that the compatibility conditions are satisfied in this case since the data vanishes and since the forcing term in the equation vanishes on $[0, T_0/2]$.

We shall solve (7.12) by iteration. We begin by fixing $\sigma = 1/4$, and let

$$m(T, w) = \sum_{|\alpha| \leq 8} \Big( \|Z^\alpha w'(T, \cdot)\|_2 + \|Z^\alpha w(T, \cdot)\|_6 \Big) + (\pi - T)^\sigma \sum_{|\alpha| \leq 5} \|Z^\alpha w(T, \cdot)\|_\infty\,. \tag{7.13}$$

We will show that $m(T, w)$ will be small for each iterate $w = w_k$ if $\varepsilon = \sup_T m(T, u_0)$ is small. Indeed, we shall show that, for such $\varepsilon$,

$$m(T, w_k) \leq C_0\,\varepsilon\,,$$

where $C_0$ is a fixed constant. We will then show that the decay estimate (7.4) holds for all $\sigma > 0$ provided it holds for $\sigma = 1/4$. The estimate that allows us to carry out the iteration is the following.

**Lemma 7.2.** *There exists constants $C_0$, $\varepsilon_0 > 0$, so that if $\varepsilon \leq \varepsilon_0$, and*

$$\sup_{0 \leq T < \pi} m(T, v) \leq (C_0 + 1)\varepsilon, \tag{7.14}$$

$$\sup_{0 \leq T < \pi} m(T, u_0) \leq \varepsilon, \tag{7.15}$$

*then the solution $w$ to the equation*

$$\begin{cases} (\Box_g + 1)w^I = (1 - \eta)\left(\sum_{j,k} \gamma^{I,jk}(T, X; v, v')\Gamma_j\Gamma_k w^I + \mathcal{G}^I(T, X; v, v')\right) \\ \qquad\qquad\qquad\qquad\qquad\qquad -[\Box_g, \eta](u_0 + w)^I, \quad I = 1, \ldots, N \\ w|_{\partial Y} = 0 \\ w(0, X) = \partial_T w(0, X) = 0, \quad X \in S^3 \setminus \mathcal{P}_0(\mathcal{K}). \end{cases} \tag{7.16}$$

*satisfies*

$$\sup_{0 \leq T < \pi} m(T, w) \leq C_0 \varepsilon. \tag{7.17}$$

We will prove the lemma in the case $\varepsilon = \varepsilon_0$, under the assumption $\varepsilon_0$ is sufficiently small. In the various estimates below, we use $C$ to denote a constant that does not depend on $u_0$ or $v$, assuming just that $\varepsilon_0$ and $(C_0 + 1)\varepsilon_0$ are sufficiently small, which we will be able to arrange.

To apply Theorem 5.1 we note that, by (7.14) and (2.17), the conditions (3.14) and (3.15) are satisfied with $\delta$ small if $(C_0 + 1)\varepsilon_0$ is small. As a result, Theorem 5.1 and the remarks following it yield

$$\sum_{|\alpha| \leq 8} \left(\|Z^\alpha w'(T, \cdot)\|_2 + \|Z^\alpha w(T, \cdot)\|_6\right)$$

$$\leq C \sup_{0 \leq S \leq T} (\pi - S)^2 \sum_{|\alpha| \leq 7} \left(\|Z^\alpha \mathcal{G}(S, \cdot)\|_2 + \|Z^\alpha[\eta, \Box_g]w(S, \cdot)\|_2 + \|Z^\alpha[\eta, \Box_g]u_0(S, \cdot)\|_2\right)$$

$$+ C \sum_{I,j,k} \left(\sup_{0 \leq S \leq T} \sum_{\substack{|\alpha_1| + |\alpha_2| \leq 8 \\ |\alpha_1|, |\alpha_2| \leq 7}} \|(Z^{\alpha_1}\gamma^{I,jk})Z^{\alpha_2}w'(S, \cdot)\|_2 \right.$$

$$+ \int_0^T (\pi - S)^{-2} \sum_{\substack{|\alpha_1| + |\alpha_2| \leq 9 \\ |\alpha_1|, |\alpha_2| \leq 8}} \left(\|(Z^{\alpha_1}\gamma^{I,jk})Z^{\alpha_2}w'(S, \cdot)\|_2 + \|(Z^{\alpha_1}\gamma^{I,jk})Z^{\alpha_2}w(S, \cdot)\|_6\right) dS\right)$$

$$+ C \int_0^T \sum_{|\alpha| \leq 8} \|Z^\alpha \mathcal{G}(S, \cdot)\|_2 \, dS + C \int_0^T \sum_{|\alpha| \leq 8} \|[Z^\alpha \Box_g, \eta]w(S, \cdot)\|_2 \, dS$$

$$+ C \int_0^T \sum_{|\alpha| \leq 8} \|[Z^\alpha \Box_g, \eta]u_0(S, \cdot)\|_2 \, dS.$$

From the fact that $[\Box_g, \eta]$ is supported near $T = 0$, and the fact that $w$ vanishes at $0$, it is easy to see that

$$\sup_{0 \leq S \leq T} \sum_{|\alpha| \leq 7} \|Z^\alpha [\Box_g, \eta] w(S, \cdot)\|_2 + \int_0^T \sum_{|\alpha| \leq 8} \|Z^\alpha [\Box_g, \eta] w(S, \cdot)\|_2 \, dS$$

$$\leq C \int_0^T \sum_{|\alpha| \leq 8} \|Z^\alpha w'(S, \cdot)\|_2 \, dS.$$

Similarly, by (7.15) one obtains

$$\sum_{|\alpha| \leq 7} \|Z^\alpha [\Box_g, \eta] u_0(T, \cdot)\|_2 + \int_0^T \sum_{|\alpha| \leq 8} \|Z^\alpha [\Box_g, \eta] u_0(S, \cdot)\|_2 \, dS \leq C \varepsilon_0.$$

Thus, if

$$I + II + III + IV = C \int_0^T \sum_{|\alpha| \leq 8} \|Z^\alpha \mathcal{G}(S, \cdot)\|_2 dS + C \sup_{0 \leq S \leq T} (\pi - S)^2 \sum_{|\alpha| \leq 7} \|Z^\alpha \mathcal{G}(S, \cdot)\|_2$$

$$+ C \sum_{I,j,k} \Bigg( \sup_{0 \leq S \leq T} \sum_{\substack{|\alpha_1| + |\alpha_2| \leq 8 \\ |\alpha_1|, |\alpha_2| \leq 7}} \|(Z^{\alpha_1} \gamma^{I,jk}) Z^{\alpha_2} w'(S, \cdot)\|_2$$

$$+ \int_0^T (\pi - S)^{-2} \sum_{\substack{|\alpha_1| + |\alpha_2| \leq 9 \\ |\alpha_1|, |\alpha_2| \leq 8}} \Big( \|(Z^{\alpha_1} \gamma^{I,jk}) Z^{\alpha_2} w'(S, \cdot)\|_2 + \|(Z^{\alpha_1} \gamma^{I,jk}) Z^{\alpha_2} w(S, \cdot)\|_6 \Big) dS \Bigg)$$

then we have

$$\sum_{|\alpha| \leq 8} \Big( \|Z^\alpha w'(T, \cdot)\|_2 + \|Z^\alpha w(T, \cdot)\|_6 \Big)$$

$$\leq C \varepsilon_0 + C \int_0^T \sum_{|\alpha| \leq 8} \|Z^\alpha w'(S, \cdot)\|_2 \, dS + I + II + III + IV.$$

Using (7.5) and (7.14) we get

$$I + II \leq C \sup_{0 \leq S \leq T} (\pi - S)^{2-\sigma} (m(S,v))^2 + C \int_0^T (\pi - S)^{-\sigma} (m(S,v))^2 \, dS \leq C (C_0 + 1)^2 \varepsilon_0^2,$$

since we are assuming that $\sigma < 1/2$. Using (7.7) (with $N = 8$) and (7.14), we also obtain

$$III \leq C (C_0 + 1) \varepsilon_0 \sup_{0 \leq S \leq T} (\pi - S)^{2-\sigma} m(S, w),$$

while (7.7) with $N = 8$ also yields

$$\int_0^T (\pi - S)^{-2} \sum_{\substack{|\alpha_1| + |\alpha_2| \leq 9 \\ |\alpha_1|, |\alpha_2| \leq 8}} \|(Z^{\alpha_1} \gamma^{I,jk}) Z^{\alpha_2} w'(S, \cdot)\|_2 \, dS$$

$$\leq C (C_0 + 1) \varepsilon_0 \int_0^T (\pi - S)^{-\sigma} m(S, w) \, dS.$$

Since similar arguments give control of the $L^6$-norms, we conclude that

$$IV \leq C\,(C_0+1)\,\varepsilon_0 \int_0^T (\pi-S)^{-\sigma} m(S,w)\,dS\,.$$

Putting these arguments together yields

$$\sum_{|\alpha|\leq 8} \left( \|Z^\alpha w'(T,\cdot)\|_2 + \|Z^\alpha w(T,\cdot)\|_6 \right) \leq C \int_0^T m(S,w)\,dS$$
$$+ C\varepsilon_0 + C\,(C_0+1)^2\,\varepsilon_0^2 + C\,(C_0+1)\,\varepsilon_0 \sup_{0\leq S\leq T} m(S,w)\,. \quad (7.18)$$

We estimate the $L^\infty$-norms occurring in the definition of $m(T,w)$ using Theorem 6.1. By (6.1), if $p > 1$ is fixed,

$$\sum_{|\alpha|\leq 5} |Z^\alpha w(T,X)| \leq C \sup_{0\leq S\leq T} \sum_{|\alpha|\leq 6} \|Z^\alpha \mathcal{G}(S,\cdot)\|_2$$
$$+ C \sup_{0\leq S\leq T} \sum_{|\alpha|\leq 6} \left( \|Z^\alpha [\Box_g,\eta] u_0(S,\cdot)\|_2 + \|Z^\alpha [\Box_g,\eta] w(S,\cdot)\|_2 \right)$$
$$+ C \sup_{0\leq S\leq T} \sum_{|\alpha|\leq 6} (\pi-S)^{-2} \|Z^\alpha \mathcal{G}(S,\cdot)\|_p$$
$$+ C \sup_{0\leq S\leq T} \sum_{|\alpha|\leq 6} (\pi-S)^{-2} \left( \|Z^\alpha [\Box_g,\eta] u_0(S,\cdot)\|_p + \|Z^\alpha [\Box_g,\eta] w(S,\cdot)\|_p \right)$$
$$+ C \sup_{0\leq S\leq T} \sum_{I,j,k} \sum_{|\alpha|\leq 6} \|Z^\alpha (\gamma^{I,jk} \Gamma_j \Gamma_k w)(S,\cdot)\|_2$$
$$+ C \sup_{0\leq S\leq T} \sum_{I,j,k} \sum_{|\alpha|\leq 6} (\pi-S)^{-2} \|Z^\alpha (\gamma^{I,jk} \Gamma_j \Gamma_k w)(S,\cdot)\|_p\,.$$

Since $\pi - S$ is bounded below on the support of $[\Box_g,\eta]$, and $p < 2$, it follows that the fourth term on the right is dominated by the second. By (7.15), the second term is in turn dominated by

$$C\varepsilon_0 + C \sup_{0\leq S\leq T} \sum_{|\alpha|\leq 6} \left( \|Z^\alpha w'(S,\cdot)\|_2 + \|Z^\alpha w(S,\cdot)\|_6 \right).$$

Thus

$$\sum_{|\alpha|\leq 5} |Z^\alpha w(T,X)| \leq C\varepsilon_0 + C \sup_{0\leq S\leq T} \sum_{|\alpha|\leq 6} \left( \|Z^\alpha w'(S,\cdot)\|_2 + \|Z^\alpha w(S,\cdot)\|_6 \right)$$
$$+ I + II + III + IV\,,$$

if

$$I + II + III + IV = C \sup_{0 \leq S \leq T} \sum_{|\alpha| \leq 6} \|Z^\alpha \mathcal{G}(S, \cdot)\|_2$$

$$+ C \sup_{0 \leq S \leq T} \sum_{|\alpha| \leq 6} (\pi - S)^{-2} \|Z^\alpha \mathcal{G}(S, \cdot)\|_p$$

$$+ C \sup_{0 \leq S \leq T} \sum_{I,j,k} \sum_{|\alpha| \leq 6} \|Z^\alpha (\gamma^{I,jk} \Gamma_j \Gamma_k w)(S, \cdot)\|_2$$

$$+ C \sup_{0 \leq S \leq T} \sum_{I,j,k} \sum_{|\alpha| \leq 6} (\pi - S)^{-2} \|Z^\alpha (\gamma^{I,jk} \Gamma_j \Gamma_k w)(S, \cdot)\|_p.$$

We first note that
$$I \leq C (C_0 + 1)^2 \varepsilon_0^2.$$

To see this, we use the fact that $\mathcal{G}$ vanishes for $R \geq 2(\pi - T)$, together with Holder's inequality, to bound each of the terms in (7.5) by

$$\sum_{|\alpha| \leq 7} \left( \|Z^\alpha v\|_6^2 + \|Z^\alpha v\|_6^3 \right),$$

and then apply (7.14). Similarly,

$$III \leq C (C_0 + 1) \varepsilon_0 \sup_{0 \leq S \leq T} m(S, w),$$

where we use the fact that $\gamma^{I,jk}$ vanishes for $R \geq 2(\pi - T)$.

Again by the fact that $\mathcal{G} = 0$ for $R \geq 2(\pi - T)$, we can use (7.5) and Hölder's inequality to conclude that

$$\begin{aligned}
\sum_{|\alpha| \leq 6} (\pi - S)^{-2} \|Z^\alpha \mathcal{G}(S, \cdot)\|_p & \\
&\leq C \sum_{|\alpha| \leq 6} \|Z^\alpha v'(S, \cdot)\|_2 \\
&\quad \times \sum_{|\alpha| \leq 3} \left( \|Z^\alpha v'(S, \cdot)\|_{2p/(2-p)} + (\pi - S)^{-1} \|Z^\alpha v(S, \cdot)\|_{L^{2p/(2-p)}(R \leq 2(\pi - S))} \right) \\
&\quad + C \sum_{|\alpha| \leq 6} \|Z^\gamma v'(S, \cdot)\|_2 \cdot \sum_{|\alpha| \leq 3} \|Z^\alpha v\|_{4p/(2-p)}^2 \\
&\quad + C(\pi - S)^{-2} \|v\|_{L^{3p}(R \leq 2(\pi - S))}^3 + C(\pi - S)^{-2} \|v\|_{L^{2p}(R \leq 2(\pi - S))}^2. \quad (7.19)
\end{aligned}$$

By Hölder's inequality and (7.14), if $|\alpha| \leq 3$ we may bound

$$\|Z^\alpha v'(S, \cdot)\|_{2p/(2-p)} + (\pi - S)^{-1} \|Z^\alpha v(S, \cdot)\|_{L^{2p/(2-p)}(R \leq 2(\pi - S))}$$

$$\leq C \left( \|Z^\alpha v'(S, \cdot)\|_2 + \|Z^\alpha v(S, \cdot)\|_6 \right)^{2/p - 1} \left( \|Z^\alpha v'\|_\infty + (\pi - S)^{-1} \|Z^\alpha v\|_\infty \right)^{2 - 2/p}$$

$$\leq C (C_0 + 1) \varepsilon_0 \left( (\pi - S)^{-2-\sigma} + (\pi - S)^{-1-\sigma} \right)^{2 - 2/p}$$

$$\leq C (C_0 + 1) \varepsilon_0 (\pi - S)^{-(2+\sigma)(2 - 2/p)}.$$

Note that, given $\sigma' > 0$, one may choose $p$ small enough so that this is less than
$$C (C_0 + 1) \varepsilon_0 (\pi - S)^{-\sigma'},$$

where $C$ depends on $\sigma'$. Next, for $p$ close to 1, we have $4p/(2-p) < 6$, and consequently

$$\Big(\sum_{|\alpha|\leq 6} \|Z^\gamma v'(S,\cdot)\|_2\Big)\Big(\sum_{|\alpha|\leq 3}\|Z^\alpha v\|_{4p/(2-p)}^2\Big) \leq C\,(C_0+1)^3 \varepsilon_0^3\,.$$

The last two terms in (7.19) are similarly estimated. By the bounds that (7.14) implies on $\|v\|_\infty$, they are dominated by

$$C\,(C_0+1)^3\,\varepsilon_0^3\,(\pi-S)^{-2-3\sigma+3/p} + C\,(C_0+1)^2\,\varepsilon_0^2\,(\pi-S)^{-2-2\sigma+3/p}\,.$$

Therefore, since we are assuming that $\sigma = 1/4$, we conclude that if $p < 12/11$, and if $(C_0+1)\varepsilon_0 \leq 1$, then these terms are dominated by $C\,(C_0+1)^2\,\varepsilon_0^2$. Consequently, we have shown that, given $\sigma' > 0$, if $p$ is close enough to 1 (depending on $\sigma'$), then

$$(\pi-T)^{\sigma'} II \leq C\,(C_0+1)^2\,\varepsilon_0^2\,.$$

Similar arguments, using (2.17), yield

$$(\pi-T)^{\sigma'} IV \leq C\,(C_0+1)\,\varepsilon_0 \sup_{0\leq S\leq T} m(S,w)\,.$$

for $p$ close enough to 1 (depending on $\sigma'$).

Combining these steps with Theorem 6.1 yields that, for any $\sigma' > 0$, there exists $C$ depending on $\sigma'$ such that

$$(\pi-T)^{\sigma'}\sum_{|\alpha|\leq 5}\|Z^\alpha w(T,\cdot)\|_\infty \leq C\varepsilon_0 + C\,(C_0+1)^2\,\varepsilon_0^2 + C\,(C_0+1)\,\varepsilon_0 \sup_{0\leq S\leq T} m(S,w)$$
$$+ C \sup_{0\leq S\leq T}\sum_{|\alpha|\leq 6}\Big(\|Z^\alpha w'(S,\cdot)\|_2 + \|Z^\alpha w(S,\cdot)\|_6\Big)\,. \quad (7.20)$$

We now take $\sigma' = \sigma = 1/4$, and using (7.18) we obtain

$$m(T,w) \leq C\,\varepsilon_0 + C\,(C_0+1)^2\,\varepsilon_0^2 + C\,(C_0+1)\,\varepsilon_0 \sup_{0\leq S\leq T} m(S,w) + C\int_0^T m(S,w)\,dS\,.$$

If we let

$$M(T,w) = \sup_{0\leq S\leq T} m(S,w),$$

then the last inequality gives

$$M(T,w) \leq C\,\varepsilon_0 + C\,(C_0+1)^2\,\varepsilon_0^2 + C\,(C_0+1)\,\varepsilon_0\,M(T,w) + C\int_0^T M(S,w)\,dS\,.$$

By first taking $C_0$ large (depending on $C$), and then taking $\varepsilon_0$ small in order that $(C_0+1)^2\varepsilon_0 \leq 1$ and $(C_0+1)\varepsilon_0$ is sufficiently small, we may absorb the third term on the right into the left hand side, and then apply Gronwall's inequality to conclude that

$$M(T,w) \leq C_0\,\varepsilon_0\,. \quad \square$$

We now apply Lemma 7.2 to show that we can solve (7.12) by iteration. We assume that $u$ satisfies (7.8). For $\delta_0$ sufficiently small, Theorem 9.4 implies that $u_0$ satisfies (7.15). We now set $w_0 = 0$ and then define $w_k$, $k = 1, 2, 3, \ldots$, inductively by requiring that $w_k$ solve (7.16) with $v = v_k = u_0 + w_{k-1}$. Since $v_1 = u_0$ satisfies (7.14), we conclude by Lemma 7.2 that for all $k$

$$\sup_{0 \leq T < \pi} m(T, w_k) \leq C_0 \, \varepsilon \, .$$

In particular, the $w_k$ are a bounded sequence in $C^5(Y : T \leq T_0)$ for any $T_0 < \pi$. We now show that the sequence $w_k$ is Cauchy in the energy norm. It then follows that $w_k$ converges to a classical solution $w$ of (7.12), so that $u = u_0 + w$ satisfies (7.4); indeed

$$\sup_{0 \leq T < \pi} m(T, u) \leq (C_0 + 1) \, \varepsilon \, ,$$

where $\varepsilon$ can be taken as a constant multiple of $\delta_0$ if $\delta_0$ is sufficiently small.

To show that $w_k - w_{k-1}$ is Cauchy in the energy norm, we note that $w_k - w_{k-1}$ has vanishing initial data, and solves the Dirichlet-wave equation

$$\big((\Box_g + 1) - \sum_{jk} \gamma^{I,jk}(T, X; v_k, v'_k) \Gamma_j \Gamma_k \big)(w_k^I - w_{k-1}^I)$$
$$= \sum_{jk} \Big(\gamma^{I,jk}(T, X; v_k, v'_k) - \gamma^{I,jk}(T, X; v_{k-1}, v'_{k-1})\Big) \Gamma_j \Gamma_k w_{k-1}^I$$
$$+ \mathcal{G}^I(T, X; v_k, v'_k) - \mathcal{G}^I(T, X; v_{k-1}, v'_{k-1}), \quad I = 1, \ldots, N.$$

Recalling that $v_k = u_0 + w_{k-1}$, we can use the estimates (2.20) and (2.21), together with the fact that $m(w_k, T) \leq C_0 \, \varepsilon_0$, and the fact that $\mathcal{G}$ and $\gamma^{I,jk}$ are supported in the set $R \leq (\pi - T)$, together with Hölder's inequality, to bound

$$\sum_{jk} \big\| \gamma^{I,jk}(T, X; v_k, v'_k) - \gamma^{I,jk}(T, X; v_{k-1}, v'_{k-1}) \big\|_2 \big\| \Gamma_j \Gamma_k w_{k-1}^I \big\|_\infty$$
$$+ \|\mathcal{G}^I(T, X; v_k, v'_k) - \mathcal{G}^I(T, X; v_{k-1}, v'_{k-1})\|_2$$
$$\leq C \, \varepsilon_0 \, (\pi - S)^{-1/2} \Big( \|w'_{k-1}(S, \cdot) - w'_{k-2}(S, \cdot)\|_2 + \|w_{k-1}(S, \cdot) - w_{k-2}(S, \cdot)\|_6 \Big).$$

By Corollary 4.2, we thus have

$$\sup_{0 \leq T < \pi} \Big( \|w'_k(T, \cdot) - w'_{k-1}(T, \cdot)\|_2 + \|w_k(T, \cdot) - w_{k-1}(T, \cdot)\|_6 \Big)$$
$$\leq C \, \varepsilon_0 \sup_{0 \leq T < \pi} \Big( \|w'_{k-1}(T, \cdot) - w'_{k-2}(T, \cdot)\|_2 + \|w_{k-1}(T, \cdot) - w_{k-2}(T, \cdot)\|_6 \Big),$$

which, for $\varepsilon_0$ small, implies that $w_k$ is Cauchy in the energy norm.

It remains to show that the solution $u$ satisfies (7.4) for all $\sigma > 0$, since the iteration yields this only for $\sigma = 1/4$. This, however, is an easy consequence of (7.20), where we take $w = v = u$, since this estimate works for all $\sigma' > 0$, where $C$ depends on $\sigma'$. □

## 8. Global Existence in Minkowski Space.

Suppose that

$$\partial_t^2 \tilde{u} - \Delta \tilde{u} = F(\tilde{u}, d\tilde{u}, d^2 \tilde{u}) \tag{8.1}$$

is an equation in $\mathbb{R}_+ \times \mathbb{R}^3 \backslash \mathcal{K}$ satisfying the hypotheses of Theorem 1.2. Then, by Proposition 2.4, the Penrose compactification transforms this to an equation on the complement of $\mathcal{P}(\mathcal{K})$ in the Einstein diamond of the form $(\Box_g + 1)u = \mathcal{F}_0(T, X; u, du, d^2 u)$, where $\mathcal{F}_0$ can be extended to all of $([0, \pi) \times S^3) \backslash \mathcal{K}_*$ so that the conditions (2.16)–(2.21) are satisfied.

We fix a cutoff function $\eta \in C^\infty([0, \pi) \times S^3)$ satisfying $\eta(T, X) = 1$ if $R \leq (\pi - T)$ and $\eta(T, X) = 0$ if $R \geq 2(\pi - T)$, such that $\Gamma^\alpha \eta = O((\pi - T)^{-|\alpha|})$. Then

$$\mathcal{F}(T, X, u, du, d^2u) = \eta \mathcal{F}_0(T, X, u, du, d^2u)$$

satisfies the hypotheses of Theorem 7.1.

Suppose that we are given Cauchy data $(\tilde{f}, \tilde{g})$ for (8.1) which satisfies the compatibility conditions of order 8, and such that the smallness condition (1.13) holds. It then follows from (2.22) that the corresponding data $(f, g)$ in $S^3 \backslash \mathcal{P}_0(\mathcal{K})$ satisfies the hypotheses of Theorem 7.1. Thus, we can solve the Dirichlet-wave equation $(\Box_g + 1)u = \mathcal{F}(T, X; u, du, d^2u)$ in $([0, \pi) \times S^3) \backslash \mathcal{K}_*$. Since $\mathcal{F} = \mathcal{F}_0$ in the Einstein diamond, the pullback of $\tilde{u} = \Omega u$ to Minkowski space gives a solution of (1.1). Using the fact that a set $0 \leq t \leq t_0$ is mapped by $\mathcal{P}$ to a set on which $(\pi - T)$ is bounded away from 0, we have the following result,

**Theorem 8.1.** *Let $\mathcal{K}$ and $F(u, du, d^2u)$ be as in Theorem 1.1. Assume further that the Cauchy data satisfies compatibility conditions of order 8, as well as the smallness condition (1.13). Then there is a solution $u = \tilde{u}$ of (1.1), such that for all $t_0 < \infty$,*

$$\partial_t^j \partial_x^\alpha u \in L^\infty([0, t_0] \times \mathbb{R}^3 \backslash \mathcal{K}), \quad j + |\alpha| \leq 5,$$

*and*

$$\partial_t^j \partial_x^\alpha \tilde{u} \in L^2_{\text{loc}}([0, \infty) \times \mathbb{R}^3 \backslash \mathcal{K}), \quad j + |\alpha| \leq 9.$$

It is possible to use energy estimates in the Minkowski space, analogous to but simpler than Theorem 5.1, to show that in fact

$$\partial_t^j \partial_x^\alpha \tilde{u} \in L^\infty_t L^2_x([0, t_0] \times \mathbb{R}^3 \backslash \mathcal{K}), \quad j + |\alpha| \leq 9,$$

for all $t_0 < \infty$.

The solution $\tilde{u}$ also verifies the decay condition (1.14). This follows from the fact that the corresponding solution $u$ in the Einstein diamond verifies $|u| \leq C (\pi - T)^{-\sigma}$, and if we pull back this estimate to Minkowski space we obtain (1.14) for the corresponding function $\tilde{u}$, as can be seen by (1.17).

To complete the proof of Theorem 1.1, we need to show that if the data $(\tilde{f}, \tilde{g})$ are smooth, and satisfy the compatibility conditions of infinite order, then the solution $\tilde{u}$ must be smooth. This fact is an easy consequence of Theorem 9.5 with $s = 8$, where we note that we can apply that theorem locally by finite propagation velocity. □

## 9. Compatibility Conditions and Local Existence.

In this section we discuss the compatibility conditions for equations of the form (1.1), as well as establish the local existence theorems necessary for section 7. The existence theorem is known for the obstacle free problem; see e.g. [6] Theorem 6.4.11. Hence we concern ourselves with local existence near the boundary, and thus work on a compact manifold $\Sigma$ with smooth boundary. For convenience, we assume that $\Sigma$ is contained in the $n$-torus $\mathbb{T}^n$, so that we may write differential operators on $\Sigma$ in terms of the $\partial_{x_i}$.

All of our arguments work on general compact Riemannian manifolds with boundary, though, by using coordinate patches.

We will use $\partial_0$ interchangeably with $\partial_t$, and let $\partial_i = \partial_{x_i}$ for $1 \leq i \leq n$. We let $\partial$ denote the full collection of $\partial_i$ for $0 \leq i \leq n$, and $\partial_x$ denote the collection with $1 \leq i \leq n$. We also use $J_k u$ to denote the collection of all spatial derivatives of $u$ up to order $k$,

$$J_k u = \{\partial_x^\alpha u : 0 \leq |\alpha| \leq k\}.$$

We consider a quasilinear Cauchy problem of the form

$$\partial_t^2 u = \Delta u + \sum_{i,j=0}^n \gamma^{ij}(t, x, J_1 u, \partial_t u)\, \partial_i \partial_j u + G(t, x, J_1 u, \partial_t u), \tag{9.1}$$

$$u(0, x) = f(x), \qquad \partial_t u(0, x) = g(x), \qquad u(t, x) = 0 \text{ if } x \in \partial \Sigma.$$

Throughout this section, we assume that $\gamma^{ij}$ and $G$ are smooth functions of their arguments, with smooth extensions across the boundary of $\Sigma$ to all of $\mathbb{T}^n$. We assume that, for all values of its arguments,

$$\sum_{i,j=0}^n |\gamma^{ij}| \leq \frac{1}{2},$$

so that the equation is hyperbolic. By dividing by $1 + \gamma^{00}$, we will also assume that

$$\gamma^{00} = 0.$$

Consequently, we may write the equation in the form

$$\partial_t^2 u = F(t, x, J_2 u, J_1 \partial_t u), \tag{9.2}$$

$$u(0, x) = f(x), \qquad \partial_t u(0, x) = g(x).$$

Given a Cauchy problem of the form (9.2), there are associated *compatibility functions* $\psi_k$. The first few are explicitly given by

$$\psi_0 = f, \quad \psi_1 = g, \quad \psi_2 = F(t, x, J_2 f, J_1 g).$$

For $k \geq 2$, we note that we may formally write

$$\partial_t^{k-2} F(t, x, J_2 u, J_1 \partial_t u)$$
$$= \sum F_{\alpha_1, j_1, \ldots, \alpha_m, j_m}(t, x, J_2 u, J_1 \partial_t u) \left(\partial_x^{\alpha_1} \partial_t^{j_1} u\right) \cdots \left(\partial_x^{\alpha_m} \partial_t^{j_m} u\right),$$

where the functions $F_{\alpha_1, j_1, \ldots, \alpha_m, j_m}$ are smooth in their arguments, and where for each term in the sum there are numbers $n_i$ with $\sum n_i \leq k - 2$, such that

$$|\alpha_i| + j_i \leq 2 + n_i, \qquad j_i \leq 1 + n_i.$$

In particular, $j_i \leq k - 1$, and we may thus recursively define the $\psi_j$ by the procedure

$$\psi_k = \sum F_{\alpha_1, j_1, \ldots, \alpha_m, j_m}(0, x, J_2 f, J_1 g) \left(\partial_x^{\alpha_1} \psi_{j_1}\right) \cdots \left(\partial_x^{\alpha_m} \psi_{j_m}\right). \tag{9.4}$$

Since $|\alpha_i| + j_i \leq k$, it follows by induction that $\psi_k$ may be written in the form

$$\psi_k = \psi_k(J_k f, J_{k-1} g),$$

meaning that $\psi_k$ may be written as some function, the form of which depends on $F$, of the variables $J_k f$ and $J_{k-1} g$.

The interpretation of the compatibility functions is that, if $u$ is a smooth solution of the Cauchy problem (9.1), then necessarily $\partial_t^k u(0, \cdot) = \psi_k(J_k f, J_{k-1} g)$.

**Lemma 9.1.** *Suppose that $f \in H^{s+1}(M)$, $g \in H^s(M)$, where $s \geq n + 2$. Then the function $\psi_k(J_k f, J_{k-1} g)$ belongs to $H^{s+1-k}(M)$, for $0 \leq k \leq s + 1$.*

We prove this inductively. Thus, assume that it holds for $k - 1$, and take $k \geq 2$, the result being trivial for $k = 0, 1$. Consider the case of $k$ odd. We note that, from the condition $\sum n_i \leq k - 2$, there is at most one index $i$ with $n_i > \frac{k-3}{2}$. Consequently, for all indices $i$ in any given term in (9.4), with at most one exception,

$$\partial_x^{\alpha_i} \psi_{j_i} \in H^{s - \frac{k-1}{2}} \subseteq H^{\frac{n+1}{2}} \cap H^{s+1-k},$$

and the last space is an algebra of functions. Also, $J_2 f$, $J_1 g \in H^{s-1} \subseteq H^{\frac{n+1}{2}} \cap H^{s+1-k}$, so that

$$F_{\alpha_1, j_1, \ldots, \alpha_m, j_m}(0, x, J_2 f, J_1 g) \in H^{\frac{n+1}{2}} \cap H^{s+1-k}.$$

The result follows, since for the remaining index $i$ we have $\partial_x^{\alpha_i} \psi_{j_i} \in H^{s+1-k}$, since $|\alpha_i| + j_i \leq k$.

For $k$ even, there is at most one index $i$ with $n_i > \frac{k-2}{2}$, and the same proof goes through, noting that $s - \frac{k}{2} \geq \frac{n+1}{2}$ if $k$ is even. □

**Definition 9.2.** *For a Dirichlet-Cauchy problem of the form (9.1), with Cauchy data $f \in H^{s+1}(\Sigma)$, $g \in H^s(\Sigma)$, we say that the compatibility conditions of order $s$ are satisfied if $\psi_j(x)$ vanishes on $\partial \Sigma$, for all $0 \leq j \leq s$.*

The compatibility conditions are thus a (possibly nonlinear) condition on the Cauchy data $f, g$, and are a necessary condition to produce solutions $u$ to the Dirichlet-Cauchy problem of regularity $s+1$. We will also use compatibility conditions for linear equations that arise from (9.10). The compatibility functions and conditions for such equations have the obvious meaning, and in fact are linear in the data $f, g$; see e.g. [8] equation (2.30).

**Lemma 9.3.** *Assume that $f \in H^{s+1}$ and $g \in H^s$, where $s \geq n$ if $n$ is odd, and $s \geq n+1$ if $n$ is even. Let $\psi_k(J_k f, J_{k-1} g)$ be the compatibility functions for the Cauchy problem (9.1). Suppose that $v(t, x)$ is a function such that, for some $T > 0$,*

$$\partial_t^j v \in C\big([0, T); H^{s+1-j}(M)\big), \quad \text{for} \quad 0 \leq j \leq s+1,$$

*and suppose that for $0 \leq k \leq s$,*

$$\partial_t^k v(0, \cdot) = \psi_k(J_k f, J_{k-1} g).$$

*Let $\overline{\psi}_j$ be the compatibility functions for the Cauchy problem*

$$\partial_t^2 u = \Delta u + \sum_{i,j=0}^{n} \gamma^{ij}(t, x, J_1 v, \partial_t v) \, \partial_i \partial_j u + G(t, x, J_1 v, \partial_t v), \qquad (9.5)$$

$$u(0, x) = f(x), \qquad \partial_t u(0, x) = g(x).$$

*Then for $0 \leq k \leq s$,*

$$\overline{\psi}_k = \psi_k(J_k f, J_{k-1} g).$$

We check this by induction, the result being immediate for $k = 0, 1, 2$. Assume thus the result holds for $0 \leq j \leq k-1$. That $\overline{\psi}_k = \psi_k$ then follows by noting that, if we apply $\partial_t^{k-2}$ to the right hand side of (9.5), set $t = 0$, and then substitute $\partial_t^j u = \psi_j$, then we obtain the same result as applying $\partial_t^{k-2}$ to the right hand side of (9.1), followed by setting $t = 0$ and $\partial_t^j u = \psi_j$. □

The main results of this section are the following two theorems.

**Theorem 9.4.** *Consider the Dirichlet-Cauchy problem* (9.1), *with data* $f \in H^{s+1}(\Sigma)$, $g \in H^s(\Sigma)$, *where* $s \geq (3n+6)/2$ *if* $n$ *is even*, $s \geq (3n+3)/2$ *if* $n$ *is odd, and* $s \geq 4$ *if* $n = 1$. *Suppose that the compatibility conditions of order* $s$ *are satisfied by the data. Then there exists* $T > 0$, *depending on* $s$ *and bounds on the norms of* $f$ *and* $g$, *such that there exists a solution* $u$ *to* (9.1) *on* $[0, T] \times \Sigma$, *which satisfies*

$$\sup_{0 \leq t \leq T} \sum_{j=0}^{s+1} \|\partial_t^j u(t, \cdot)\|_{H^{s+1-j}} \leq C < \infty.$$

*Furthermore, if* $\bigl(\|f\|_{H^{s+1}} + \|g\|_{H^s}\bigr) \leq 1$, *and* $G(t, x, 0, 0) = 0$, *then there exists* $C$ *and* $T$ *independent of* $f$ *and* $g$, *so that the solution exists for* $0 \leq t \leq T$, *and satisfies*

$$\sup_{0 \leq t \leq T} \sum_{j=0}^{s+1} \|\partial_t^j u(t, \cdot)\|_{H^{s+1-j}} \leq C \bigl(\|f\|_{H^{s+1}} + \|g\|_{H^s}\bigr). \tag{9.6}$$

Theorem 9.4 does not yield existence of $C^\infty$ solutions, since $T$ may depend on $s$. However, the following result together with the above does imply local existence of solutions of arbitrarily high smoothness for (9.1). We will also use the next theorem to establish existence of global $C^\infty$ solutions for our original equation (1.1).

**Theorem 9.5.** *Suppose that the conditions of Theorem 9.4 are satisfied by the integer* $s$, *and suppose that* $u$ *is a solution to* (9.1) *on an interval* $[0, T']$, *such that*

$$\sup_{0 \leq t \leq T'} \sum_{j=0}^{s+1} \|\partial_t^j u(t, \cdot)\|_{H^{s+1-j}} < \infty.$$

*Suppose that* $m > s$, *that* $f \in H^{m+1}(\Sigma)$, $g \in H^m(\Sigma)$, *and that the compatibility conditions of order* $m$ *are satisfied. Then*

$$\sup_{0 \leq t \leq T'} \sum_{j=0}^{m+1} \|\partial_t^j u(t, \cdot)\|_{H^{m+1-j}} < \infty.$$

The proof of Theorems 9.4 and 9.5 will be based on a priori estimates for solutions $u$ to the linearized Cauchy problem (9.5). We let

$$M_{s+1}(v, t) = \sum_{|\alpha| \leq s+1} \|\partial^\alpha v(t, \cdot)\|_{L^2(\Sigma)}.$$

In the following lemma, we assume that $s$ satisfies the conditions of Theorem 9.4, although the proof works for $s$ as in Klainerman's argument for the obstacle-free case, on which our argument is based. We refer to the treatment on page 117 of Hörmander [6].

**Lemma 9.6.** *Let $u$ be a solution to the equation (9.5), where we assume that*
$$\partial_t^j u,\, \partial_t^j v \in C\big([0,T); H^{s+1-j}(\Sigma)\big), \quad \text{for} \quad 0 \leq j \leq s+1,$$
*and that $u(t,x) = 0$ for $x \in \partial\Sigma$. Suppose also that*
$$\sup_{0 \leq t \leq T} \sum_{|\alpha| \leq 2} \|\partial^\alpha u(t,\cdot)\|_\infty < M,$$
$$\sup_{0 \leq t \leq T} \sum_{|\alpha| \leq 2} \|\partial^\alpha v(t,\cdot)\|_\infty < M.$$
*Then there exists a constant $C$, independent of $u$ and $v$, such that*
$$M_{s+1}(u,t) \leq C\, e^{CMt} \bigg( M_{s+1}(u,0) + C[M] \int_0^t \big( M_{s+1}(u,r) + M_{s+1}(v,r) \big)\, dr$$
$$+ \int_0^t C[M_s(v,r)]\, (M_s(u,r) + 1)\, dr \bigg) + C[M_s(v,t)]\, (M_s(u,t) + 1). \quad (9.7)$$
*Here, $C[\,\cdot\,]$ denotes a constant that depends on the quantity inside the brackets.*

Let $V$ denote a collection of $n+1$ vector fields on $\mathbb{R} \times \Sigma$ which are tangent to $\partial\Sigma$, and which span the Lie algebra of all such vector fields. Then, if $|\alpha| \leq s$, the function $V^\alpha u$ is an $H^1$ solution to the following equation
$$\big(\partial_t^2 - \Delta - \sum_{i,j=0}^n \gamma^{ij}(t,x,J_1 v,\partial_t v)\, \partial_i \partial_j \big)(V^\alpha u)$$
$$= \big[V^\alpha, \Box - \sum_{ij} \gamma^{ij} \partial_i \partial_j \big] u + V^\alpha G(t,x,J_1 v, \partial_t v),$$
such that $V^\alpha u$ vanishes on $\partial\Sigma$. To begin, we bound
$$\big\| [\Box, V^\alpha]\, u \big\|_2 \leq C\, M_{s+1}(u,t).$$
Next, we write
$$\big[\gamma^{ij}(t,x,J_1 v, \partial_t v),\, V^\alpha\big] \partial_i \partial_j u = \sum_{\substack{|\alpha_1|+|\alpha_2|\leq s+2 \\ |\alpha_1|\leq s,\, |\alpha_2|\leq s+1}} b_{\alpha_1,\alpha_2}(t,x)\, \big(\partial^{\alpha_1} \gamma^{ij}(t,x,J_1 v, \partial_t v)\big)\, \big(\partial^{\alpha_2} u\big),$$
where the $b_{\alpha_1,\alpha_2}$ are smooth functions. By considering the terms $|\alpha_1| = s$ or $|\alpha_2| = s+1$, and then the remaining terms, we may bound the $L^2$ norm by
$$C[M]\, \big(M_{s+1}(u,t) + M_{s+1}(v,t)\big) + C[M_s(v,t)]\, M_s(u,t).$$
We may also bound
$$\|V^\alpha G(t,x,J_1 v, \partial_t v)\|_2 \leq C[M]\, M_{s+1}(v,t) + C[M_s(v,t)].$$
By energy estimates (see [6], Proposition 6.3.2 for the obstacle-free version), and the fact that $|\partial_t \gamma^{ij}| \leq C\, M$, we have the following bounds,
$$\|\partial V^\alpha u(t,\cdot)\|_2 \leq C\, e^{CMt} \bigg( M_{s+1}(u,0) + C[M] \int_0^t \big(M_{s+1}(u,r) + M_{s+1}(v,r)\big)\, dr$$
$$+ \int_0^t C[M_s(v,r)]\, (M_s(u,r) + 1)\, dr \bigg). \quad (9.8)$$

We now work in geodesic normal coordinates near $\partial\Sigma$ such that $x_1$ is the normal direction. We thus may write
$$\Box = -\partial_1^2 + \Box',$$
where $\Box'$ involves derivatives of order at most 1 in the $x_1$ direction. Then by equation (9.5) we can write
$$\partial^\alpha \partial_1^2 u = \partial^\alpha\Big( \big(1+\gamma^{11}\big)^{-1} \big(\Box' u - \sum_{i,j\neq(1,1)} \gamma^{ij}\,\partial_i\partial_j u - G\big)\Big).$$

For each multi-index $\alpha$, the expression on the right involves one lower power in $\partial_1$ than the left hand side. Furthermore, for $|\alpha| \leq s-1$, any term on the right in which $u$ is differentiated $s+1$ times involves zero derivatives falling on $\gamma^{ij}$, and thus can be estimated using the fact that $|\gamma^{ij}| \leq \frac{1}{2}$. We thus have
$$\|\partial^\alpha \partial_1^2 u\|_2 \leq C \sum_{|\beta|\leq s+1} \|\partial^\beta u\|_2 + C\big[M_s(v,t)\big]\big(M_s(u,t)+1\big),$$
where the sum is over $\beta$ of strictly lower order in $\partial_1$ than the left hand side. Since (9.8) gives control over derivatives of order at most 1 in $\partial_1$, a simple induction in the order of $\partial_1$ completes the proof of the lemma. $\square$

The above estimate will be used to prove Theorem 9.5. For Theorem 9.4, we need a variation which can be iterated for small $T$.

**Lemma 9.7.** *Let $u$ be a solution to the equation* (9.5), *where we assume that*
$$\partial_t^j u,\, \partial_t^j v \in C\big([0,T); H^{s+1-j}(\Sigma)\big), \quad \text{for} \quad 0 \leq j \leq s+1,$$
*and that $u(t,x) = 0$ for $x \in \partial\Sigma$. Let*
$$M_{s+1}(v) = \sup_{0 \leq t \leq T} \sum_{|\alpha|\leq s+1} \|\partial^\alpha v(t,\cdot)\|_2$$
*and suppose that $M_{s+1}(v) < \infty$. Then there exists a constant $C$, independent of $u$ and $v$, such that for $0 \leq t < T$,*
$$M_{s+1}(u,t) \leq C\,e^{CM_{s+1}(v)t}\left( M_{s+1}(u,0) + C\big[M_{s+1}(v)\big] \int_0^t \big(M_{s+1}(u,r)+1\big)\,dr \right)$$
$$+ T^{\frac{1}{3}} C\big[M_{s+1}(v)\big]\big(M_{s+1}(u,t)+1\big). \quad (9.9)$$

To prove this, we begin by letting
$$\mathcal{B} = \{x \in \Sigma : \operatorname{dist}(x,\partial\Sigma) \leq 2t\,\}.$$
Since $\sum_{ij} |\gamma^{ij}| < \frac{1}{2}$, it follows that the complement of $\mathcal{B}$ is causal, in the sense that it contains the domain of influence for each of its points. On the open set $\mathcal{B}^c$, for any multi-index $\alpha$ with $|\alpha| \leq s$, the following holds,
$$\big(\partial_t^2 - \Delta - \sum_{i,j=0}^n \gamma^{ij}(t,x,J_1 v,\partial_t v)\,\partial_i\partial_j\big)(\partial^\alpha u)$$
$$= -\big[\partial^\alpha, \textstyle\sum_{ij} \gamma^{ij}\partial_i\partial_j\big]u + \partial^\alpha G(t,x,J_1 v,\partial_t v).$$
The $L^2$ norm of the right hand side over all of $\Sigma$ is bounded by
$$C\big[M_{s+1}(v,t)\big]\big(M_{s+1}(u,t)+1\big).$$

Since $\gamma^{ij} \in C^2$, we may use domain of dependence arguments to conclude that

$$\sum_{|\alpha| \leq s+1} \|\partial^\alpha u(t, \cdot)\|_{L^2(\mathcal{B}_t^c)} \leq C\, e^{CMt} \left( M_{s+1}(u, 0) + C[M_{s+1}(v)] \int_0^t \left( M_{s+1}(u, r) + 1 \right) dr \right).$$

It remains to estimate the norms over $\mathcal{B}_t$. By the arguments leading to (9.8), we may bound

$$\|\partial V^\alpha u(t, \cdot)\|_2 \leq C\, e^{CMt} \left( M_{s+1}(u, 0) + C[M_{s+1}(v)] \int_0^t \left( M_{s+1}(u, r) + 1 \right) dr \right)$$

for the vector fields $V$ that are tangent to $\partial \Sigma$. Following the proof of Lemma 9.6, we work in geodesic normal coordinates near $\partial \Sigma$ such that $x_1$ is the normal direction, so that we can write, for $|\alpha| \leq s - 1$,

$$\partial^\alpha \partial_1^2 u = \partial^\alpha \left( \left(1 + \gamma^{11}\right)^{-1} \left( \Box' u - \sum_{i,j \neq (1,1)} \gamma^{ij} \partial_i \partial_j u - G \right) \right),$$

where the right hand side involves one lower power in $\partial_1$ than the left hand side. We may thus bound

$$\|\partial^\alpha \partial_1^2 u(t, \cdot)\|_{L^2(\mathcal{B}_t)} \leq C \sum_{|\beta| \leq s+1} \|\partial^\beta u(t, \cdot)\|_2$$
$$+ C[M_{s+1}(v)] \left( \sum_{i,j} \sum_{\substack{|\alpha_1| + |\alpha_2| \leq s+1 \\ |\alpha_1| \leq s-1, |\alpha_2| \leq s}} \left\| \partial^{\alpha_1} \gamma^{ij}\, \partial^{\alpha_2} u \right\|_{L^2(\mathcal{B}_t)} + \sum_{|\theta| \leq s-1} \|\partial^\theta G\|_{L^2(\mathcal{B}_t)} \right),$$

where $\beta$ is of lower order in $\partial_1$ than the left hand side. In the second sum, a term with $|\alpha_2| \geq s + 1 - \frac{n}{2}$ may be dominated by

$$C[M_{s+1}(v)]\, \|\partial^{\alpha_2} u(t, \cdot)\|_{L^2(\mathcal{B}_t)},$$

since in this case $|\alpha_1| \leq \frac{n}{2} < s - \frac{n}{2}$. By Holder's inequality and Sobolev embedding, we may bound

$$\|\partial^{\alpha_2} u(t, \cdot)\|_{L^2(\mathcal{B}_t)} \leq C\, T^{\frac{1}{3}} \|\partial^{\alpha_2} u(t, \cdot)\|_{L^6(\mathcal{B}_t)} \leq C\, T^{\frac{1}{3}}\, M_{s+1}(u, t).$$

Terms with $|\alpha_1| \geq s - \frac{n}{2}$ are bounded by

$$C\, M_{s+1}(u, t)\, \|\partial^{\alpha_1} \gamma^{ij}\|_{L^2(\mathcal{B}_t)} \leq T^{\frac{1}{3}}\, C[M_{s+1}(v)]\, M_{s+1}(u, t)$$

by similar arguments. For the same reasons,

$$\sum_{|\theta| \leq s-1} \|\partial^\theta G\|_{L^2(\mathcal{B}_t)} \leq T^{\frac{1}{3}}\, C[M_{s+1}(v)].$$

Thus,

$$\|\partial^\alpha \partial_1^2 u(t, \cdot)\|_{L^2(\mathcal{B}_t)} \leq C \sum_{|\beta| \leq s+1} \|\partial^\beta u(t, \cdot)\|_2 + T^{\frac{1}{3}}\, C[M_{s+1}(v)] \left( M_{s+1}(u, t) + 1 \right),$$

where $\beta$ is of lower order in $\partial_1$ than the left hand side. Induction on the order of $\partial_1$ now completes the proof. $\square$

We will produce a solution to the Cauchy problem (9.1) by iteration. The first step is showing that solutions to the linearized equation (9.5) exist, after which we may apply Lemma 9.7 to obtain a priori bounds which iterate for small $T$. Our existence result for

50                    MARKUS KEEL, HART F. SMITH, AND CHRISTOPHER D. SOGGE

(9.5) is a simple extension of results of Ikawa [8]. In particular, Theorems 1 and 2 of that paper together imply the following result.

**Theorem 9.8.** *Consider the linear equation*

$$\partial_t^2 u = \Delta u + \sum_{i,j=0}^n \gamma^{ij}(t,x)\, \partial_i \partial_j u + G(t,x), \tag{9.10}$$

$$u(0,x) = f(x), \qquad \partial_t u(0,x) = g(x), \qquad u(t,x) = 0 \ \text{if}\ x \in \partial\Sigma,$$

*where* $\sum_{i,j} \|\gamma^{ij}\|_\infty \leq \frac{1}{2}$.

*Assume that* $\gamma^{ij}(t,x) \in C^k([0,T] \times \Sigma)$, *that* $f \in H^k(\Sigma)$, $g \in H^{k-1}(\Sigma)$, *and that*

$$\partial_t^j G \in C([0,T]; H^{k-2-j}(\Sigma)), \quad 0 \leq j \leq k-2,$$

$$\partial_t^{k-1} G \in L^1([0,T]; L^2(\Sigma)).$$

*Then if the compatibility conditions of order* $k-1$ *are satisfied, there exists a solution* $u$ *to* (9.10) *with* $\partial_t^j u \in C([0,T]; H^{k-j}(\Sigma))$ *for* $0 \leq j \leq k$, *and furthermore*

$$\sup_{|\alpha| \leq k} \|\partial^\alpha u(t,\cdot)\|_2$$

$$\leq C \left( \|f\|_{H^k} + \|g\|_{H^{k-1}} + \sup_{0 \leq r \leq t} \sup_{|\alpha| \leq k-2} \|\partial^\alpha G(r,\cdot)\|_2 + \int_0^t \|\partial_r^{k-1} G(r,\cdot)\|_2\, dr \right).$$

The last inequality is not explicitly stated in [8], but follows immediately from the proof of Theorem 2 of that paper. The constant $C$ depends on $\Sigma$, $k$, and the $C^k$ norm of $\gamma^{ij}$, but not the data. We note that the compatibility conditions of order $k-1$ are well defined for any $k$, since the equation is linear.

We will extend this theorem to the case that the $\gamma^{ij}$ satisfy Sobolev regularity conditions. For this, we need the following elementary elliptic regularity result, the proof of which we include for completeness.

**Lemma 9.9.** *Suppose that* $\gamma^{ij}(x) \in H^m(\Sigma)$, *where* $m > 2 + \frac{n}{2}$, *and* $\sum_{ij=1}^n \|\gamma^{ij}\|_\infty \leq \frac{1}{2}$. *Let* $u \in H^1(\Sigma)$ *satisfy the equation*

$$\Delta u(x) + \sum_{i,j=1}^n \gamma^{ij}(x)\, \partial_i\, \partial_j u(x) \;=\; F(x), \quad u(x) = 0 \ \text{if}\ x \in \partial\Sigma,$$

*Then if* $k \leq m$, *and* $F \in H^k(\Sigma)$, *it follows that* $u \in H^{k+2}(\Sigma)$, *and*

$$\|u\|_{H^{k+2}} \leq C_\gamma\, \|F\|_{H^k},$$

*where the constant* $C_\gamma$ *depends on* $\|\gamma^{ij}\|_{H^m}$, *but not* $F$.

We begin by noting that the conditions imply that $\gamma^{ij} \in C^2(\overline{\Sigma})$, so that the result holds for $k=0$ classically. We thus assume that the result is true for $k$ replaced by $k-1$, and show that it holds for $k$.

Let $V$ be a smooth vector field tangent to $\partial\Sigma$. Then $Vu$ is an $H^1$ solution to the equation

$$\Delta(Vu) + \sum_{i,j=1}^n \gamma^{ij}\, \partial_i\, \partial_j(Vu) \;=\; \left[\Delta + \sum_{ij}\gamma^{ij}\, \partial_i\, \partial_j,\, V\right] u + VF, \qquad (Vu)|_{\partial\Sigma} = 0.$$

Since $\gamma^{ij}$ and $\partial_x \gamma^{ij}$ are both multipliers on the space $H^{k-1}(\Sigma)$, we may bound
$$\left\| [\Delta + \sum_{ij} \gamma^{ij} \partial_i \partial_j, V] u \right\|_{H^{k-1}} \leq C \, \|\gamma^{ij}\|_{H^s} \, \|u\|_{H^{k+1}}.$$
By the induction hypothesis, we thus have
$$\|Vu\|_{H^{k+1}} \leq C_\gamma \|F\|_{H^k},$$
for $V$ smooth and tangent to $\partial \Sigma$. By working in local coordinates for which $\partial_i$ is tangent to $\partial \Sigma$ for $2 \leq i \leq n$, it follows that we control $\partial^\alpha u$ for all $|\alpha| = k+2$ with the exception of $\partial_1^{k+2} u$. We now use the elliptic equation for $u$ to express $\partial_1^{k+2} u$ in terms of derivatives involving at most $k+1$ factors of $\partial_1$, completing the proof. □

Our extension of Ikawa's result produces solutions $u$ of regularity $s+1$ provided that $s$ is sufficiently large so that there exists an integer $k > 2 + \frac{n}{2}$ with
$$k < s - \frac{n}{2}, \qquad 2k > s + \frac{n}{2}.$$
If $n$ is even, this requires $s \geq (3n+6)/2$, in which case $k = s - 1 - \frac{n}{2}$ works. If $n \geq 3$ is odd, this requires $s \geq (3n+3)/2$, in which case $k = s - \frac{n+1}{2}$ works. If $n = 1$, then $s \geq 4$ and $k = s - 1$ works.

**Theorem 9.10.** *Consider the Dirichlet-Cauchy problem (9.10), where $\sum_{i,j} \|\gamma^{ij}\|_\infty \leq \frac{1}{2}$. Suppose that*
$$\gamma^{ij} \in C^j\big([0,T]; H^{s-j}(\Sigma)\big), \quad 0 \leq j \leq s,$$
*with $s$ as above. Suppose also that*
$$\partial_t^j G \in C\big([0,T]; H^{s-1-j}(\Sigma)\big), \quad 0 \leq j \leq s-1,$$
$$\partial_t^s G \in L^1\big([0,T]; L^2(\Sigma)\big).$$
*If $f \in H^{s+1}(\Sigma)$, $g \in H^s(\Sigma)$, and the compatibility conditions of order $s$ are satisfied, then equation (9.10) has a solution $u$ such that $\partial_t^j u \in C\big([0,T]; H^{s+1-j}(\Sigma)\big)$ for $0 \leq j \leq s+1$. Furthermore*
$$\sup_{|\alpha| \leq s+1} \|\partial^\alpha u(t, \cdot)\|_2$$
$$\leq C \left( \|f\|_{H^{s+1}} + \|g\|_{H^s} + \sup_{0 \leq r \leq t} \sup_{|\alpha| \leq s-1} \|\partial^\alpha G(r, \cdot)\|_2 + \int_0^t \|\partial_r^s G(r, \cdot)\|_2 \, dr \right).$$

To begin, we note that $\gamma^{ij} \in C^k\big([0,T] \times \Sigma\big)$, where $k$ depends on $s$ as above, so that Theorem 9.8 guarantees solutions of regularity $k$. To show that this solution is actually of regularity $s+1$, we follow Ikawa [8] and formally pose $w = \partial_t^{s+1-k} u$. We thus seek a solution to the following equation,

$$\partial_t^2 w - \Delta w - \sum_{i,j=0}^n \gamma^{ij} \partial_i \partial_j w \qquad (9.11)$$
$$= \sum_{i,j=0}^n \sum_{m=0}^{s-k} \binom{s+1-k}{m} \big(\partial_t^{s+1-k-m} \gamma^{ij}\big)\big(\partial_t^m \partial_i \partial_j u\big) + \partial_t^{s+1-k} G,$$
$$w(0,x) = \psi_{s+1-k}(x), \qquad \partial_t w(0,x) = \psi_{s+2-k}(x), \qquad w(t,x) = 0 \text{ if } x \in \partial \Sigma,$$

subject to the condition that

$$u(t,x) = \psi_0(x) + t\,\psi_1(x) + \cdots + \frac{t^{s-k}}{(s-k)!}\psi_{s-k}(x) + \int_0^t \frac{(t-r)^{s-k}}{(s-k)!}\,w(r,x)\,dr\,. \tag{9.12}$$

We begin by establishing estimates that will allow us to solve this equation by iteration. Suppose then that (9.12) is replaced by the condition

$$u(t,x) = \psi_0(x) + t\,\psi_1(x) + \cdots + \frac{t^{s-k}}{(s-k)!}\psi_{s-k}(x) + \int_0^t \frac{(t-r)^{s-k}}{(s-k)!}\,\widetilde{w}(r,x)\,dr\,, \tag{9.13}$$

and let $w$ be the solution to (9.11) of regularity $k$ guaranteed by Theorem 9.4. We seek bounds on $w$ in terms of $\widetilde{w}$.

Consider first the quantity

$$\sum_{m=0}^{s-k} \sup_{|\alpha_1|+|\alpha_2|\leq k-2} \left\|\left(\partial^{\alpha_1}\partial_t^{s+1-k-m}\gamma^{ij}\right)\left(\partial^{\alpha_2}\partial_t^m \partial_i\partial_j u\right)(t,\cdot)\right\|_2.$$

We may use the fact that $\gamma^{ij} \in C^k$ to bound the terms for which $|\alpha_1|+s+1-k-m \leq k$ by

$$C \sum_{m=0}^{s-k} \sup_{|\alpha|\leq k} \|\partial^\alpha \partial_t^m u(t,\cdot)\|_2\,.$$

Since the order of differentiation in $x$ is at most $k$, and the total order of differentiation at most $s$, by (9.13) this is in turn dominated (for bounded $t$) by

$$\sum_{j=0}^{s} \left\|\partial_t^j u(0,\cdot)\right\|_{H^{s-j}} + \sup_{|\alpha|\leq k} \int_0^t \left\|\partial^\alpha \widetilde{w}(r,\cdot)\right\|_2 dr\,. \tag{9.14}$$

If $|\alpha_1|+s+1-k-m > k$, then $|\alpha_2|+2+m < s+1-k$, and thus we may bound the remaining terms by

$$C \sup_{|\alpha|\leq s-k} \|\partial^\alpha u(t,\cdot)\|_\infty \leq C \sup_{|\alpha|\leq k} \|\partial^\alpha u(t,\cdot)\|_2$$

where we use the fact that $s-k < k - \frac{n}{2}$. Consequently, these terms are also dominated by the quantity (9.14).

Next consider the quantity

$$\sum_{\substack{m_1+m_2\leq s \\ m_2 \leq s-1}} \int_0^t \left\|\left(\partial_t^{m_1}\gamma^{ij}\right)\left(\partial_t^{m_2}\partial_i\partial_j u\right)(r,\cdot)\right\|_2 dr\,.$$

If $m_1 \leq k$, we may bound this by

$$C \sup_{m\leq s-1} \int_0^t \left\|\partial_t^m \partial_i\partial_j u(r,\cdot)\right\|_2 dr\,.$$

Since $\widetilde{w} = \partial_t^{s+1-k}u$, it follows that at most $k$ derivatives hit $\widetilde{w}$, and consequently this term is dominated by (9.14). Finally, if $m_1 > k$, then $m_2 \leq s-1-k$, and since $k > 2+\frac{n}{2}$

we may bound this by

$$C \sup_{m \leq s-1-k} \int_0^t \left\|\partial_t^m \partial_i \partial_j u(r,\cdot)\right\|_\infty dr \leq C \sup_{\substack{m \leq s-1-k \\ |\alpha| \leq k}} \int_0^t \left\|\partial_t^m \partial^\alpha u(r,\cdot)\right\|_2 dr,$$

which is similarly bounded by (9.14).

It follows by Theorem 9.8, that if $w$ satisfies (9.11), where $u$ is given by (9.13), then

$$\sup_{|\alpha| \leq k} \|\partial^\alpha w(t,\cdot)\| \leq C \Big( \sum_{j=0}^{s+1} \|\partial_t^j u(0,\cdot)\|_{H^{s+1-j}} + \sup_{|\alpha| \leq k} \int_0^t \|\partial^\alpha \widetilde{w}(r,\cdot)\|_2 dr$$
$$+ \sup_{0 \leq r \leq t} \sup_{|\alpha| \leq s-1} \|\partial^\alpha G(r,\cdot)\|_2 + \int_0^t \|\partial_t^s G(r,\cdot)\|_2 dr \Big). \quad (9.15)$$

Now consider the sequence of functions produced by the following iterative procedure,

$$\partial_t^2 w_{l+1} - \Delta w_{l+1} - \sum_{i,j=0}^n \gamma^{ij} \partial_i \partial_j w_{l+1}$$
$$= \sum_{i,j=0}^n \sum_{m=0}^{s-k} \binom{s+1-k}{m} (\partial_t^{s+1-k-m} \gamma^{ij})(\partial_t^m \partial_i \partial_j u_l) + \partial_t^{s+1-k} G(t,x),$$
$$w_{l+1}(0,x) = \psi_{s+1-k}(x), \quad \partial_t w_{l+1}(0,x) = \psi_{s+2-k}(x), \quad w_{l+1}(t,x) = 0 \text{ if } x \in \partial\Sigma,$$

where

$$u_l(t,x) = \psi_0(x) + t\,\psi_1(x) + \cdots + \frac{t^{s-k}}{(s-k)!}\psi_{s-k}(x) + \int_0^t \frac{(t-r)^{s-k}}{(s-k)!} w_l(r,x)\, dr,$$

and we set $w_0 \equiv 0$.

By (9.15), it follows that for each $l$

$$\sup_{|\alpha| \leq k} \|\partial^\alpha w_l(t,\cdot)\| < \infty.$$

Furthermore, for $l \geq k$, it is easy to see that

$$\partial_t^m (u_{l+1} - u_l)(0,x) = 0, \quad \text{if} \quad 0 \leq m \leq s+1.$$

We thus can apply (9.15) to the equation

$$\partial_t^2 (w_{l+1} - w_l) - \Delta(w_{l+1} - w_l) - \sum_{i,j=0}^n \gamma^{ij} \partial_i \partial_j (w_{l+1} - w_l)$$
$$= \sum_{i,j=0}^n \sum_{m=0}^{s-k} \binom{s+1-k}{m} (\partial_t^{s+1-k-m} \gamma^{ij}) \partial_t^m \partial_i \partial_j (u_l - u_{l-1}),$$

to obtain that, for $l \geq k$,

$$\sup_{|\alpha| \leq k} \|\partial^\alpha (w_{l+1} - w_l)(t,\cdot)\| \leq C \sup_{|\alpha| \leq k} \int_0^t \|\partial^\alpha (w_l - w_{l-1})(r,\cdot)\|_2 dr,$$

and hence
$$\sup_{|\alpha|\leq k} \|\partial^\alpha(w_{l+1}-w_l)(t,\cdot)\| \leq K\,\frac{(Ct)^{l-k}}{(l-k)!}\,,$$
where
$$K \leq C\left(\|f\|_{H^{s+1}} + \|g\|_{H^s} + \sup_{0\leq r\leq t}\sup_{|\alpha|\leq s-1}\|\partial^\alpha G(r,\cdot)\|_2 + \int_0^t \|\partial_r^s G(r,\cdot)\|_2\,dr\right).$$

It follows that the sequence $w_l$ converges to a limit $w$ such that $\partial_t^m w \in C([0,T];H^{k-m})$ for $0\leq m\leq k$. We define $u$ by equation (9.12), and following Ikawa we see that $u$ is a solution to (9.10), such that $\partial_t^{s+1-m}u \in C([0,T];H^m)$ for $0\leq m\leq k$, and furthermore

$$\sum_{m=0}^k \|\partial_t^{s+1-m}u(t,\cdot)\|_{H^m}$$
$$\leq C\left(\|f\|_{H^{s+1}} + \|g\|_{H^s} + \sup_{0\leq r\leq t}\sup_{|\alpha|\leq s-1}\|\partial^\alpha G(r,\cdot)\|_2 + \int_0^t \|\partial_r^s G(r,\cdot)\|_2\,dr\right). \quad (9.16)$$

We now establish bounds on the higher spatial derivatives by elliptic regularity. Suppose that we have shown $\partial_t^{s+1-m}u \in C([0,T];H^m)$ for $0\leq m\leq p$, where $p$ is some integer with $k\leq p\leq s+1$, and that (9.16) holds with $k$ replaced by $p$. We write

$$\left(\Delta + \sum_{i,j=1}^n \gamma^{ij}(t,x)\,\partial_i\,\partial_j\right)\!\left(\partial_t^{s-p}u(t,\cdot)\right) = 2\sum_{j=1}^n\sum_{l=0}^{s-p}\binom{s-p}{l}(\partial_t^l\gamma^{0j})(\partial_j\partial_t^{s+1-p-l}u)(t,\cdot)$$
$$+ \partial_t^{s+2-p}u(t,\cdot) - \sum_{i,j=1}^n\sum_{l=1}^{s-p}\binom{s-p}{l}(\partial_t^l\gamma^{ij})(\partial_i\partial_j\partial_t^{s-p-l}u)(t,\cdot) + \partial_t^{s-p}G(t,\cdot). \quad (9.17)$$

The $H^{p-1}$ norm of the right hand side involves terms of the form
$$\sum_{|\alpha_1|+|\alpha_2|\leq p-1}(\partial_x^{\alpha_1}\partial_t^l\gamma^{ij})(\partial_x^{\alpha_2}\partial_i\partial_j\partial_t^{s-p-l}u)(t,\cdot)\,,$$

where either $i=0$ or $l\geq 1$. Consider such terms for which $|\alpha_1|+l\leq k$. Since $\gamma^{ij}\in C^k$, these can be controlled by

$$\sum_{m=0}^{s+1-p}\|\partial_t^m u(t,\cdot)\|_{H^p} + \sum_{m=0}^{s-1-p}\|\partial_t^m u(t,\cdot)\|_{H^{p+1}}\,.$$

On the other hand, if $|\alpha_1|+l>k$, then $|\alpha_2|\leq p+l-k-2\leq s-k-2$. Since $p\geq k$, it follows that $p+k>s+\frac{n}{2}$, hence

$$\|\partial_x^{\alpha_2}\partial_i\partial_j\partial_t^{s-p-l}u(t,\cdot)\|_\infty \leq \sum_{m=0}^{s+1-p}\|\partial_t^m u(t,\cdot)\|_{H^p}\,.$$

Consequently, the $H^{p-1}$ norm of the right hand side of (9.17) is bounded by

$$C\int_0^t \|\partial_t^{s-p}u(r,\cdot)\|_{H^{p+1}}\,dr + C\sup_{0\leq r\leq t}\sum_{m=0}^p \|\partial_t^{s+1-m}u(r,\cdot)\|_{H^m}$$
$$+ C\left(\|f\|_{H^{s+1}} + \|g\|_{H^s} + \sup_{|\alpha|\leq s-1}\|\partial^\alpha G(t,\cdot)\|_2\right),$$

where we are using the fact that we may bound

$$\sum_{m=0}^{s+1} \|\partial_t^{s+1-m} u(0,\cdot)\|_{H^m} \leq C\left(\|f\|_{H^{s+1}} + \|g\|_{H^s}\right),$$

since the compatibility functions $\psi_j$ are linear expressions in $f$ and $g$. We may thus use Lemma 9.9 to conclude that (9.16) holds with $k$ replaced by $p+1$. By the continuity of the right hand side of (9.17), it is easy to see that $\partial_t^{s-p} u \in C([0,T]; H^{p+1})$, completing the proof of the Theorem 9.10. □

**Proof of Theorem 9.4.** We produce a solution to (9.1) by iteration. Thus, define the sequence of functions $u_l$ by letting

$$\partial_t^2 u_{l+1} = \Delta u_{l+1} + \sum_{i,j=0}^{n} \gamma^{ij}(t,x,J_1 u_l, \partial_t u_l)\,\partial_i\partial_j u_{l+1} + G(t,x,J_1 u_l, \partial_t u_l),$$

$$u_{l+1}(0,x) = f(x), \qquad \partial_t u_{l+1}(0,x) = g(x), \qquad u_{l+1}(t,x) = 0 \text{ if } x \in \partial\Sigma.$$

We take $u_0$ be the solution of the nonlinear problem (9.1) without Dirichlet conditions, where the data $f, g$ are extended across $\partial\Sigma$. The existence of $u_0$ on some interval $[0, T']$, where $T'$ depends only on bounds for the norms of $f$ and $g$, follows by [6] Theorem 6.4.11. Since $\partial_t^k u_0(0,x) = \psi_k(x)$, it follows from Lemma 9.3 that the compatibility functions and conditions are the same at each step of the iteration as for the nonlinear problem (9.1), and hence the existence of the sequence $u_l$ follows by Theorem 9.10.

We now show that there exists $M < \infty$ and $T > 0$ such that such that

$$M_{s+1}(u_l) = \sup_{0\leq t\leq T}\sum_{|\alpha|\leq s+1} \|\partial^\alpha u_l(t,\cdot)\|_2 \leq M \tag{9.18}$$

for all values of $l$. We let $M = 8C\left(M_{s+1}(u_0) + 1\right)$, where

$$M_{s+1}(u_0) = \sup_{0\leq t\leq T'}\sum_{|\alpha|\leq s+1} \|\partial^\alpha u_0(t,\cdot)\|_2,$$

and establish (9.18) by induction. Thus, assume that (9.18) holds for $l$, where $T$ is small enough so that

$$C[M]\,T^{\frac{1}{3}} \leq \frac{1}{2}, \qquad CMT \leq \frac{1}{2}, \qquad 2\,C\,e^{1/2}\,C[M]\,T \leq \frac{1}{2},$$

where the various constants are as in (9.9). Then, by (9.9), we have

$$M_{s+1}(u_{l+1},t) \leq 2\,C\,e^{1/2}\left(M_{s+1}(u_0) + 1 + C[M]\int_0^t M_{s+1}(u_{l+1},r)\,dr\right),$$

where we assumed that $C \geq 1$, and we used the fact that

$$M_{s+1}(u_{l+1},0) = M_{s+1}(u_0,0) \leq M_{s+1}(u_0).$$

By Gronwall's lemma, we thus have

$$M_{s+1}(u_{l+1}) \leq 2\,C\,e^{1/2}\bigl(M_{s+1}(u_0)+1\bigr)e^{2Ce^{1/2}C[M]T} \leq M.$$

We conclude by showing that the sequence $u_l$ is Cauchy in the energy norm. By weak compactness it then follows that $u_l$ converges to a solution $u$ of (9.1) that satisfies the conditions of Theorem 9.4, completing the proof of the first part of the theorem.

By subtracting successive equations, we obtain

$$\partial_t^2(u_{l+1} - u_l) - \Delta(u_{l+1} - u_l) - \sum_{i,j=0}^n \gamma^{ij}(t, x, J_1 u_l, \partial_t u_l) \, \partial_i \partial_j (u_{l+1} - u_l)$$

$$= \sum_{i,j=0}^n \left(\gamma^{ij}(t, x, J_1 u_l, \partial_t u_l) - \gamma^{ij}(t, x, J_1 u_{l-1}, \partial_t u_{l-1})\right) \partial_i \partial_j u_l$$

$$+ G(t, x, J_1 u_l, \partial_t u_l) - G(t, x, J_1 u_{l-1}, \partial_t u_{l-1}).$$

Since we have uniform bounds for the $C^2$ norm of $u_l$ for all $l$, and the functions $\gamma^{ij}$ and $G$ are smooth in their arguments, we can bound the $L^2$ norm of the right hand side by

$$C \sum_{|\alpha| \leq 1} \left\| \partial^\alpha (u_l - u_{l-1}) \right\|_2.$$

Since the Cauchy data of $u_{l+1} - u_l$ vanishes, we can apply the energy inequality to obtain

$$\sum_{|\alpha| \leq 1} \left\| \partial^\alpha (u_{l+1} - u_l)(t, \cdot) \right\|_2 \leq C \int_0^t \sum_{|\alpha| \leq 1} \left\| \partial^\alpha (u_l - u_{l-1})(r, \cdot) \right\|_2 dr,$$

and hence

$$\sup_{0 \leq t \leq T} \sum_{|\alpha| \leq 1} \left\| \partial^\alpha (u_{l+1} - u_l)(t, \cdot) \right\|_2 \leq M \frac{(CT)^l}{l!}.$$

It remains to show that if $\left( \|f\|_{H^{s+1}} + \|g\|_{H^s} \right) \leq 1$, and $G(t, x, 0, 0) = 0$, then we may take $C$ and $T$ independent of $f$ and $g$ so that (9.6) holds. To see this, let $u$ be the solution to (9.1), and observe that we have uniform bounds on $M_{s+1}(u, t)$ for $0 \leq t \leq T$, independent of $f$ and $g$. We thus have

$$\sum_{|\alpha| \leq s} \|\partial^\alpha G(t, x, J_1 u, \partial_t u)\|_2 \leq C \, M_{s+1}(u, t),$$

for some constant $C$. Consequently, we may replace (9.9) by the following inequality,

$$M_{s+1}(u, t) \leq C \, e^{Ct} \left( M_{s+1}(u, 0) + \int_0^t M_{s+1}(u, r) \, dr \right) + C \, T^{\frac{1}{3}} M_{s+1}(u, t).$$

We take $T$ small so that $C \, T^{\frac{1}{3}} \leq \frac{1}{2}$, and apply the Gronwall Lemma to obtain

$$\sup_{0 \leq t \leq T} M_{s+1}(u, t) \leq 2 \, C \, e^{CT} e^{2Ce^{CT}} M_{s+1}(u, 0).$$

Since $\partial_t^k u(0, x) = \psi_k(J_k f, J_{k-1} g)$, we have the bounds

$$M_{s+1}(u, 0) \leq \sum_{k=0}^{s+1} \|\psi_k(J_k f, J_{k-1} g)\|_{H^{s+1-k}},$$

where the functions $\psi_k(J_k f, J_{k-1} g)$ are defined recursively by (9.4). Since $G(t, x, 0, 0) = 0$, it follows that $F(t, x, 0, 0) = 0$, and the proof of Lemma 9.1 shows inductively that

$$\sum_{k=0}^{s+1} \|\psi_k(J_k f, J_{k-1} g)\|_{H^{s+1-k}} \leq C \left( \|f\|_{H^{s+1}} + \|g\|_{H^s} \right),$$

provided that $\left( \|f\|_{H^{s+1}} + \|g\|_{H^s} \right) \leq 1$. $\square$

**Proof of Theorem 9.5.** By Theorem 9.4, we have local existence of solutions of regularity $m+1$ given that the compatibility conditions of order $m$ are satisfied. Since the compatibility conditions propagate, it suffices to show that, given a solution to (9.1) such that the quantity

$$\sum_{|\alpha|\leq m+1} \|\partial^\alpha u(t,\cdot)\|_2$$

is locally bounded for $0 \leq t < T'$, and such that

$$\sup_{0\leq t\leq T'} \sum_{|\alpha|\leq s+1} \|\partial^\alpha u(t,\cdot)\|_2 < \infty, \tag{9.19}$$

then it follows that

$$M_{m+1}(u) = \sup_{0\leq t< T'} \sum_{|\alpha|\leq m+1} \|\partial^\alpha u(t,\cdot)\|_2 < \infty.$$

We apply Lemma 9.6 in the case $u = v$. By (9.19), we are given uniform bounds on the $C^2$ norm $M$ of $u$. We conclude that, for $k \geq s+1$, and $0 \leq t < T'$,

$$M_{k+1}(u,t) \leq C\, e^{CMT'} \left( M_{k+1}(u,0) + 2\, C[M] \int_0^t M_{k+1}(u,r)\, dr \right) + C\bigl[T', M_k(u)\bigr],$$

and consequently

$$M_{k+1}(u) \leq \Bigl( C\, e^{CMT'} M_{k+1}(u,0) + C\bigl[T', M_k(u)\bigr] \Bigr) e^{C\, e^{2CMT'} C[M] T'}.$$

The proof now follows by induction, given that $M_k(u)$ is bounded for $k = s+1$. □

Department of Mathematics, California Institute of Technology, Pasadena, CA 91125

Department of Mathematics, University of Washington, Seattle, WA 98195

Department of Mathematics, The Johns Hopkins University, Baltimore, MD 21218